\documentclass[a4paper,11pt,draft]{article}
\usepackage{amssymb,amsmath,amsfonts,amsthm}
\newcommand\dslash{d\llap {\raisebox{.9ex}{$\scriptstyle-\!$}}}
\newtheorem{theorem}{Theorem}[section]
\newtheorem{corollary}{Corollary}[section]
\newtheorem{lemma}{Lemma}[section]
\newtheorem{proposition}{Proposition}[section]
\newtheorem{definition}{Definition}[section]

\newtheorem{remark}{Remark}
\newtheorem{example}{Example}

\def\ni{\noindent}
\def\amp{{\mathcal{A}}^{m}_{l}}
\def\ampN{{\mathcal{A}}^{m}_{l,N}}
\def\ampM{{\mathcal{A}}^{m}_{l,M}}
\def\S{{\mathcal{S}}(\mathbb{R}^n)}

\def\Sp{{\mathcal{S}^{'}}(\mathbb{R}^n)}
\def\ampg{\overline{\mathcal{S}}^{m}_{\Lambda,\rho}}
\def\gts{\mathcal{G}_{\tau,\mathcal{S}}(\mathbb{R}^n)}
\def\gss{\mathcal{G}^\infty_{\mathcal{S}}(\mathbb{R}^n)}
\def\ampc{\overline{\mathcal{S}}^{m}_{\Lambda,\rho,0}}
\def\amprN{\overline{\mathcal{S}}^{m}_{\Lambda,\rho,N}}
\def\symrN{{\mathcal{S}}^{m}_{\Lambda,\rho,N}}

\begin{document}
\catcode`\@=11


  \renewcommand{\theequation}{\thesection.\arabic{equation}}
  \renewcommand{\section}%
  {\setcounter{equation}{0}\@startsection {section}{1}{\z@}{-3.5ex plus -1ex
   minus -.2ex}{2.3ex plus .2ex}{\Large\bf}}
\title{\bf Pseudo-differential operators in algebras of generalized functions and global hypoellipticity}
\author{Claudia Garetto\\
[0.2cm]
Dipartimento di Matematica, Universit\`a di Torino,\\
via Carlo Alberto 10, 10123 Torino, Italia\\[0.2cm]
\texttt{garettoc@dm.unito.it}\\ }
\date{ }
\maketitle
\begin{abstract}
\ni The aim of this work is to develop a global calculus for pseudo-dif\-fe\-ren\-tial operators acting on suitable algebras of generalized functions. In particular, a condition of global hypoellipticity of the symbols gives a result of regularity for the corresponding pseudo-differential equations. This calculus and this frame are proposed as tools for the study in Colombeau algebras of partial differential equations globally defined on $\mathbb{R}^n$.
\end{abstract}
\bf{Key words}\rm: Colombeau algebras, pseudo-differential operators, global hypoellipticity\\[0.2cm]
AMS mathematics subject classification (1991): 46F30, 47G30, 35H05 
\section{Introduction}
Colombeau's theory of generalized functions has been developed in connection with nonlinear problems, see \cite{biag, ober, ros}, but it is also important for linear problems \cite{ned}. In particular we are interested in the development of the pseudo-differential calculus in the frame of generalized functions and tempered ge\-ne\-ra\-li\-zed functions \cite{gram, nedI, ned, piliT, pili}.\\
In Section 2 of this paper we recall the basic notions of Colombeau's theory, Colombeau-Fourier transformation and weak equality \cite{biag, colo, coloII, hoer, ned, ober, pili}, considering in place of the more usual $\mathcal{G}_{\tau}(\mathbb{R}^n)$ \cite{colo}, a new algebra of generalized functions $\gts$ containing $\mathcal{S}^{'}(\mathbb{R}^n)$.
More precisely, fol\-lo\-wing arguments similar to Radyno \cite{anto, rad}, in our construction we substitute the ideal ${\mathcal{N}_\tau(\mathbb{R}^n)}$ of $\mathcal{G}_\tau(\mathbb{R}^n)$ by the smaller ${\mathcal{N}}_{\mathcal{S}}(\mathbb{R}^n)$, modelled on $\mathcal{S}(\mathbb{R}^n)$. 
The choice of $\gts$ is motivated by the existence of a subalgebra of $\mathcal{S}$-regular generalized functions, denoted by $\gss$, for which we prove the equality $\gss\cap\mathcal{S}^{'}(\mathbb{R}^n)=\mathcal{S}(\mathbb{R}^n)$.
This result of regularity is inspired by the local version given by Oberguggenberger \cite{ober}, and a first global investigation, involving $\mathcal{S}^{'}(\mathbb{R}^n)$ and $\mathcal{O}_M(\mathbb{R}^n)$, proposed by H\"{o}rmann \cite{hoer}.\\
\ni In Sections 3 and 4 we collect some preliminary arguments necessary for the definition of pseudo-differential operators acting on $\gts$ and their global calculus. In detail we consider symbols and amplitudes depending on a parameter $\epsilon\in(0,1]$, whose global estimates involve a weight function $\Lambda(x,\xi)$ \cite{bog, bogII}, and powers of $\epsilon$ with negative exponent \cite{nedI, ned, piliT, pili}.\\
\ni The definition of pseudo-differential operators on $\gts$ with cor\-re\-spon\-ding mapping pro\-per\-ties, a characterization of operators with $\mathcal{S}$-regular kernel and the comparison with the classical pseudo-differential operators on $\mathcal{S}^{'}(\mathbb{R}^n)$ and $\mathcal{S}(\mathbb{R}^n)$ are the topics of Sections 5 and 6. Star\-ting points for our reasoning are \cite{nedI, ned, piliT, pili} and for the classical theory \cite{bog, bogII, shu}.\\
\ni The relationship with the work \cite{nedI, ned, piliT, pili} of Nedeljkov, Pilipovi\'{c} and Scarpal\'{e}zos, where a definition of generalized pseudo-differential operators acting on $\mathcal{G}_\tau(\mathbb{R}^n)$ is given, is evident in the estimates on $\epsilon$, while the elements of novelty are the choice of the weight function $\Lambda$ in place of the standard $\langle\xi\rangle$ and the generalization, to the Colombeau setting based on $\gts$ and $\gss$, of the classical global calculus.\\
\ni In particular, in order to deal with the composition formula, we present in Section 7 a ge\-ne\-ra\-li\-za\-tion of the well-known Weyl symbols. By means of this technical tool, we prove that the composition of two pseudo-differential o\-pe\-ra\-tors acting on $\gts$ coincides with the action of a pseudo-differential operator in the weak or g.t.d. (generalized tempered distributions) sense.\\
\ni The introduction of suitable sets of hypoelliptic and elliptic symbols in Section 8, allows us to obtain, through the construction of a parametrix, an interesting result of regularity, modelled on the  classical statement that if  $A$ is a pseudo-differential operator with hypoelliptic symbol, $u\in\mathcal{S}^{'}(\mathbb{R}^n)$, $f\in\mathcal{S}(\mathbb{R}^n)$, then $Au=f$ implies $u\in\mathcal{S}(\mathbb{R}^n)$ \cite{bog, bogII}.\\
\ni Finally, let us emphasize that applications of Colombeau algebras to the study of PDE's in  \cite{biag, colo, coloII, gram, ober} concern mainly local problems; our calculus and our definitions of $\gts$ and $\gss$ are proposed as tools for the study of equations globally defined on $\mathbb{R}^n$. At the end of the paper we give some examples, considering above all partial differential operators of the form $P=\sum_{(\alpha,\beta)\in\mathcal{A}}c_{\alpha,\beta}x^\alpha D^\beta$, where $\mathcal{A}$ is a finite subset of multi-indices in $\mathbb{N}^{2n}$ and the coefficients $c_{\alpha,\beta}$ are Colombeau generalized numbers.
\section{Basic notions}
In this section we recall the definitions and results needed from the theory of Colombeau ge\-ne\-ra\-li\-zed functions. Since we do not motivate the constructions and we do not provide proofs, we refer for details to \cite{biag, colo, coloII, gro, ned, ober, piliT, pili, ros}. We begin with considering the simplified Colombeau algebra on an open subset $\Omega$ of $\mathbb{R}^n$. This is obtained as a factor algebra in the differential algebra $\mathcal{E}[\Omega]$ of all the sequences $(u_\epsilon)_{\epsilon\in(0,1]}$ of smooth functions $u_\epsilon\in\mathcal{C}^\infty(\Omega)$. In the sequel we use in place of $(u_\epsilon)_{\epsilon\in(0,1]}$ the simpler notation $(u_\epsilon)_{\epsilon}$.
\begin{definition}
We call moderate the elements of $\mathcal{E}[\Omega]$ such that for all $K\subset\subset\Omega$, for all $\alpha\in\mathbb{N}^n$, there exists $N\in\mathbb{N}$, such that
\[
\sup_{\epsilon\in(0,1]}\epsilon^{N}\Vert \partial^\alpha u_\epsilon\Vert_{L^{\infty}(K)}<\infty .
\]
\end{definition}
\noindent
The set of these elements is a differential algebra denoted by ${\mathcal{E}}_M(\Omega)$.
\begin{definition}
We call negligible the elements of $\mathcal{E}[\Omega]$ such that for all $K\subset\subset\Omega$, for all $\alpha\in\mathbb{N}^n$ and for all $q\in\mathbb{N}$ 
\[
\sup_{\epsilon\in(0,1]}\epsilon^{-q}\Vert \partial^\alpha u_\epsilon\Vert_{L^{\infty}(K)}<\infty .
\]
\end{definition}
\noindent
$\mathcal{N}(\Omega)$, the set of these elements, is an ideal of ${\mathcal{E}}_M(\Omega)$, closed with respect to derivatives.\\ The \it{Colombeau algebra of generalized functions}\rm\ is $\mathcal{G}(\Omega)={\mathcal{E}_M(\Omega)}/{\mathcal{N}(\Omega)}$.\\
Let $\varphi\in\mathcal{S}(\mathbb{R}^n)$ with
\begin{equation}
\int_{\mathbb{R}^n}\hskip-5pt\varphi(x)dx=1\quad\ \text{and}\ \quad \int_{\mathbb{R}^n}\hskip-5pt x^\alpha\varphi(x)dx=0
\end{equation}
for all $\alpha\in\mathbb{N}^n\ ,\ \alpha\neq 0$. As usual we put
\begin{equation}
\varphi_\epsilon(x)=\epsilon^{-n}\varphi(\frac{x}{\epsilon}) .
\end{equation}
Using the convolution product we can easily define the following embedding of $\mathcal{E}^{'}(\Omega)$ into $\mathcal{G}(\Omega)$:
\begin{equation}
\imath_o:\mathcal{E}^{'}(\Omega)\rightarrow\mathcal{G}(\Omega):w\rightarrow((w\ast\varphi_\epsilon)\vert_{\Omega})_\epsilon +\mathcal{N}(\Omega) .
\end{equation}
This map can be extended to an embedding $\imath$ of $\mathcal{D}^{'}(\Omega)$, employing the sheaf properties of $\mathcal{G}(\Omega)$ and a suitable partition of unity, and it renders $\mathcal{C}^\infty(\Omega)$ a subalgebra. Moreover in $\mathcal{G}(\Omega)$ the derivatives extend the usual ones in the sense of distributions.\\
In order to talk of tempered generalized functions, we introduce the following subalgebras of $\mathcal{E}[\mathbb{R}^n]$.
\begin{definition}
$\mathcal {E}_{\tau}(\mathbb{R}^n)$ is the set of all elements $(u_\epsilon)_\epsilon$ belonging to $\mathcal{E}[\mathbb{R}^n]$ with the following property: for all $\alpha\in\mathbb{N}^n$ there exists $N\in\mathbb{N}$ such that
\[
\sup_{\epsilon\in(0,1]}\epsilon^{N}\Vert\langle x\rangle^{-N}\partial^\alpha u_\epsilon\Vert_{L^\infty(\mathbb{R}^n)}<\infty .
\]
\end{definition}
\begin{definition}
$\mathcal{N}_{\tau}(\mathbb{R}^n)$ is the set of all elements $(u_\epsilon)_\epsilon$ belonging to $\mathcal{E}[\mathbb{R}^n]$ with the following property: for all $\alpha\in\mathbb{N}^n$ there exists $N\in\mathbb{N}$ such that for all $q\in\mathbb{N}$
\[
\sup_{\epsilon\in(0,1]}\epsilon^{-q}\Vert\langle x\rangle^{-N}\partial^\alpha u_\epsilon\Vert_{L^\infty(\mathbb{R}^n)}<\infty .
\]
\end{definition}
\noindent
The \it{Colombeau algebra of tempered generalized functions}\rm\ is by definition the factor $\mathcal{G}_\tau(\mathbb{R}^n)=\mathcal {E}_{\tau}(\mathbb{R}^n)/\mathcal{N}_{\tau}(\mathbb{R}^n)$. It is a differential algebra containing the space of distributions $\mathcal{S}^{'}(\mathbb{R}^n)$, where the derivatives extend the usual ones on $\mathcal{S}^{'}(\mathbb{R}^n)$ and 
\[
\mathcal{O}_{C}(\mathbb{R}^n)=\{f\in\mathcal{C}^\infty(\mathbb{R}^n):\, \exists N\in\mathbb{N}:  \forall\alpha\in\mathbb{N}^n ,\ \ \Vert\langle x\rangle^{-N}\partial^\alpha f\Vert_{L^\infty(\mathbb{R}^n)}<\infty\}
\]
is a subalgebra. In particular the embedding is done by
\[
\imath:\mathcal{S}^{'}(\mathbb{R}^n)\rightarrow\mathcal{G}_\tau(\mathbb{R}^n):w\rightarrow (w\ast\varphi_\epsilon)_\epsilon+\mathcal{N}_\tau(\mathbb{R}^n) .
\]
If $f\in\mathcal{O}_C(\mathbb{R}^n)$, we can define the constant embedding $\sigma(f)=(f)_\epsilon+\mathcal{N}_\tau(\mathbb{R}^n)$, and in this way the important equality $\imath(f)=\sigma(f)$ holds in $\mathcal{G}_\tau(\mathbb{R}^n)$.\\
\ni The constants of $\mathcal{G}(\mathbb{R}^n)$ or respectively $\mathcal{G}_\tau(\mathbb{R}^n)$, constitute the \it{algebra of Co\-lom\-beau ge\-ne\-ra\-li\-zed complex numbers $\overline{\mathbb{C}}$}.\rm\ It contains $\mathbb{C}$ and is defined as the factor $\overline{\mathbb{C}}=\mathcal{E}_{o,M}/\mathcal{N}_o$, where $\mathcal{E}_o=\mathbb{C}^{(0,1]}$ and 
\begin{equation}
\begin{split}
\mathcal{E}_{o,M}&=\{ (h_\epsilon)_\epsilon\in\mathcal{E}_o :\, \exists N\in\mathbb{N} :\ \sup_{\epsilon\in(0,1]}\epsilon^{N}|h_\epsilon|<\infty\} ,\\
\mathcal{N}_o&=\{ (h_\epsilon)_\epsilon\in\mathcal{E}_o :\, \forall q\in\mathbb{N} ,\  \sup_{\epsilon\in(0,1]}\epsilon^{-q}|h_\epsilon|<\infty\} .
\end{split}
\end{equation}
The generalized complex numbers allow us to integrate on $\mathbb{R}^n$ an arbitrary element of $\mathcal{G}_\tau(\mathbb{R}^n)$. In the sequel we denote a representative of $u\in\mathcal{G}_\tau(\mathbb{R}^n)$ with $(u_\epsilon)_\epsilon\in\mathcal{E}_\tau(\mathbb{R}^n)$.
\begin{proposition}
Let $u\in\mathcal{G}_\tau(\mathbb{R}^n)$. Then
\[
\int_{\mathbb{R}^n}\hskip-4pt u(x)dx:=\biggl(\int_{\mathbb{R}^n}\hskip-4pt u_\epsilon(x)\widehat{\varphi_\epsilon}(x)dx\biggr)_\epsilon +\mathcal{N}_o ,
\]
where $\varphi_\epsilon$ is as in (2.2), is a well-defined element of $\overline{\mathbb{C}}$, called the integral of $u$ over $\mathbb{R}^n$.
\end{proposition}
\noindent
The article ``the'' in the expression ``the integral of $u$ over $\mathbb{R}^n$'' has been used cum grano salis, since different variants of integrals are used and studied earlier, e.g., in \cite{hoer, ned}.
We collect now the main properties of the integral in $\mathcal{G}_\tau(\mathbb{R}^n)$:
\begin{itemize}
\item[i)]let $u\in\mathcal{G}_\tau(\mathbb{R}^n)$ and $f\in\mathcal{S}(\mathbb{R}^n)$. Then for all $\alpha\in\mathbb{N}^n$
\[
\int_{\mathbb{R}^n}\hskip-4pt\partial^\alpha u(x)\imath(f)(x)dx=(-1)^{|\alpha|}\int_{\mathbb{R}^n}\hskip-4pt u(x)\imath(\partial^\alpha f)(x)dx ;
\]
\item[ii)]let $w\in\mathcal{S}^{'}(\mathbb{R}^n)$ and $f\in\S$. Then $\displaystyle\int_{\mathbb{R}^n}\hskip-4pt\imath(w)(x)\imath(f)(x)dx=(\langle w,f\rangle)_\epsilon+\mathcal{N}_o ;$
\item[iii)]let $f\in\mathcal{S}(\mathbb{R}^n)$. Then $\displaystyle\int_{\mathbb{R}^n}\hskip-4pt\imath(f)(x)dx=\bigg(\int_{\mathbb{R}^n} f(x)dx\bigg)_\epsilon+\mathcal{N}_o .$
\end{itemize}
Using this definition of integral we introduce in $\mathcal{G}_\tau(\mathbb{R}^n)$ a weak equality.
\begin{definition}
$u$ and $v$ in $\mathcal{G}_\tau(\mathbb{R}^n)$ are equal in the sense of generalized tempered distribution (or weak sense) iff for all $f\in\mathcal{S}(\mathbb{R}^n)$
\[
\int_{\mathbb{R}^n}\hskip-4pt(u-v)(x)\imath(f)(x)dx=0 .
\]
\end{definition}
\noindent
We write $u=_{g.t.d.}v$. Then $=_{g.t.d.}$ is an equivalence relation compatible with the linear struc\-tu\-re and the derivatives in $\mathcal{G}_\tau(\mathbb{R}^n)$. From ii) it follows that  $\mathcal{S}^{'}(\mathbb{R}^n)$ is a subspace of the factor $\mathcal{G}_\tau(\mathbb{R}^n)/=_{g.t.d.}$.\\
There exists in $\mathcal{G}_\tau(\mathbb{R}^n)$ a natural definition of \it{Colombeau-Fourier transform and anti-transform}\rm. 
\begin{definition}
Let $u\in\mathcal{G}_\tau(\mathbb{R}^n)$. The Colombeau-Fourier transform of u is given by the representative 
\[
\mathcal{F}_\varphi u_\epsilon(\xi)=\int_{\mathbb{R}^{n}}\hskip-7pt e^{-iy\xi} u_\epsilon(y)\widehat{\varphi_\epsilon}(y)dy .
\]
The Colombeau-Fourier anti-transform of $u$ is given by the representative
\[
\mathcal{F}^{\ast}_\varphi u_\epsilon(y)=\int_{\mathbb{R}^{n}}\hskip-7pt e^{iy\xi} u_\epsilon(\xi)\widehat{\varphi_\epsilon}(\xi)\dslash\xi ,
\]
where $\dslash\xi=(2\pi)^{-n}d\xi$.
\end{definition}
\ni One can easily prove that the previous definition makes sense: in this way ${\mathcal{F}}_\varphi$, respectively  $\mathcal{F}_\varphi^{\ast}$, defines a linear map from ${\mathcal{G}}_\tau(\mathbb{R}^n)$ into $\mathcal{G}_\tau(\mathbb{R}^n)$. Moreover the following properties hold:
\begin{itemize}
\item[i)]${\mathcal{F}}_\varphi$ and $\mathcal{F}_\varphi^{\ast}$ extend the classical transformations on $\mathcal{S}(\mathbb{R}^n)$; in other words for all $f\in\mathcal{S}(\mathbb{R}^n)$ $\mathcal{F}_\varphi(\imath(f))=\imath(\hat{f})$ and $\mathcal{F}_\varphi^{\ast}(\imath(f))=\imath(\check{f})$;
\item[ii)] for all $u\in\mathcal{G}_\tau(\mathbb{R}^n)$ and for all $f\in\mathcal{S}(\mathbb{R}^n)$ 
\[
\begin{split}
\int_{\mathbb{R}^n}\hskip-8pt\mathcal{F}_\varphi u(x)\imath(f)(x)dx&=\int_{\mathbb{R}^n}\hskip-8pt u(x)\imath(\hat{f})(x)dx ,\\
\int_{\mathbb{R}^n}\hskip-8pt \mathcal{F}^{\ast}_\varphi u(x)\imath(f)(x)dx&=\int_{\mathbb{R}^n}\hskip-8pt u(x)\imath(\check{f})(x)dx ;
\end{split}
\]
\item[iii)] for all $w\in\mathcal{S}^{'}(\mathbb{R}^n)$,  $\mathcal{F}_\varphi(\imath(w))=_{g.t.d.}\imath(\hat{w})$ and $\mathcal{F}^{\ast}_\varphi(\imath(w))=_{g.t.d.}\imath(\check{w})$;
\item[iv)] in general $u\neq\mathcal{F}^{\ast}_\varphi\mathcal{F}_\varphi u\neq \mathcal{F}_\varphi\mathcal{F}^{\ast}_\varphi u$ but $u=_{g.t.d.}\mathcal{F}^{\ast}_\varphi\mathcal{F}_\varphi u=_{g.t.d.}\mathcal{F}_\varphi\mathcal{F}^{\ast}_\varphi u ;$
\item[v)] for all $u\in\mathcal{G}_\tau(\mathbb{R}^n)$, $\alpha\in\mathbb{N}^n$, $\mathcal{F}_\varphi(\imath(y^\alpha)u)=i^{|\alpha|}\partial^\alpha\mathcal{F}_\varphi u$ and $\mathcal{F}^{\ast}_\varphi\hskip-2pt(\imath(y^\alpha)u)$\\
$=(-i)^{|\alpha|}\partial^\alpha\mathcal{F}^{\ast}_\varphi u ;$
\item[vi)] for all $u\in\mathcal{G}_\tau(\mathbb{R}^n)$ and $\alpha\in\mathbb{N}^n$, $(-i)^{|\alpha|}\mathcal{F}_\varphi(\partial^\alpha u)=_{g.t.d.}\imath(y^\alpha)\mathcal{F}_\varphi u$ and\\
$i^{|\alpha|}\mathcal{F}^{\ast}_\varphi(\partial^\alpha u)=_{g.t.d.}\imath(y^\alpha)\mathcal{F}^{\ast}_\varphi u$, while the equalities in $\mathcal{G}_\tau(\mathbb{R}^n)$ are not true.
\end{itemize}
As a consequence in order to obtain the usual properties of Fourier transform and an\-ti\-tran\-sform, we consider the definition of $\mathcal{F}_\varphi$ and $\mathcal{F}^{\ast}_\varphi$ on the factor $\mathcal{G}_\tau(\mathbb{R}^n)/=_{g.t.d.}$.\\
We conclude this section by reporting some results concerning regularity theory. The star\-ting point for regularity theory and microlocal analysis in Colombeau algebras of generalized functions was the introduction of the subalgebra $\mathcal{G}^\infty(\Omega)$ of $\mathcal{G}(\Omega)$ by Oberguggenberger in \cite{ober}.
\begin{definition}
$\mathcal{G}^\infty(\Omega)$ is the set of all $u\in\mathcal{G}(\Omega)$ having a representative $(u_\epsilon)_\epsilon\in\mathcal{E}_M(\Omega)$ with the following property: for all $K\subset\subset\Omega$ there exists  $N\in\mathbb{N}$ such that for all $\alpha\in\mathbb{N}^n$
\begin{equation}
\sup_{\epsilon\in(0,1]}\epsilon^{N}\Vert\partial^\alpha u_\epsilon\Vert_{L^\infty(K)}<\infty .
\end{equation}
\end{definition}
\noindent In \cite{ober}, th.25.2, the identity $\mathcal{G}^\infty(\Omega)\cap\mathcal{D}^{'}(\Omega)=\mathcal{C}^\infty(\Omega)$ is proved. We introduce now a suitable algebra of $\mathcal{S}$-regular generalized functions on $\mathbb{R}^n$. At first we consider another differential algebra containing $\mathcal{S}^{'}(\mathbb{R}^n)$ and $\mathcal{S}(\mathbb{R}^n)$, where  the ideal satisfies an estimate of rapidly decreasing type.
\begin{definition}
We denote by $\mathcal{G}_{\tau,\mathcal{S}}(\hskip-2pt\mathbb{R}^n\hskip-2pt )$ the factor $\mathcal {E}_{\tau}(\mathbb{R}^n)/\mathcal{N}_{\mathcal{S}}(\mathbb{R}^n)$, where $\mathcal{N}_{\mathcal{S}}(\mathbb{R}^n)$ is the set of all $(u_\epsilon)_\epsilon\in\mathcal{E}[\mathbb{R}^n]$ fulfilling the following condition:
\begin{equation}
\begin{array}{cc}
\forall\alpha,\beta\in\mathbb{N}^n,\ \forall q\in\mathbb{N},\\[0.2cm]
\displaystyle\sup_{\epsilon\in(0,1]}\epsilon^{-q}\Vert x^\alpha\partial^\beta u_\epsilon\Vert_{L^\infty(\mathbb{R}^n)}<\infty .
\end{array}
\end{equation}
\end{definition}
\noindent Since $\mathcal{N}_{\mathcal{S}}(\mathbb{R}^n)\subset\mathcal{N}_\tau(\mathbb{R}^n)$, $\mathcal{S}^{'}(\mathbb{R}^n)$ is a subspace of $\mathcal{G}_{\tau,\mathcal{S}}(\mathbb{R}^n)$, and for all $f\in\mathcal{S}(\mathbb{R}^n)$, $(f-f\ast\varphi_\epsilon)_\epsilon\in\mathcal{N}_{\mathcal{S}}(\mathbb{R}^n)$ implies the embedding as subalgebra of ${\mathcal{S}}(\mathbb{R}^n)$ into $\mathcal{G}_{\tau,\mathcal{S}}(\mathbb{R}^n)$.
For simplicity we continue to denote the class of $w\in\mathcal{S}^{'}(\mathbb{R}^n)$ in $\mathcal{G}_{\tau,\mathcal{S}}(\mathbb{R}^n)$ with $\imath(w)$. Obviously Proposition 2.1 and the corresponding properties of integral hold with $\gts$ in place of $\mathcal{G}_\tau(\mathbb{R}^n)$.
\begin{proposition}
$\mathcal{F}_\varphi$ and $\mathcal{F}_\varphi^{\ast}$ map $\gts$ into $\gts$.
\end{proposition}
\begin{proof}
It suffices to prove that $(u_\epsilon)_\epsilon\in\mathcal{N}_{\mathcal{S}}(\mathbb{R}^n)$ implies $(\mathcal{F}_\varphi u_\epsilon)_\epsilon\in\mathcal{N}_{\mathcal{S}}(\mathbb{R}^n)$ and $(u_\epsilon)_\epsilon\in\mathcal{N}_{\mathcal{S}}(\mathbb{R}^n)$ implies $(\mathcal{F}^{\ast}_\varphi u_\epsilon)_\epsilon\in\mathcal{N}_{\mathcal{S}}(\mathbb{R}^n)$.
\end{proof}
\ni All the properties of the Colombeau-Fourier transform mentioned above hold in $\gts$ as well. 
\begin{definition}
An element $u\in\mathcal{G}_{\tau,{\mathcal{S}}}(\mathbb{R}^n)$ is called $\mathcal{S}$-regular (or $u\in\mathcal{G}^\infty_{\mathcal{S}}(\mathbb{R}^n)$) if it has a representative $(u_\epsilon)_\epsilon$ such that
\begin{equation}
\label{sreg}
\begin{array}{cc}
\exists N\in\mathbb{N}:\ \forall\alpha,\beta\in\mathbb{N}^n ,\\[0.2cm]
\displaystyle\sup_{\epsilon\in(0,1]}\epsilon^{N}\Vert x^\alpha\partial^\beta u_\epsilon\Vert_{L^\infty(\mathbb{R}^n)}<\infty .
\end{array}
\end{equation}
\end{definition}
\noindent We observe that any representative of an $\mathcal{S}$-regular generalized function, satisfies \eqref{sreg}, because for the elements of $\mathcal{N}_{\mathcal{S}}(\mathbb{R}^n)$ this property holds. In this way if we set
\[
\mathcal{E}^\infty_{\mathcal{S}}(\mathbb{R}^n)=\displaystyle\bigcup_{N\in\mathbb{N}}\mathcal{E}^{N}_{\mathcal{S}}(\mathbb{R}^n) ,
\]
where 
\[
\mathcal{E}^{N}_{\mathcal{S}}(\mathbb{R}^n)
=\displaystyle\{(u_\epsilon)_\epsilon\in\mathcal{E}[\mathbb{R}^n]:\ \forall\alpha,\beta\in\mathbb{N}^n,\ \sup_{\epsilon\in(0,1]}\epsilon^{N}\Vert x^\alpha\partial^\beta u_\epsilon\Vert_{L^\infty(\mathbb{R}^n)}<\infty\} ,
\]
we can define $\gss$ as the factor $\mathcal{E}^\infty_{\mathcal{S}}(\mathbb{R}^n)/\mathcal{N}_{\mathcal{S}}(\mathbb{R}^n)$. Obviously if $u=(u_\epsilon)_\epsilon +{\mathcal{N}}_{\mathcal{S}}(\mathbb{R}^n)\in{\mathcal{G}}^\infty_{\mathcal{S}}(\mathbb{R}^n)$, then the generalized function $(u_\epsilon)_\epsilon +\mathcal{N}(\mathbb{R}^n)$ belongs to ${\mathcal{G}}^\infty(\mathbb{R}^n)$. We observe that if $(u_\epsilon)_\epsilon\in{\mathcal{E}}^\infty_{\mathcal{S}}(\mathbb{R}^n)$, then $(u_\epsilon\widehat{\varphi_\epsilon}-u_\epsilon)_\epsilon\in{\mathcal{N}}_{\mathcal{S}}(\mathbb{R}^n)$. This result allows us to eliminate the mollifier $\varphi$ in Proposition 2.1 and Definition 2.6, in the case of $u\in\gss$. More precisely, for $u\in{\mathcal{G}}^\infty_{\mathcal{S}}(\mathbb{R}^n)$, ${\mathcal{F}}_\varphi u=(\widehat{u_\epsilon})_\epsilon +{\mathcal{N}}_{\mathcal{S}}(\mathbb{R}^n)$ and ${\mathcal{F}}^{\ast}_\varphi u=(\check{u_\epsilon})_\epsilon +{\mathcal{N}}_{\mathcal{S}}(\mathbb{R}^n)$.  
\begin{proposition}
${\mathcal{F}}_\varphi$ and $\mathcal{F}_\varphi^{\ast}$ map $\gss$ into $\gss$.
\end{proposition}
\begin{proof}
It remains to prove that if $(u_\epsilon)_\epsilon\in\mathcal{E}^\infty_{\mathcal{S}}(\mathbb{R}^n)$,  $({\mathcal{F}}_\varphi u_\epsilon)_\epsilon$ belongs to $\mathcal{E}^\infty_{\mathcal{S}}(\mathbb{R}^n)$. In particular we show that $(u_\epsilon)_\epsilon\in{\mathcal{E}}^N_{\mathcal{S}}(\mathbb{R}^n)$ implies $(\mathcal{F}_\varphi u_\epsilon)_\epsilon\in\mathcal{E}^{N}_{\mathcal{S}}(\mathbb{R}^n)$. In fact since, for all $y\in\mathbb{R}^n$ and $\epsilon\in(0,1]$, $|\partial^\gamma(-iy)^\beta\partial^{\alpha-\gamma} u_\epsilon(y)|$ is estimated by a constant multiplied by $\langle y\rangle^{-n-1}\epsilon^{-N}$, we have that
\begin{equation}
|\xi^{\alpha}\partial^\beta\mathcal{F}_\varphi u_\epsilon(\xi)|\le \sum_{\gamma\le\alpha}c_{\gamma}\int_{\mathbb{R}^n}\hskip-7pt|\partial^\gamma(-iy)^\beta\partial^{\alpha-\gamma} u_\epsilon(y)|dy\hskip1pt\le c \epsilon^{-N}.
\end{equation}
Analogously we obtain that $(\mathcal{F}^{\ast}_\varphi u_\epsilon)_\epsilon\in\mathcal{E}^\infty_{\mathcal{S}}(\mathbb{R}^n)$.
\end{proof}
\ni We conclude this section investigating the intersection of $\gss$ with $\mathcal{S}^{'}(\mathbb{R}^n)$. Inspired by \cite{ober}, Theorem 25.2, and \cite{hoer}, Theorem 16, we obtain the following result.
\begin{theorem}
\[
\gss\cap\mathcal{S}^{'}(\mathbb{R}^n)=\S .
\]
\end{theorem}
\begin{proof}
The inclusion $\S\subseteq\gss\cap\mathcal{S}^{'}(\mathbb{R}^n)$ is clear. 
Let $w\in\mathcal{S}^{'}(\mathbb{R}^n)$. We assume that $\imath(w)$ belongs to $\gss$. Denoting $w\ast\varphi_\epsilon$ by $w_\epsilon$, as a consequence of Proposition 1.2.21 in \cite{gro}, $(w_\epsilon)_\epsilon +{\mathcal{N}}(\mathbb{R}^n)\in{\mathcal{G}}^\infty(\mathbb{R}^n)\cap\mathcal{D}^{'}(\mathbb{R}^n)$, and then, from Theorem 25.2 in \cite{ober}, we already know that $w$ is a smooth function on $\mathbb{R}^n$. Moreover, since $(\widehat{w}\widehat{\varphi_\epsilon})_\epsilon$ belongs to $\S$ for every $\epsilon$ and $\widehat{\varphi}(0)=1$, taking $\epsilon$ as small as we want, we conclude that $\widehat{w}\in\mathcal{C}^\infty(\mathbb{R}^n)$. Now from the definition of ${\mathcal{G}}^\infty_{\mathcal{S}}(\mathbb{R}^n)$
\begin{equation}
\begin{array}{cc}
\exists N\in\mathbb{N}\ :\ \forall\alpha,\beta\in\mathbb{N}^n,\ \exists c>0:\ \forall\epsilon\in(0,1],\ \forall x\in\mathbb{R}^n ,\\[0.2cm]
|x^\alpha\partial^\beta w_\epsilon(x)|\le c\epsilon^{-N} .
\end{array}
\end{equation}
(2.9) implies the following statement:
\begin{equation}
\begin{array}{cc}
\exists N\in\mathbb{N}:\ \forall m\in\mathbb{N},\ \forall\alpha\in\mathbb{N}^n,\ \exists c>0:\ \forall\epsilon\in(0,1],\ \forall x\in\mathbb{R}^n ,\\[0.2cm]
|\langle x\rangle^{n+1}\Delta^{m}(x^\alpha w_\epsilon(x))|\le c\epsilon^{-N} .
\end{array}
\end{equation}
From (2.10) it follows that
\begin{equation}
\begin{array}{cc}
\exists N\in\mathbb{N}:\ \forall m\in\mathbb{N},\ \forall\alpha\in\mathbb{N}^n ,\ \exists c>0:\ \forall\epsilon\in(0,1] ,\\[0.2cm]
\Vert\Delta^{m}(x^\alpha w_\epsilon)\Vert_{L^{1}(\mathbb{R}^n)}\le c\epsilon^{-N} .
\end{array}
\end{equation}
Using Fourier transform
we conclude that
\begin{equation}
\Vert(\Delta^{m}(x^\alpha w_\epsilon))\widehat{\ }\ \Vert_{L^\infty(\mathbb{R}^n)}\le c\epsilon^{-N} ,
\end{equation}
thus 
\begin{equation}
\Vert\ |\xi|^{2m}\partial^\alpha(\widehat{w}\widehat{\varphi_\epsilon})\Vert_{L^\infty(\mathbb{R}^n)} \le c\epsilon^{-N} .
\end{equation}
We want to prove that (2.13) implies the following statement
\begin{equation}
\begin{array}{cc}
\forall m\in\mathbb{N},\ \forall\alpha\in\mathbb{N}^n ,\\[0.2cm]
\Vert\ |\xi|^{2m}\partial^\alpha\widehat{w}\Vert_{L^\infty(\mathbb{R}^n)}<\infty .
\end{array}
\end{equation}
If (2.14) holds we obtain our claim. We argue by induction. At first we verify that
\begin{equation}
\forall m\in\mathbb{N},\qquad\qquad\qquad \Vert |\xi|^{2m}\widehat{w}\Vert_{L^\infty(\mathbb{R}^n)}<\infty .\qquad\qquad\qquad
\end{equation}
Assume to the contrary that $\Vert\ |\xi|^{2\overline{m}}\widehat{w}\Vert_{L^\infty(\mathbb{R}^n)}=\infty$ for some $\overline{m}$. There exists a sequence $\{{\xi_j}\}\subset\mathbb{R}^n$, $|{\xi_j}|\to +\infty$, such that 
\begin{equation}
|{\xi_j}|^{2\overline{m}}|\widehat{w}({\xi_j})|\to +\infty .
\end{equation}
Since (2.13) holds for arbitrary $m$, we have in particular that 
\begin{equation}
|\xi_j|^{2\overline{m}}|\widehat{w}(\xi_j)||\widehat{\varphi}(\epsilon\xi_j)|\le c\epsilon^{-N}|\xi_j|^{-N-1} .
\end{equation}
Since $\widehat{\varphi}(0)=1$, there is $r>0$ such that $|\widehat{\varphi}(\xi)|\ge 1/2$ when $|\xi|\le r$. We define
\begin{equation}
\epsilon_j=r|\xi_j|^{-1} .
\end{equation}
Then $|\widehat{\varphi}(\epsilon_j\xi_j)|\ge 1/2$, and (2.17) implies that
\begin{equation}
|\xi_j|^{2\overline{m}}|\widehat{w}(\xi_j)|\frac{1}{2}\le c r^{-N}|\xi_j|^{-1} .
\end{equation}
This contradicts (2.16) because $|\xi_j|\to +\infty$.\\
\noindent In order to complete this proof we assume that
\begin{equation}
\begin{array}{cc}
\forall m\in\mathbb{N},\ \forall\alpha\in\mathbb{N}^n, |\alpha|\le k,\\[0.2cm]
\Vert|\xi|^{2m}\partial^\alpha\widehat{w}\Vert_{L^\infty(\mathbb{R}^n)}<\infty .
\end{array}
\end{equation}
We want to prove that
\begin{equation}
\begin{array}{cc}
\forall m\in\mathbb{N},\ \forall\alpha\in\mathbb{N}^n,\ |\alpha|\le k+1,\\[0.2cm]
\Vert|\xi|^{2m}\partial^\alpha\widehat{w}\Vert_{L^\infty(\mathbb{R}^n)}<\infty .
\end{array}
\end{equation}
As before we suppose that 
\begin{equation}
\exists \overline{m}\in\mathbb{N},\ \exists\overline{\alpha}\in\mathbb{N}^n,\ |\overline{\alpha}|\le k+1:\qquad \Vert|\xi|^{2\overline{m}}\partial^{\overline{\alpha}}\widehat{w}\Vert_{L^\infty(\mathbb{R}^n)}=\infty .
\end{equation}
Then we find a sequence $\{{\xi_j}\}_{j\in\mathbb{N}}\subset\mathbb{R}^n$, such that $|{\xi_j}|\to +\infty$ and 
\begin{equation} |{\xi_j}|^{2\overline{m}}|\partial^{\overline{\alpha}}\widehat{w}({\xi_j})|\to +\infty .
\end{equation}
We choose $m^{'}-s=\overline{m}$, with $2s>N$. From (2.13) we obtain
\begin{equation}
\begin{split}
&|{\xi_j}|^{2m^{'}}|\partial^{\overline{\alpha}}(\widehat{w}\widehat{\varphi_\epsilon})({\xi_j})|=\\
&|{\xi_j}|^{2m^{'}}\biggl|\sum_{\beta<\overline{\alpha}}\binom{\overline{\alpha}}{\beta}\partial^\beta\widehat{w}({\xi_j})\partial^{\overline{\alpha}-\beta}\widehat{\varphi}(\epsilon{\xi_j})\epsilon^{|\overline{\alpha}-\beta|}+\partial^{\overline{\alpha}}\widehat{w}({\xi_j})\widehat{\varphi}(\epsilon{\xi_j})\biggr|\le c\epsilon^{-N} .
\end{split}
\end{equation}
By induction hypothesis, all terms involving a derivative of order $\beta<\overline{\alpha}$ are bounded by a constant times $\epsilon^{-N}$. We arrive at
\begin{equation}
|\xi_j|^{2\overline{m}}|\partial^{\overline{\alpha}}\widehat{w}(\xi_j)\widehat{\varphi}(\epsilon\xi_j)|\le c^{'}\epsilon^{-N}|\xi_j|^{-2s} ,
\end{equation}
which leads to a contradiction with (2.23) as before.
\end{proof}
\section{Oscillatory integrals}
We describe the meaning and the most important properties of the integral, depending on a  real parameter $\epsilon$, of the type
\[
\int_{\mathbb{R}^n}e^{i\omega(x)}a_\epsilon(x)\ dx ,
\]
with phase function $\omega$ and amplitude $a_\epsilon$, satisfying suitable assumptions. In many proofs, we refer to \cite{bog, bogII, buz, saint} for details.
\begin{definition}
$\omega\in{\mathcal{C^\infty}}({\mathbb{R}^n}\setminus0)$ is a phase function of order $k>0$, $\omega\in{\Phi}^k(\mathbb{R}^n)$ for short, if it is real valued, positively homogeneous of order $k$, i.e. $\omega(tx)=t^k \omega(x)$ for $t>0$, and
\[
\nabla w(x)\neq 0\ \ for\ x\neq 0 .
\]
\end{definition}
\begin{definition}
Let $m\in\mathbb{R}$ and $l\in\mathbb{R}$. We denote by ${{\mathcal{A}}}^{m}_{l}(\mathbb{R}^n)$ or ${\mathcal{A}}^{m}_l$ for short, the set of all generalized amplitudes $(a_\epsilon)_{\epsilon\in(0,1]}\in\mathcal{E}[\mathbb{R}^n]$, satisfying the following requirement: for all $\alpha\in\mathbb{N}^n$, there exists $N\in\mathbb{N}$ such that
\[
\sup_{\epsilon\in(0,1]}\epsilon^{N}\Vert\langle x\rangle^{l|\alpha|-m}\partial^\alpha a_\epsilon\Vert_{L^\infty(\mathbb{R}^n)}<\infty ,
\]
where $\langle x\rangle=(1+|x|^2)^{\frac{1}{2}}$.
\end{definition}
\noindent We observe that for fixed $\epsilon\in(0,1]$, $a_\epsilon(x)$ belongs to the class considered in \cite{bog, buz}. It is immediate to verify that ${\mathcal{A}}^{m}_l$ is a linear space with the following properties:
\begin{itemize}
\item[i)] if $m\le m^{'}$ and $l\ge l^{'}$, then $\amp\subset{{\mathcal{A}}}^{m^{'}}_{l^{'}}$;
\item[ii)] if $(a_\epsilon)_\epsilon\in\amp$ and $(b_\epsilon)_\epsilon\in{{\mathcal{A}}}^{m^{'}}_{l^{'}}$, then $(a_\epsilon b_\epsilon)_\epsilon\in{{\mathcal{A}}}^{m+m^{'}}_{\min(l,l^{'})}$;
\item[iii)] if $(a_\epsilon)_\epsilon\in\amp$ then for all $\alpha\in\mathbb{N}^n$, $(\partial^\alpha a_\epsilon)_\epsilon\in{\mathcal{A}}^{m-l|\alpha|}_l$.
\end{itemize}
\begin{definition}
Let $N\in\mathbb{N}$. $\ampN(\mathbb{R}^n)$ or $\ampN$ for short, is the set of generalized amplitudes in $\amp$ such that for all $\alpha\in\mathbb{N}^n$  
\[
\sup_{\epsilon\in(0,1]}\epsilon^{N}\Vert\langle x\rangle^{l|\alpha|-m}\partial^\alpha a_\epsilon\Vert_{L^\infty(\mathbb{R}^n)}<\infty .
\]
The elements of $\cup_{N}{\mathcal{A}}^{m}_{l,N}$ are called regular generalized amplitudes.
\end{definition}
\ni It is clear that every ${\mathcal{A}}^{m}_{l,N}$ is a linear subspace of ${\mathcal{A}}^{m}_{l}$. Moreover
\begin{itemize}
\item[i)] if $m\le m^{'}$ , $l\ge l^{'}$ and $N\le N^{'}$ then $\ampN\subset{{\mathcal{A}}}^{m^{'}}_{l^{'},N^{'}}$;
\item[ii)] if $(a_\epsilon)_\epsilon\in\ampN$ and $(b_\epsilon)_\epsilon\in{{\mathcal{A}}}^{m^{'}}_{l^{'},N^{'}}$, then $(a_\epsilon b_\epsilon)_\epsilon\in{{\mathcal{A}}}^{m+m^{'}}_{\min(l,l^{'}),N+N^{'}}$;
\item[iii)] if $(a_\epsilon)_\epsilon\in\ampN$, for all $\alpha\in\mathbb{N}^n$, $(\partial^\alpha a_\epsilon)_\epsilon\in{\mathcal{A}}^{m-l|\alpha|}_{l,N}$.
\end{itemize}
The classes in \cite{bog, buz} are subsets of ${\mathcal{A}}^{m}_{l,0}$ with elements not depending on $\epsilon$.
Before stating the theorem on definition of oscillatory integral, we recall a useful lemma proved in \cite{bog}.
\begin{lemma}
Let $\omega\in\Phi^{k}(\mathbb{R}^n)$,  $a\in\mathcal{C}^\infty(\mathbb{R}^n\setminus 0)$ and $\chi\in\mathcal{C}^\infty_c(\mathbb{R}^n)$ vanishing in a neighbourhood of the origin. Then for every $N\in\mathbb{Z}^{+}$ there exists $c_N>0$, depending only on $\omega$ and $\chi$, such that for every $\mu>0$
\[
\biggr|\int_{\mathbb{R}^n}e^{i\mu^{k}w(y)}a(\mu y)\chi(y)dy\biggl|\le c_{N}\mu^{-k{N}}\ \displaystyle\sup_{|\alpha|\le{N}}\sup_{y\in supp\ \chi} \mu^{|\alpha|}|\partial^{\alpha}a(\mu y)| . 
\]
\end{lemma}
\begin{theorem}
Let $(a_\epsilon)_\epsilon\in\amp$, $\omega\in\Phi^{k}$ with $1-k<l\le 1$. Let $\psi$ be an arbitrary function of ${\mathcal{S}}({\mathbb{R}}^{n})$ with $\psi(0)=1$ and $\phi$ any function in ${\mathcal{C}}^{\infty}_c({\mathbb{R}}^n)$ such that $\phi(x)=1$ for $|x|\le 1$ and $\phi(x)=0$ for $|x|\ge 2$. Then for all $\epsilon\in(0,1]$ 
\[
\displaystyle\lim_{h\rightarrow 0^{+}} \int_{\mathbb{R}^n}e^{i\omega(x)}a_\epsilon(x)\psi(hx)dx= \displaystyle\lim_{j\rightarrow +\infty} \int_{\mathbb{R}^n}e^{i\omega(x)}a_\epsilon(x)\phi(2^{-j}x)dx ,
\]
i.e. the two limits exist in $\mathbb{C}$ and have the same value $I(\epsilon)$. Moreover
\begin{equation}
\begin{array}{cc}
\exists\overline{N}\in\mathbb{N},\ \overline{N}\ge\frac{m+n+1}{l+k-1}:\ \forall N\ge\overline{N},\  \exists M\in\mathbb{N}\ :\ \forall\epsilon\in(0,1],\\[0.3cm]
|I(\epsilon)|\le c\Vert a_\epsilon\Vert_{N}\le c^{'}\epsilon^{-{M}} ,
\end{array}
\end{equation}
where $\displaystyle\Vert a_\epsilon\Vert_{N}=\sup_{|\alpha|\le N}\Vert\langle x\rangle^{l|\alpha|-m}\partial^\alpha a_\epsilon\Vert_{L^\infty(\mathbb{R}^n)}$, $c$ does not depend on $\psi,\phi,a,\epsilon$, and $c^{'}$ does not depend on $\psi,\phi,\epsilon$.
\end{theorem}
\begin{remark}
If $l=0$ and $k=2$ we can choose $\overline{N}$ as the least integer greater than $m+n+2$. (See the proof of Theorem 0.1 in \cite{bog}, p. 14-16)
\end{remark}
\begin{remark} If $(a_\epsilon)_\epsilon\in\ampM$, the second line of (3.1) is valid, with $M$, co\-ming from the definition of regular generalized amplitude, independent of $N$.
\end{remark}
\ni Let $(a_\epsilon)_\epsilon\in\amp$, $\omega\in\Phi^{k}$ with $1-k<l\le 1$. For fixed $\epsilon$, we recall that the definition of oscillatory integral is given by 
\[
\begin{split}
\int_{\mathbb{R}^n}e^{i\omega(x)}a_\epsilon(x)dx&:=\lim_{h\rightarrow 0^{+}}\int_{\mathbb{R}^n}e^{i\omega(x)}a_\epsilon(x)\psi(hx)dx\\
&=\lim_{j\rightarrow+\infty}\int_{\mathbb{R}^n}e^{i\omega(x)}a_\epsilon(x)\phi(2^{-j}x)dx .
\end{split}
\]
Theorem 3.1 shows that the net $\displaystyle\biggl(\int_{\mathbb{R}^n}e^{i\omega(x)}a_\epsilon(x)dx\biggr)_\epsilon$ belongs to $\mathcal{E}_{o,M}$. We conclude this section with some useful pro\-per\-ties.
\begin{proposition}
Let $\omega\in\Phi^{k}(\mathbb{R}^n_x)$ a polynomial phase function, $(a_\epsilon)_{\epsilon}\in\mathcal{E}[\mathbb{R}^{2n}]$. We assume that:
\begin{itemize}
\item[i)] $(a_\epsilon(x,y))_\epsilon\in\amp(\mathbb{R}^n_x)$ with $1-k<l\le 1$;
\item[ii)] $\forall\beta\in\mathbb{N}^n$, $\exists m(\beta)\in\mathbb{R}$:\ $(\partial^\beta_y a_\epsilon(x,y))_\epsilon\in{\mathcal{A}}^{m(\beta)}_{l}(\mathbb{R}^n_x)$;
\item[iii)] $\forall\alpha,\ \beta\in\mathbb{N}^n$, $\exists M\in\mathbb{N}$:\ $\forall r>0 ,$
\[
\displaystyle\sup_{x\in\mathbb{R}^n, |y|\le r, \epsilon\in(0,1]}\hskip-25pt\epsilon^{M}\langle x\rangle^{l|\alpha|-m(\beta)}|\partial^\alpha_x\partial^\beta_y a_\epsilon(x,y)|<\infty .
\]
\end{itemize}
Then for all $\epsilon\in(0,1]$,\ $\displaystyle b_\epsilon(y)=\int_{\mathbb{R}^n}\hskip-7pt e^{i\omega(x)}a_\epsilon(x,y)dx\in\mathcal{C}^\infty(\mathbb{R}^n_y)$, for every $\beta\in\mathbb{N}^n$,\ $\displaystyle \partial^\beta b_\epsilon(y)=\int_{\mathbb{R}^n}\hskip-7pt e^{i\omega(x)}\partial^\beta_y a_\epsilon(x,y)dx $ and in addition $(b_\epsilon)_\epsilon\in\mathcal{E}_M(\mathbb{R}^n)$.
\end{proposition}
\begin{proof}
Choosing $\phi\in\mathcal{C}^\infty_c(\mathbb{R}^n)$ under the hypothesis of Theorem 3.1, we write for $\epsilon\in(0,1]$ and $y\in\mathbb{R}^n$
\begin{equation}
b_\epsilon(y)=\int_{\mathbb{R}^n}\hskip-7pt e^{i\omega(x)}a_\epsilon(x,y)dx=\lim_{j\rightarrow+\infty}\int_{\mathbb{R}^n}\hskip-7pt e^{i\omega(x)}a_\epsilon(x,y)\phi(2^{-j}x)dx .
\end{equation}
We define $\displaystyle b_{j,\epsilon}(y)\hskip-2pt=\hskip-4pt\int_{\mathbb{R}^n}\hskip-7pt e^{i\omega(x)}a_\epsilon(x,y)\phi(2^{-j}x)dx$ and we have $\displaystyle b_\epsilon(y)\hskip-2pt =\hskip-4pt \lim_{j\rightarrow +\infty}b_{j,\epsilon}(y)$.
From the hy\-po\-the\-sis $b_{j,\epsilon}\in\mathcal{C}^\infty(\mathbb{R}^n)$ and for all $\beta\in\mathbb{N}^n$
\[
\displaystyle\partial^\beta b_{j,\epsilon}(y)=\int_{\mathbb{R}^n}\hskip-7pt e^{i\omega(x)}\partial^\beta_y a_\epsilon(x,y)\phi(2^{-j}x)dx .
\]
In order to obtain the assertion, it suffices to show that for $\epsilon\in(0,1]$ and arbitrary $\beta\in\mathbb{N}^n$, $\{\partial^\beta b_{j,\epsilon}\}_j$ converges uniformly on compact sets of $\mathbb{R}^n$. Using Lemma 3.1, with $\chi(x)=\phi(x)-\phi(2x)$, we conclude that 
\[
\begin{array}{cc}
\forall N\in\mathbb{Z}^{+},\ \exists c_{N}>0:\ \forall y\in\mathbb{R}^n,\ \forall\epsilon\in(0,1],\\[0.3cm]
|\partial^\beta b_{j,\epsilon}(y)-\partial^\beta b_{j-1,\epsilon}(y)|\le c_{N}(\omega,\chi)2^{j(n-Nk)}\hskip-10pt\displaystyle\sup_{|\alpha|\le N, x\in supp\ \chi}\hskip-15pt 2^{j|\alpha|}|\partial^\alpha_x\partial^\beta_y a_\epsilon(2^{j}x,y)| .
\end{array}
\]
From hypothesis ii) and iii), we have that
\begin{equation}
\begin{array}{cc}
\forall N\in\mathbb{Z}^{+},\ N\ge\frac{m(\beta)+n+1}{l+k-1},\ \exists M\in\mathbb{N}:\ \forall r>0, \forall y\in\mathbb{R}^n,\ |y|\le r,\\[0.3cm]
\begin{split}
|\partial^\beta b_{j,\epsilon}(y)-\partial^\beta b_{j-1,\epsilon}(y)|&\le c_N\epsilon^{-M}2^{j(n-Nk)}\sup_{|\alpha|\le N,1/2\le |x|\le 2}\hskip-10pt 2^{j|\alpha|}\langle 2^{j}x\rangle^{m(\beta)-l|\alpha|}\\[0.2cm]
&\le c_{N}^{'}\epsilon^{-M}2^{j((1-l-k)N+m(\beta)+n)}\le c_{N}^{'}\epsilon^{-M}2^{-j} ,
\end{split}
\end{array}
\end{equation}
where the constants do not depend on $\epsilon$ and $j$. This result completes the proof.
\end{proof}
\ni We summarize in the following, without proofs, other properties used in this paper.
\begin{proposition}
Let $(a_\epsilon)_\epsilon\in\amp$ and $(b_\epsilon)_\epsilon\in{\mathcal{A}}^{p}_l$. Let $\omega\in\Phi^k$ with $1-k<l\le 1$ be a polynomial function. Then $(e^{-i\omega}\partial^{\alpha}(e^{i\omega}b_\epsilon))_\epsilon\in{\mathcal{A}}^{p-l|\alpha|}_l$ and
\[
\int_{\mathbb{R}^n}e^{i\omega(x)}\partial^\alpha a_\epsilon(x)b_\epsilon(x)dx= (-1)^{|\alpha|}\int_{\mathbb{R}^n}e^{i\omega(x)}a_\epsilon(x)e^{-i\omega(x)}\partial^\alpha(e^{i\omega(x)}b_\epsilon(x))dx .
\]
\end{proposition}
\begin{proposition}
Let $(a_\epsilon)_\epsilon\in{\mathcal{A}}^{m}_l(\mathbb{R}^n_x\times\mathbb{R}^p_y)$, $\omega\in\Phi^k(\mathbb{R}^n)$, $\eta\in\Phi^{k}(\mathbb{R}^p)$ with $1-k<l\le 1$. Suppose further that $\omega$ and $\eta$ are polynomials. Then we have $\biggl(\displaystyle\int_{\mathbb{R}^n}\hskip-7pt e^{i\omega(x)}a_\epsilon(x,y)dx\biggr)_{\hskip-4pt\epsilon}\in\amp(\mathbb{R}^p_y)$ ,\ $\biggl(\displaystyle\int_{\mathbb{R}^p}\hskip-7pt e^{i\eta(y)}a_\epsilon(x,y)dy\biggr)_{\hskip-4pt\epsilon}\in\amp(\mathbb{R}^n_x)$ and 
\[
\begin{split}
\int_{\mathbb{R}^{n+p}}\hskip-7pt e^{i(\omega(x)+\eta(y))}a_\epsilon(x,y)dx\ dy&=\\
\int_{\mathbb{R}^p}\hskip-7pt e^{i\eta(y)}\int_{\mathbb{R}^n}\hskip-7pt e^{i\omega(x)}a_\epsilon(x,y)dx\ dy &=\int_{\mathbb{R}^n}\hskip-7pt e^{i\omega(x)}\int_{\mathbb{R}^p}\hskip-7pt e^{i\eta(y)}a_\epsilon(x,y)dy\ dx .
\end{split}
\]
\end{proposition}
\begin{proposition}
Let $(a_\epsilon)_\epsilon\in\amp$ with $-1<l\le 1$. Then  
\[
\int_{\mathbb{R}^{2n}}\hskip-10pt e^{-iy\eta}a_\epsilon(y)\dslash y d\eta=\int_{\mathbb{R}^{2n}}\hskip-10pt e^{-iy\eta}a_\epsilon(\eta)\dslash y d\eta =a_\epsilon(0) . 
\]
\end{proposition}
\section{Symbols and amplitudes}
We introduce in this section the symbols and amplitudes, involving weight functions in their estimates, used in our definition of pseudo-differential operator acting on $\mathcal{G}_{\tau,\mathcal{S}}(\mathbb{R}^n)$. In the following for $f$ and $g$, real functions on $\mathbb{R}^{p}$, we write $f(z)\prec g(z)$ on $A\subseteq\mathbb{R}^{p}$, if there exists a positive constant $c$ such that $f(z)\le cg(z)$, for all $z\in A$. Here the constant $c$ may depend on parameters, indices, etc. possibly appearing in the expression of $f$ and $g$, but not on $z\in A$. If $f(z)\prec g(z)$ and $g(z)\prec f(z)$, we write $f(z)\sim g(z)$. We refer to \cite{beals, bog, bogII} for details in the classical setting. 
\begin{definition}
A continuous real function $\Lambda(z)$ on $\mathbb{R}^{2n}$ is a weight function iff
\begin{itemize}
\item[i)] there exists $\mu>0$ such that\quad $\langle z\rangle^{\mu}\prec\Lambda(z)\prec\langle z\rangle$\quad on $\mathbb{R}^{2n}$; 
\item[ii)] $\Lambda(z)\sim\Lambda(\zeta)$\quad on $A=\{(z,\zeta):\ |\zeta-z|\le\mu\Lambda(z)\}$;
\item[iii)] for all $t$ in $\mathbb{R}^{2n}$,\quad $\Lambda(tz)\prec\Lambda(z)$\quad on $\mathbb{R}^{2n}$,
\end{itemize}
where $tz=(t_1z_1,t_2z_2,..., t_{2n}z_{2n})$.
\end{definition} 
\noindent We recall that:\\
- from ii) it follows that $\Lambda(z)$ is temperate, i.e.
\begin{equation}
\Lambda(z)\prec\Lambda(\zeta)\langle z-\zeta\rangle ;
\end{equation}
- combining ii) and iii) we obtain that for all $t^{'}$ and $t^{''}$ in $\mathbb{R}^{2n}$
\begin{equation}
\Lambda(t^{'}z+t^{''}\zeta)\prec\Lambda(\zeta)\langle z-\zeta\rangle ;
\end{equation}
- from (4.1) it also follows that for any $s\in\mathbb{R}$
\begin{equation}
\Lambda(z)^{s}\prec\Lambda(\zeta)^{s}\langle z-\zeta\rangle^{|s|}
\end{equation}
and more precisely for $s<0$
\begin{equation}
\Lambda(z)^{s}\prec(1+\Lambda(\zeta)\langle z-\zeta\rangle^{-1})^s .
\end{equation}
In the next proposition we combine the preceding estimates in a more general form.
\begin{proposition}
Define $\displaystyle\tilde{\lambda}_s(x,y,\xi)=\begin{cases}\Lambda(x,\xi)^{s}\langle x-y\rangle^{s} , &\text{$s\ge 0 ,$}\\
(1+\Lambda(x,\xi)\langle x-y\rangle^{-1})^{s} , &\text{$s<0$} .
\end{cases}$\\ Then
\[
\Lambda(v^{'}x+v^{''}y,\xi)^{s}\prec\min(\tilde{\lambda}_s(x,y,\xi),\tilde{\lambda}_s(y,x,\xi)) ,
\]
for all $s\in\mathbb{R}$, $v^{'}, v^{''}\in\mathbb{R}^{n}$, provided $(x,y)\rightarrow (v^{'}x+v^{''}y,x-y)$ is an isomorphism on $\mathbb{R}^{2n}$.
\end{proposition}
\noindent Let us observe that starting from i) and ii), one can always find $\tilde\Lambda(z)\in\mathcal{C}^\infty(\mathbb{R}^{2n})$, with $\tilde{\Lambda}(z)\sim\Lambda(z)$, satisfying i), ii) and the property
\begin{equation}
\forall\gamma\in\mathbb{N}^{2n},\quad\quad\qquad |\partial^\gamma\tilde{\Lambda}(z)|\prec\tilde{\Lambda}^{1-|\gamma|}(z) .\quad\qquad\qquad\qquad
\end{equation}
In this way we do not lose generality if we assume that for the weight function (4.5) holds. We recall that if ${\mathcal{P}}\subset({\mathbb{R}}^{+}_0)^{2n}$ is a complete polyhedron with set of vertices $V({\mathcal{P}})$, and the estimate $0<\mu_0=(\min_{\gamma\in V(\mathcal{P})\setminus\{0\}}|\gamma|)\le\mu_1=(\max_{\gamma\in V(\mathcal{P})\setminus\{0\}}|\gamma|)\le\mu$ holds, where $\mu$ is the formal order of the polyhedron (for definition and details, see \cite{bog}, p. 20-22), then $\Lambda_{\mathcal{P}}(z)^{\frac{1}{\mu}}:=\big(\sum_{\gamma\in V(\mathcal{P})}z^{2\gamma}\big)^{\frac{1}{2\mu}}$ is an example of a weight function. For instance if $\mathcal{P}$ is the triangle in $\mathbb{R}^2$ of vertices $(k,0), (0,h), (0,0)$, $k,h\in\mathbb{N}\setminus\{0\}$, and $z=(x,\xi)$ then we get the weight function $\Lambda(z)=(1+x^{2k}+\xi^{2h})^{\frac{1}{2\max(k,h)}}$. Now we can introduce our sets of symbols.
\begin{definition}
Let $m\in\mathbb{R}$, $\rho\in(0,1]$, $\Lambda(z)$ be a weight function and $z=(x,\xi)\in\mathbb{R}^{2n}$. We denote by $\mathcal{S}^{m}_{\Lambda,\rho}(\mathbb{R}^{2n})$ or $\mathcal{S}^{m}_{\Lambda,\rho}$ for short, the set of symbols $(a_\epsilon)_\epsilon\in\mathcal{E}[\mathbb{R}^{2n}]$ fulfilling the condition:
\[
\begin{array}{cc}
\forall\alpha\in\mathbb{N}^{2n},\ \exists N\in\mathbb{N},\ \exists c>0:\ \forall\epsilon\in(0,1],\ \forall z\in\mathbb{R}^{2n},\\[0.3cm]
|\partial^\alpha a_\epsilon(z)|\le c\Lambda(z)^{m-\rho|\alpha|}\epsilon^{-N} .\end{array}
\]
\end{definition}
\begin{remark}
It follows from Definition 4.2 that
\begin{equation}
|\partial^\alpha a_\epsilon(z)|\le c\Lambda(z)^{m-\rho|\alpha|}\epsilon^{-N}\le c^{'}\langle z\rangle^{m_+-\mu\rho|\alpha|}\epsilon^{-N} ,\end{equation}
where $m_+=\max(0,m)$. In other words $\mathcal{S}^{m}_{\Lambda,\rho}\subset\mathcal{A}^{m_+}_{\mu\rho}$.
\end{remark}
\begin{definition}
Let $m\in\mathbb{R}$, $\rho\in(0,1]$, $\Lambda$ be a weight function and $N\in\mathbb{N}$. We denote by $\mathcal{S}^{m}_{\Lambda,\rho,N}$, the subset of the elements of $\mathcal{S}^{m}_{\Lambda,\rho}$ satisfying the property:
\[
\begin{array}{cc}
\forall\alpha\in\mathbb{N}^{2n}, \exists c>0:\ \forall\epsilon\in(0,1],\ \forall z\in\mathbb{R}^{2n},\\[0.3cm]
|\partial^\alpha a_\epsilon(z)|\le c\Lambda(z)^{m-\rho|\alpha|}\epsilon^{-N} .\end{array}
\]
The symbols of $\cup_{N}{\mathcal{S}}^{m}_{\Lambda,\rho,N}$ are called regular.
\end{definition}
\ni Let us emphasize that in Definition 4.3 the integer $N$ is independent of $\alpha$.\\
\noindent In \cite{bog, bogII} the authors consider as symbols the elements of ${\mathcal{S}}^{m}_{\Lambda,\rho,0}$ which do not depend on $\epsilon$. Moreover every $(a_{{\epsilon}})_\epsilon\in\mathcal{S}^{m}_{\Lambda,\rho}$, for fixed $\epsilon$, can be considered as a symbol in \cite{bog, bogII}. The presence of the parameter $\epsilon$ tending to $0$ in the definition of symbols, resembles the constructions of semiclassical analysis, (see \cite{shu}, pages 432-448, and \cite{mar}). However, it is simple to see that if, for instance, $a(x,\xi)$ is a classical symbol with weight $\Lambda(x,\xi)=\langle(x,\xi)\rangle$ and $\rho=1$, the semiclassical variant $a_\epsilon(x,\xi)=a(x,\epsilon\xi)$ is not, in general, a regular symbol in the sense of Definition 4.3. In fact for $|\beta|$ large enough
\[
|\partial^\beta_x a_\epsilon(x,\xi)|\le c\langle(x,\epsilon\xi)\rangle^{m-|\beta|}\le c\epsilon^{m-|\beta|}\langle (x,\xi)\rangle^{m-|\beta|} .
\]
\begin{definition}
Let $m\in\mathbb{R}$, $\rho\in(0,1]$ and $\Lambda$ be a weight function. We denote by $\mathcal{N}^{m}_{\Lambda,\rho}$, the subset of the elements of $\mathcal{S}^{m}_{\Lambda,\rho}$ satisfying the property:
\[
\begin{array}{cc}
\forall\alpha\in\mathbb{N}^{2n},\ \forall q\in\mathbb{N},\ \exists c>0:\ \forall\epsilon\in(0,1],\ \forall z\in\mathbb{R}^{2n},\\[0.3cm]
|\partial^\alpha a_\epsilon(z)|\le c\Lambda(z)^{m-\rho|\alpha|}\epsilon^{q} .
\end{array}
\]
We call them negligible symbols of order $m$.
\end{definition}
\ni Now we turn to the extension of the classical theory and list some basic results, the proofs of which are elementary and thus omitted.
\begin{proposition}
For a fixed weight function $\Lambda$, we have that:
\begin{itemize}
\item[i)] if $m\le m^{'}$, $\rho\ge\rho^{'}$, $N\le N^{'}$, then  
$\mathcal{S}^{m}_{\Lambda,\rho,N}\subset\mathcal{S}^{m^{'}}_{\Lambda,\rho^{'},N^{'}}$;
\item[ii)] if $(a_\epsilon)_\epsilon\in{\mathcal{S}}^{m}_{\Lambda,\rho,N}$and $(b_\epsilon)_\epsilon\in\mathcal{S}^{m^{'}}_{\Lambda,\rho,N^{'}}$ then $(a_\epsilon b_\epsilon)_\epsilon\in\mathcal{S}^{m+m^{'}}_{\Lambda,\rho,N+N^{'}}$ and\\
$(a_\epsilon +b_\epsilon)_\epsilon\in\mathcal{S}^{\max(m,m^{'})}_{\Lambda,\rho,\max(N,N^{'})};$
\item[iii)] if $(a_\epsilon)_\epsilon\in\mathcal{S}^{m}_{\Lambda,\rho,N}$ then for all $\alpha\in\mathbb{N}^{2n}$, $(\partial^\alpha a_\epsilon)_\epsilon\in\mathcal{S}^{m-\rho|\alpha|}_{\Lambda,\rho,N};$
\item[iv)] if $(a_\epsilon)_\epsilon\in\mathcal{S}^{m}_{\Lambda,\rho,N}$ then for all $\zeta\in\mathbb{R}^{2n}$ $(T_\zeta a_\epsilon(z))_\epsilon=(a_\epsilon(z-\zeta))_\epsilon\in\mathcal{S}^{m}_{\Lambda,\rho,N}$.
\end{itemize}
The previous four statements hold, without the third subscripts $N$, $N'$ etc., for the elements of $\mathcal{S}^{m}_{\Lambda,\rho}$ and $\mathcal{N}^{m}_{\Lambda,\rho}$ respectively.
\end{proposition} 
\begin{definition}
We call smoothing symbols the elements of
\[
\mathcal{S}^{-\infty}=\bigcup_{N\in\mathbb{N}}\mathcal{S}^{-\infty}_{\Lambda,\rho,N}=\bigcup_{N\in\mathbb{N}}\bigcap_{m\in\mathbb{R}}\mathcal{S}^{m}_{\Lambda,\rho,N} .
\]
\end{definition}
\begin{proposition}
$(a_\epsilon)_\epsilon$ is a smoothing symbol iff there exists $N\in\mathbb{N}$ such that for all $\alpha, \beta\in\mathbb{N}^{2n}$
\[
\sup_{\epsilon\in(0,1]}\epsilon^{N}\Vert z^\alpha\partial^\beta a_\epsilon\Vert_{L^\infty(\mathbb{R}^{2n})}<+\infty .
\]
\end{proposition}
\begin{proof}
In order to prove the necessity of the condition, it suffices to choose $m\in\mathbb{R}$, with $\mu m<-|\alpha|$ and to observe that
\[
|\partial^\beta a_\epsilon(z)|\le c\Lambda(z)^{m-\rho|\beta|}\epsilon^{-N}\le c^{'}\langle z\rangle^{\mu m}\epsilon^{-N} ,
\]
where $N$ does not depend on $m$ and $\beta$. On the contrary suppose that $(a_\epsilon)_\epsilon$ satisfies an estimate of rapidly decreasing type. For arbitrary weight function $\Lambda$, $\rho\in(0,1]$ and $m\in\mathbb{R}$, we have, if $m-\rho|\beta|\le 0$
\[
|\partial^\beta a_\epsilon(z)|\le c\langle z\rangle^{m-\rho|\beta|}\epsilon^{-N}\le c^{'}\Lambda(z)^{m-\rho|\beta|}\epsilon^{-N} ,
\]
and if $m-\rho|\beta|>0$
\[
|\partial^\beta a_\epsilon(z)|\le c\langle z\rangle^{\mu(m-\rho|\beta|)}\epsilon^{-N}\le c^{'}\Lambda(z)^{m-\rho|\beta|}\epsilon^{-N} ,
\]
where $\mu$ depends on $\Lambda$.
\end{proof}
\ni As a consequence of this proposition, the definition of $\mathcal{S}^{-\infty}$ is independent of the weight function and $\rho\in(0,1]$.
\begin{remark}
$\displaystyle{\mathcal{N}}^{-\infty}:=\displaystyle  {\bigcap}_{m\in\mathbb{R}}{\mathcal{N}}^{m}_{\Lambda,\rho}$ is characterized by the following statement:
\[
\begin{array}{cc}
\forall\alpha ,\beta\in\mathbb{N}^{2n},\ \forall q\in\mathbb{N},\\[0.2cm]
\displaystyle\sup_{\epsilon\in(0,1]}\epsilon^{-q}\Vert z^\alpha\partial^\beta a_\epsilon\Vert_{L^\infty({\mathbb{R}^{2n}})}<+\infty .
\end{array}
\]
\end{remark}
\ni We record now some examples of symbols. It is obvious that if $a(z)$ is a symbol of the type considered in \cite{bog, bogII}, $a_\epsilon(z):=a(z)\epsilon^b$, $b\in\mathbb{R}$, is a regular symbol according to Definition 4.3. Other examples are given by the following polynomials.
\begin{proposition}
Let $\Lambda$ be an arbitrary weight function and let $a_\epsilon(z)=\sum_{\alpha\in\mathcal{A}}c_{\alpha,\epsilon}z^\alpha$ be a polynomial with coefficients in $\mathcal{E}_o$, where $\mathcal{A}$ is a finite subset of $\mathbb{N}^{2n}$. There exists $r\in\mathbb{R}$ depending on $\Lambda$ and $\mathcal{A}$ such that the following statements hold:
\begin{itemize}
\item[i)] if every $(c_{\alpha,\epsilon})_\epsilon\in\mathcal{E}_{o,M}$ then $(a_\epsilon)_\epsilon\in\mathcal{S}^r_{\Lambda,1,N}$ for a suitable $N\in\mathbb{N}$ depending on the coefficients $(c_{\alpha,\epsilon})_\epsilon$;
\item[ii)] if every $(c_{\alpha,\epsilon})_\epsilon\in\mathcal{N}_o$ then $(a_\epsilon)_\epsilon\in\mathcal{N}^r_{\Lambda,1}$.
\end{itemize}
\end{proposition}
\begin{proof}
We begin by recalling that from Definition 4.1, $\langle z\rangle^\mu\prec\Lambda(z)\prec\langle z\rangle$ for $\mu >0$. In the sequel we write $k=\max_{\alpha\in\mathcal{A}}|\alpha|$ and $N=\max_{\alpha\in\mathcal{A}} N_\alpha$, where $|c_{\alpha,\epsilon}|\le c_\alpha\epsilon^{-N_\alpha}$ for all $\epsilon\in(0,1]$. Therefore,
\begin{equation}
|a_\epsilon(z)|=|\sum_{\alpha\in\mathcal{A}}c_{\alpha,\epsilon}z^\alpha|\le c_1\epsilon^{-N}\sum_{\alpha\in\mathcal{A}}\langle z\rangle^{|\alpha|}\le c_2\epsilon^{-N}\Lambda(z)^{\frac{k}{\mu}},\qquad\quad z\in\mathbb{R}^{2n},\ \epsilon\in(0,1],
\end{equation}
and for $\gamma\in\mathbb{N}^{2n}$, $\gamma\neq 0$,
\begin{equation}
|\partial^\gamma a_\epsilon(z)|\le c_1\epsilon^{-N}\sum_{\substack{\alpha\in\mathcal{A}\\ \alpha\ge\gamma}}\langle z\rangle^{|\alpha|-|\gamma|}\le c_2\epsilon^{-N}\Lambda(z)^{\frac{k}{\mu}-|\gamma|},\qquad\quad z\in\mathbb{R}^{2n},\ \epsilon\in(0,1].
\end{equation}
In conclusion $(a_\epsilon)_\epsilon\in\mathcal{S}^r_{\Lambda,1,N}$ with $r=k/\mu$. Analogously we can prove that $(c_{\alpha,\epsilon})_\epsilon\in\mathcal{N}_o$, for every $\alpha\in\mathcal{A}$, implies $(a_\epsilon)_\epsilon\in\mathcal{N}^r_{\Lambda,1}$.
\end{proof}
\ni As a consequence of Proposition 4.4 we can associate to the polynomial $\sum_{\alpha\in\mathcal{A}}c_\alpha z^\alpha\in{\overline{\mathbb{C}}}[z]$ the class $(\sum_{\alpha\in\mathcal{A}}c_{\alpha,\epsilon}z^\alpha)_\epsilon+\mathcal{N}^{r}_{\Lambda,1}$ as an element of the factor ${\mathcal{S}}^{r}_{\Lambda,1,N}/{\mathcal{N}}^{r}_{\Lambda,1}$. Let us consider the more general set of polynomials $\sum_{\beta\in\mathcal{B}}c_\beta(x)\xi^\beta\in\mathcal{G}(\mathbb{R}^n)[\xi]$ where $\mathcal{B}$ is a finite subset of $\mathbb{N}^n$. We obtain the following result.
\begin{proposition}
Let $\Lambda(x,\xi)$ be an arbitrary weight function.\\
If $b_\epsilon(x,\xi)=\sum_{\beta\in\mathcal{B}}c_{\beta,\epsilon}(x)\xi^\beta$, with $(c_{\beta,\epsilon})_\epsilon\in\mathcal{E}_{M}(\mathbb{R}^n)$, belongs to $\mathcal{S}^{r}_{\Lambda,\rho}$ then every $(c_{\beta,\epsilon})_\epsilon$ is a polynomial in $x$ with coefficients in $\mathcal{E}_{o,M}$.\\
If $b_\epsilon(x,\xi)=\sum_{\beta\in\mathcal{B}}c_{\beta,\epsilon}(x)\xi^\beta$, with $(c_{\beta,\epsilon})_\epsilon\in\mathcal{N}(\mathbb{R}^n)$, belongs to $\mathcal{N}^{r}_{\Lambda,\rho}$ then every $(c_{\alpha,\epsilon})_\epsilon$ is a polynomial in $x$ with coefficients in $\mathcal{N}_{o}$.
\end{proposition}
\begin{proof}
For every fixed $\epsilon\in(0,1]$ we are under the assumptions of Proposition 1.2 in \cite{bog}. Therefore, for $|\gamma|>r/\rho$, for all $\beta\in\mathcal{B}$ and $x\in\mathbb{R}^n$, $\partial^\gamma c_{\beta,\epsilon}(x)=0$, i.e. $c_{\beta,\epsilon}(x)=\sum_{|\gamma|\le r/\rho}\frac{\partial^\gamma c_{\beta,\epsilon}(0)}{\gamma !}x^\gamma$. At this point the conclusion is obvious.
\end{proof}
\ni From the previous proposition we have that if $\sum_{\beta\in\mathcal{B}}c_{\beta}(x)\xi^\beta\in\mathcal{G}(\mathbb{R}^n)[\xi]$ is an element of the factor $\mathcal{S}^{r}_{\Lambda,\rho}/\mathcal{N}^{r}_{\Lambda,\rho}$ then every $c_\beta(x)$ belongs to $\overline{\mathbb{C}}[x]$.\\[0.2cm]
\ni Formal series and asymptotic expansions play a basic role in the classical theory of pseudo-differential operators. In the following we generalize these concepts to our context.
\begin{definition}
We denote by $F{\mathcal{S}}^{m}_{\Lambda,\rho,N}$ the set of formal series $\sum_{j=0}^{\infty}(a_{j,\epsilon})_\epsilon$, such that for all $j\in\mathbb{N}$, $(a_{j,\epsilon})_\epsilon\in\mathcal{S}^{m_j}_{\Lambda,\rho,N_j}$ and the following assumptions hold:\\
- $m_0=m$, the sequence $\{m_j\}_{j\in\mathbb{N}}$ is decreasing with $m_j\rightarrow -\infty$;\\
- for all $j$, $N_j\le N$.
\end{definition}
\begin{definition}
Let $(a_\epsilon)_\epsilon\in\mathcal{E}[\mathbb{R}^{2n}]$ and $\sum_{j=0}^{\infty}(a_{j,\epsilon})_\epsilon\in F{\mathcal{S}}^{m}_{\Lambda,\rho,N}$.\\ \ni$\sum_{j=0}^{\infty}(a_{j,\epsilon})_\epsilon$ is the asymptotic expansion of $(a_\epsilon)_\epsilon$ and we write $(a_\epsilon)_\epsilon\sim\sum_{j=0}^{\infty}(a_{j,\epsilon})_\epsilon$ iff
\[
\forall r\ge1,\quad\quad\qquad \biggl(a_\epsilon-\sum_{j=0}^{r-1}a_{j,\epsilon}\biggr)_\epsilon\in\mathcal{S}^{m_r}_{\Lambda,\rho,N}.\qquad\qquad\qquad
\]
\end{definition}
\begin{theorem}
For any $\sum_{j=0}^{\infty}(a_{j,\epsilon})_\epsilon\in F{\mathcal{S}}^{m}_{\Lambda,\rho,N}$ there exists $(a_\epsilon)_\epsilon\in\mathcal{S}^{m}_{\Lambda,\rho,N}$ such that $(a_\epsilon)_\epsilon\sim\sum_{j=0}^{\infty}(a_{j,\epsilon})_\epsilon$. Moreover if $(a^{'}_\epsilon)_\epsilon\sim\sum_{j=0}^{\infty}(a_{j,\epsilon})_\epsilon$ then the difference $(a_{\epsilon}-a^{'}_\epsilon)_\epsilon$ belongs to $\mathcal{S}^{-\infty}_{\Lambda,\rho,N}$.
\end{theorem}
\begin{proof}
We start by considering $\psi\in\mathcal{C}^{\infty}(\mathbb{R})$ such that $\psi(t)=0$ for $t\le 1$ and $\psi(t)=1$ for $t\ge 2$. We define for $j\in\mathbb{N}$, $\lambda_j\in\mathbb{R}^{+}$
\begin{equation}
b_{j,\epsilon}(z)=\psi(\lambda_j\Lambda(z))a_{j,\epsilon}(z) .
\end{equation}
Our aim is to verify that for a suitable decreasing sequence of strictly positive numbers $\lambda_j$, such that $\lambda_j\rightarrow 0$, the following sum
\begin{equation}
(a_\epsilon)_\epsilon=\sum_{j=0}^{\infty}(b_{j,\epsilon})_\epsilon
\end{equation}
defines an element of $\mathcal{S}^{m}_{\Lambda,\rho,N}$ and it has $\sum_{j=0}^{\infty}(a_{j,\epsilon})_\epsilon$ as asymptotic expansion. We observe first that the sum in (4.10) is locally finite, since
\begin{equation}
supp\ b_{j,\epsilon}\subseteq supp\ \psi(\lambda_j\Lambda)\subseteq D_j=\{z\in\mathbb{R}^{2n}:\ \lambda_j\Lambda(z)\ge 1\}
\end{equation}
and, by definition of the weight function,
\begin{equation}
\mathbb{R}^{2n}\setminus D_j\subseteq\{z\in\mathbb{R}^{2n}:\ |z|\le c\lambda_j^{-1/\mu}\} .
\end{equation} 
Moreover for all $j$, $D_{j+1}\subseteq D_j$ and for $k\neq 0$, $supp\ {\psi}^{(k)}(\lambda_j\Lambda(z))$ is contained in the region $\{z\in\mathbb{R}^{2n}:\ 1\le\lambda_j\Lambda(z)\le 2\}$.\\
In order to obtain the claim we show that for all $j\in\mathbb{N}$, $ (b_{j,\epsilon})_\epsilon\in\mathcal{S}^{m_j}_{\Lambda,\rho,N}$, 
and more precisely for all $\gamma\in\mathbb{N}^{2n}$
\begin{equation}
|\partial^\gamma b_{j,\epsilon}(z)|\le c_{j,\gamma}\Lambda(z)^{m_j-\rho|\gamma|}\epsilon^{-N} ,
\end{equation}
with $c_{j,\gamma}$ independent of $\lambda_j$. To prove (4.13) we observe that
\begin{equation}
\partial^\gamma b_{j,\epsilon}(z)=\sum_{k\le|\gamma|}\psi^{(k)}(\lambda_j\Lambda(z))\tilde{b}_{j,k,\gamma,\epsilon}(z) ,
\end{equation}
where $({\tilde{b}}_{j,k,\gamma,\epsilon})_\epsilon\in\mathcal{S}^{m_j-\rho|\gamma|}_{\Lambda,\rho,N}$ with estimates independent of $\lambda_j$. In fact, if $|\gamma|=0$, (4.14) holds with estimates of required type, and by induction, if (4.14) is valid for $|\gamma|=h$, we have
\[
\begin{split}
\partial_{z_j}\partial^\gamma b_{j,\epsilon}(z)&\hskip-2pt=\hskip-5pt\sum_{k\le|\gamma|}\hskip-4pt\psi^{(k+1)}(\lambda_j\Lambda(z))\lambda_j\partial_{z_j}\Lambda(z){\tilde{b}}_{j,k,\gamma,\epsilon}(z)\\
&+\sum_{k\le|\gamma|}\hskip-4pt\psi^{(k)}(\lambda_j\Lambda(z))\partial_{z_j}\tilde{b}_{j,k,\gamma,\epsilon}(z) ,
\end{split}
\]
where in the second sum, $(\partial_{z_j}\tilde{b}_{j,k,\gamma,\epsilon})_\epsilon$ belongs to ${\mathcal{S}}^{m_j-\rho|\gamma+1|}_{\Lambda,\rho,N}$, with estimates independent of $\lambda_j$. We obtain the same result for the first sum, considering that the weight function is an element of ${\mathcal{S}}^{1}_{\Lambda,1,0}$ and as a consequence for $z\in supp\ \psi^{(k)}(\lambda_j\Lambda)$, we can write $\lambda_j\partial^\beta\partial_{z_j}\Lambda(z)\prec\lambda_j\Lambda(z)^{1-\rho(1+|\beta|)}\prec\Lambda(z)^{-\rho(1+|\beta|)}$.\\
\ni
We prove now that for every $r\ge 1$, $\displaystyle
\big(a_{\epsilon}-\sum_{j=0}^{r-1}a_{j,\epsilon}\big)_\epsilon\in\mathcal{S}^{m_r}_{\Lambda,\rho,N}$. This will give $(a_\epsilon)_\epsilon\in\mathcal{S}^{m}_{\Lambda,\rho,N}$ and $(a_{\epsilon})_\epsilon\sim\sum_{j=0}^{\infty}(a_{j,\epsilon})_\epsilon$. We write
\[
\partial^{\gamma}\biggr(a_{\epsilon}-\sum_{j=0}^{r-1}a_{j,\epsilon}\biggl)=\partial^\gamma\sum_{j=0}^{r-1}(b_{j,\epsilon}-a_{j,\epsilon})+\sum_{j=r}^{j=s}\partial^\gamma b_{j,\epsilon}+\sum_{j=s+1}^{\infty}\partial^\gamma b_{j,\epsilon} .
\]
The first term of the right-hand side is an element of $\mathcal{S}^{m_r-\rho|\gamma|}_{\Lambda,\rho,N}$ since $(a_{j,\epsilon})_\epsilon\in\mathcal{S}^{m}_{\Lambda,\rho,N}$ and $\psi-1\in\mathcal{C}^\infty_{c}(\mathbb{R})$; in the second term $(\partial^{\gamma}b_{j,\epsilon})_\epsilon\in\mathcal{S}^{m_j-\rho|\gamma|}_{\Lambda,\rho,N}\subset\mathcal{S}^{m_r-\rho|\gamma|}_{\Lambda,\rho,N}$. It remains to estimate the third term of the right-hand side. We assume $s\ge|\gamma|$ such that $m_s\le m_r-1$, and remembering (4.13), we choose $\lambda_j$ satisfying the following condition
\begin{equation}
c_{j,\gamma}\lambda_j\le 2^{-j},\qquad{\rm{for}}\ |\gamma|\le j.
\end{equation}
We obtain then for $j>s$, by using (4.11),
\begin{equation}
\begin{split}
|\partial^\gamma b_{j,\epsilon}(z)|&\le c_{j,\gamma}\Lambda(z)^{m_j-\rho|\gamma|}\epsilon^{-N}\le 2^{-j}\lambda_j^{-1}\Lambda(z)^{m_j-\rho|\gamma|}\epsilon^{-N}\\
&\le 2^{-j}\lambda_j^{-1}\Lambda(z)^{-1}\Lambda(z)^{m_r-\rho|\gamma|}\epsilon^{-N}\le 2^{-j}\Lambda(z)^{m_r-\rho|\gamma|}\epsilon^{-N} .
\end{split}
\end{equation}
Suppose finally that $(a^{'}_\epsilon)_\epsilon\sim\sum_{j=0}^{\infty}(a_{j,\epsilon})_\epsilon$. Then, for every $r\ge 1$
\[
(a_\epsilon-a^{'}_\epsilon)_\epsilon=\biggl(a_\epsilon-\sum_{j=0}^{r-1}a_{j,\epsilon}\biggr)_\epsilon-\biggr(a^{'}_\epsilon-\sum_{j=0}^{r-1}a_{j,\epsilon}\biggl)_\epsilon\in\mathcal{S}^{m_r}_{\Lambda,\rho,N}
\]
and $\displaystyle (a_\epsilon-a^{'}_\epsilon)_\epsilon\in \cap_{m}\mathcal{S}^{m}_{\Lambda,\rho,N}=\mathcal{S}^{-\infty}_{\Lambda,\rho,N}$
\end{proof}
\ni In Definition 4.6 and Theorem 4.1 we can omit the assumption $m_{j+1}\le m_j$. In this case the meaning of $(a_\epsilon)_\epsilon\sim\sum_{j=0}^{\infty}(a_{j,\epsilon})_\epsilon$ is $\big(a_\epsilon-\sum_{j=0}^{r-1}a_{j,\epsilon}\big)_\epsilon\in\mathcal{S}^{\overline{m}_r}_{\Lambda,\rho,N}$ for every $r\ge 1$, where $\overline{m}_r=\max_{j\ge r}(m_j)$. 
\begin{definition}
Let $m\in\mathbb{R}$, $\rho\in(0,1]$ and $\Lambda$ be a weight function. We denote by $\overline{\mathcal{S}}^{m}_{\Lambda,\rho}(\mathbb{R}^{3n})$ or $\overline{\mathcal{S}}^{m}_{\Lambda,\rho}$ for short, the set of all amplitudes $(a_\epsilon(x,y,\xi))_\epsilon\in\mathcal{E}[\mathbb{R}^{3n}]$ fulfilling the condition
\begin{equation}
\begin{array}{cc}
\forall\alpha, \beta, \gamma\in\mathbb{N}^n,\ \exists N\in\mathbb{N},\ \exists c>0:\ \forall\epsilon\in(0,1],\ \forall(x,y,\xi)\in\mathbb{R}^{3n},\\[0.3cm]
|\partial^\alpha_\xi\partial^\beta_x\partial^\gamma_y a_\epsilon(x,y,\xi)|\le c \lambda_{m,m^{'},\alpha,\beta,\gamma}(x,y,\xi)\epsilon^{-N} ,\end{array}
\end{equation}
with
\begin{equation}
\lambda_{m,m^{'},\alpha,\beta,\gamma}(x,y,\xi)=\Lambda(x,\xi)^{m}\langle x-y\rangle^{m^{'}}\big(1+\Lambda(x,\xi)\langle x-y\rangle^{-m^{'}}\big)^{-\rho|\alpha+\beta+\gamma|} ,
\end{equation}
for a suitable $m^{'}\in\mathbb{R}$ independent of derivatives.\\
\ni Interchanging $\exists N\in\mathbb{N}$ with $\forall\alpha,\beta,\gamma\in\mathbb{N}^n$ in $(4.17)$ we define the subset $\overline{\mathcal{S}}^{m}_{\Lambda,\rho,N}$ of $\overline{\mathcal{S}}^{m}_{\Lambda,\rho}$. The amplitudes of $\cup_N\overline{\mathcal{S}}^{m}_{\Lambda,\rho,N}$ are called regular.
\end{definition}
\ni The smooth amplitudes discussed in \cite{bog, bogII} are elements of $\overline{\mathcal{S}}^{m}_{\Lambda,\rho,0}$ and for fixed $\epsilon$, $a_\epsilon(x,y,\xi)$ can be considered as in \cite{bog, bogII}.
\begin{definition}
Let $m\in\mathbb{R}$, $\rho\in(0,1]$ and $\Lambda$ be a weight function. We denote by $\overline{\mathcal{N}}^{m}_{\Lambda,\rho}$ the subset of $\overline{\mathcal{S}}^{m}_{\Lambda,\rho}$ of all the elements satisfying the property
\begin{equation}
\begin{array}{cc}
\forall\alpha, \beta, \gamma\in\mathbb{N}^n,\ \forall q\in\mathbb{N},\ \exists c>0:\ \forall\epsilon\in(0,1],\ \forall(x,y,\xi)\in\mathbb{R}^{3n},\\[0.3cm]
|\partial^\alpha_\xi\partial^\beta_x\partial^\gamma_y a_\epsilon(x,y,\xi)|\le c\lambda_{m,m^{'},\alpha,\beta,\gamma}(x,y,\xi)\epsilon^{q} .\end{array}
\end{equation}
We call the elements of $\overline{\mathcal{N}}^{m}_{\Lambda,\rho}$ negligible amplitudes of order $m$.
\end{definition}
\ni In the sequel we collect, without proofs, some useful results.
\begin{proposition}
\label{simmetria}
The estimate in (4.17) is equivalent to each one of the fol\-lo\-wing two for suitable values of $m^{'}\in\mathbb{R}$:
\[
\begin{split}
|\partial^\alpha_\xi\partial^\beta_x\partial^\gamma_y a_\epsilon(x,y,\xi)|&\le c  \lambda_{m,m^{'},\alpha,\beta,\gamma}(y,x,\xi)\epsilon^{-N} ,\\[0.3cm]
|\partial^\alpha_\xi\partial^\beta_x\partial^\gamma_y a_\epsilon(x,y,\xi)|&\le c\min\{\lambda_{m,m^{'},\alpha,\beta,\gamma}(x,y,\xi), \lambda_{m,m^{'},\alpha,\beta,\gamma}(y,x,\xi)\}\epsilon^{-N} ,
\end{split}
\]
where $c$ does not depend on $\epsilon$.
\end{proposition}
\ni It is simple to prove that if $(b_\epsilon(x,y,\xi))_\epsilon\in\overline{\mathcal{S}}^{m}_{\Lambda,\rho}$ then the symbol $(a_\epsilon(x,\xi))_\epsilon=(b_\epsilon(x,x,\xi))_\epsilon$ belongs to $\mathcal{S}^{m}_{\Lambda,\rho}$, and as a consequence of Proposition \ref{simmetria}, we have that $(b_\epsilon(y,x,\xi))_\epsilon\in\overline{\mathcal{S}}^{m}_{\Lambda,\rho}$. The same conclusion is true with $\overline{\mathcal{S}}^{m}_{\Lambda,\rho,N}$ and $\mathcal{S}^{m}_{\Lambda,\rho,N}$ (or $\overline{\mathcal{N}}^{m}_{\Lambda,\rho}$ and $\mathcal{N}^{m}_{\Lambda,\rho}$) in place of $\overline{\mathcal{S}}^{m}_{\Lambda,\rho}$ and $\mathcal{S}^{m}_{\Lambda,\rho}$ respectively.
\begin{proposition}
Let $(a_\epsilon)_\epsilon\in{\mathcal{S}}^{m}_{\Lambda,\rho}$. Then for every $\tau\in\mathbb{R}^n$, $b_\epsilon(x,y,\xi)=a_\epsilon((1-\tau)x+\tau y,\xi)$ defines an element of $\overline{\mathcal{S}}^{m}_{\Lambda,\rho}$. The same results hold with  $\mathcal{S}^{m}_{\Lambda,\rho,N}$ and $\overline{\mathcal{S}}^{m}_{\Lambda,\rho,N}$ (or $\mathcal{N}^{m}_{\Lambda,\rho}$ and $\overline{\mathcal{N}}^{m}_{\Lambda,\rho}$) in place of $\mathcal{S}^{m}_{\Lambda,\rho}$ and $\overline{\mathcal{S}}^{m}_{\Lambda,\rho}$ respectively.
\end{proposition}
\ni Choosing $\tau=0$ we have the inclusions ${\mathcal{S}}^{m}_{\Lambda,\rho}\subset{\overline{\mathcal{S}}}^{m}_{\Lambda,\rho}$, ${\mathcal{S}}^{m}_{\Lambda,\rho,N}\subset{\overline{\mathcal{S}}}^{m}_{\Lambda,\rho,N}$ and ${\mathcal{N}}^{m}_{\Lambda,\rho}\subset{\overline{\mathcal{N}}}^{m}_{\Lambda,\rho}$.
\begin{proposition}
The following elementary properties hold:
\begin{itemize}
\item[i)] if $m\le m^{'}$, $\rho\ge\rho^{'}$, $N\le N^{'}$, then  $\overline{\mathcal{S}}^{m}_{\Lambda,\rho,N}\subset\overline{\mathcal{S}}^{m^{'}}_{\Lambda,\rho^{'},N^{'}}$;
\item[ii)] if $(a_\epsilon)_\epsilon\in\overline{\mathcal{S}}^{m}_{\Lambda,\rho,N}$and $(b_\epsilon)_\epsilon\in\overline{\mathcal{S}}^{m^{'}}_{\Lambda,\rho,N^{'}}$ then $(a_\epsilon b_\epsilon)_\epsilon\in\overline{\mathcal{S}}^{m+m^{'}}_{\Lambda,\rho,N+N^{'}}$ and $(a_\epsilon +b_\epsilon)_\epsilon\in\overline{\mathcal{S}}^{\max(m,m^{'})}_{\Lambda,\rho,\max(N,N^{'})}$;
\item[iii)] if $(a_\epsilon)_\epsilon\in\overline{\mathcal{S}}^{m}_{\Lambda,\rho,N}$ then for all $\alpha,\beta,\gamma\in\mathbb{N}^{n}$,  $(\partial^\alpha_\xi\partial^\beta_x\partial^\gamma_y a_\epsilon)_\epsilon\in \overline{\mathcal{S}}^{m-\rho|\alpha+\beta+\gamma|}_{\Lambda,\rho,N}$;
\item[iv)] if $(a_\epsilon)_\epsilon\in\overline{\mathcal{S}}^{m}_{\Lambda,\rho,N}$ then for all $(x^{'},y^{'},\xi^{'})\in\mathbb{R}^{3n}$ $(T_{(x^{'},y^{'},\xi^{'})} a_\epsilon(x,y,\xi))_\epsilon=(a_\epsilon(x-x^{'},y-y^{'},\xi-\xi^{'}))_\epsilon\in\overline{\mathcal{S}}^{m}_{\Lambda,\rho,N}$.
\end{itemize}
Moreover the same statements are valid for amplitudes in $\overline{\mathcal{S}}^{m}_{\Lambda,\rho}$ and $\overline{\mathcal{N}}^{m}_{\Lambda,\rho}$ respectively.
\end{proposition}
\section{Pseudo-differential operators acting on $\mathcal{G}_{\tau,\mathcal{S}}(\mathbb{R}^n)$}
In the classical theory an integral
\begin{equation}
Au(x)=\int_{\mathbb{R}^{2n}}\hskip-10pt e^{i(x-y)\xi}a(x,y,\xi)u(y)\hskip2pt dy\dslash\xi ,
\end{equation}
where $a$ is an amplitude as in \cite{bog, bogII} and $u\in\mathcal{S}(\mathbb{R}^n)$, defines a continuous map from ${\S}$ to $\S$, which can be extended as a continuous map from $\Sp$ to $\Sp$.\\ Formally we obtain a pseudo-differential operator acting on $\mathcal{G}_{\tau,\mathcal{S}}(\mathbb{R}^n)$, by substitution of $a$ and $u$ with $a_\epsilon$ and $u_\epsilon(y)\widehat{\varphi_\epsilon}(y)$ respectively, where $(a_\epsilon)_\epsilon\in\overline{\mathcal{S}}^{m}_{\Lambda,\rho}$, $(u_\epsilon)_\epsilon\in\mathcal {E}_{\tau}(\mathbb{R}^n)$ and $\varphi$ is a fixed mollifier. We begin with the following propositions.
\begin{proposition}
Let $(a_\epsilon)_\epsilon\in\ampg$, $(u_\epsilon)_\epsilon\in\mathcal {E}_{\tau}(\mathbb{R}^n)$ and $\varphi$ be a mollifier. Then for every $x\in\mathbb{R}^n$
\begin{equation}
(b_{\epsilon}(x,y,\xi))_\epsilon=(e^{ix\xi}a_\epsilon(x,y,\xi)u_\epsilon(y)\widehat{\varphi_\epsilon}(y))_\epsilon\in\mathcal{A}^{\nu}_{0}(\mathbb{R}^n_y\times\mathbb{R}^n_\xi) ,
\end{equation}
where $\nu=m_++m^{'}_+$.
\end{proposition}
\begin{proof}
According to Definition 3.2, we have to estimate $\partial^\alpha_\xi\partial^\beta_y b_\epsilon(x,y,\xi)$. This is a finite sum of terms of the type
\begin{equation}
c_{\alpha^{'},\beta^{'},\gamma}e^{ix\xi}(ix)^{\alpha^{'}}\partial^{\alpha-\alpha^{'}}_\xi\partial^{\beta^{'}}_y a_\epsilon(x,y,\xi)\partial^\gamma_y u_\epsilon(y)\partial_y^{\beta-\beta^{'}-\gamma}\hat{\varphi}(\epsilon y)\epsilon^{|\beta-\beta^{'}-\gamma|} ,
\end{equation}
which can be estimated using the definition of $(a_\epsilon)_\epsilon\in\ampg$ and $(u_\epsilon)_\epsilon\in\mathcal {E}_{\tau}(\mathbb{R}^n)$, by
\begin{equation}
c(\alpha,\beta,a,u,\varphi)\langle x\rangle^{|\alpha|}\Lambda(x,\xi)^{m}\langle x-y\rangle^{m^{'}}\langle y\rangle^{N_u(\beta)}\langle y\rangle^{-N_u(\beta)}\epsilon^{-N_a(\alpha,\beta)}\epsilon^{-2N_u(\beta)} .
\end{equation}
In (5.4) one negative power $\epsilon^{-N_u(\beta)}$ comes from the estimate of $|\partial^\gamma u_\epsilon(y)|$, and the other comes from $|\partial^{\beta-\beta^{'}-\gamma}_y\hat{\varphi}(\epsilon y)|\prec\langle\epsilon y\rangle^{-N_u(\beta)}\prec\epsilon^{-N_u(\beta)}\langle y\rangle^{-N_u(\beta)}$, for every $\epsilon\in(0,1]$. Since 
\begin{equation}
\begin{split}
\langle x\rangle^{|\alpha|}\Lambda(x,\xi)^{m}\langle x-y\rangle^{m^{'}}&\prec\langle x\rangle^{|\alpha|}\langle(x,\xi)\rangle^{m_+}\langle x-y\rangle^{m^{'}_+}\\
&\prec\langle x\rangle^{|\alpha|}\langle x\rangle^{m_+}\langle\xi\rangle^{m_+}\langle x\rangle^{m^{'}_+}\langle y\rangle^{m^{'}_+}\\
&\prec\langle x\rangle^{|\alpha|+m_++m^{'}_+}\langle(y,\xi)\rangle^{m_++m^{'}_+} ,
\end{split}
\end{equation}
we obtain our assertion, i.e.
\begin{equation}
|\partial^\alpha_\xi\partial^\beta_y b_\epsilon(x,y,\xi)|\le c\langle x\rangle^{|\alpha|+m_++m^{'}_+}\langle(y,\xi)\rangle^{m_++m^{'}_+}\epsilon^{-N_a(\alpha,\beta)-2N_u(\beta)} ,
\end{equation}
with $c$ independent of $\epsilon$.
\end{proof}
\ni As a consequence of this result we can define the oscillatory integral
\begin{equation}
\int_{\mathbb{R}^{2n}}\hskip-10pt e^{i(x-y)\xi}a_\epsilon(x,y,\xi)u_\epsilon(y)\widehat{\varphi_\epsilon}(y)\hskip2pt dy\dslash\xi ,
\end{equation}
where $-y\xi$ is the phase function of order 2 and $(b_\epsilon(x,y,\xi))_\epsilon\in\mathcal{A}^{\nu}_{0}(\mathbb{R}^n_y\times\mathbb{R}^n_\xi)$ the amplitude. Moreover we observe that for fixed ${\epsilon}\in(0,1]$, this integral can be interpreted as (5.1), since $u_{{\epsilon}}(y)\widehat{\varphi_{{\epsilon}}}(y)$ belongs to $\S$ and $a_{{\epsilon}}(x,y,\xi)$ is an amplitude as in \cite{bog, bogII}. One easily verifies that for all $\epsilon\in(0,1]$, (5.7) can be written as an iterated integral in $dy$ and $d\xi$. Finally the integral (5.7) defines for every $\epsilon\in(0,1]$, a smooth function on $\mathbb{R}^n$, since it satisfies the assumptions of Proposition 3.1.
\begin{proposition}
Under the previous hypothesis, for all $\beta\in\mathbb{N}^n$
\[
(\partial^\beta_x b_\epsilon(x,y,\xi))_\epsilon=(\partial^\beta_x\{e^{ix\xi}a_\epsilon(x,y,\xi)u_\epsilon(y)\widehat{\varphi_\epsilon}(y)\})_\epsilon\in\mathcal{A}_0^{\nu+|\beta|}(\mathbb{R}^n_y\times\mathbb{R}^n_\xi) 
\]
and for all $\alpha,\gamma\in\mathbb{N}^n$, there exists $M\in\mathbb{N}$ such that for every $r>0$
\[
\sup_{(y,\xi)\in\mathbb{R}^{2n},|x|\le r,\epsilon\in(0,1]}\epsilon^{M}\langle(y,\xi)\rangle^{-(\nu+|\beta|)}|\partial^\alpha_\xi\partial^\gamma_y(\partial^\beta_x b_\epsilon(x,y,\xi))|<\infty .
\]
\end{proposition}
\begin{proof}
At first we consider for arbitrary $\alpha,\gamma\in\mathbb{N}^n$, $\partial^\alpha_\xi\partial^\gamma_y(\partial^\beta_x b_\epsilon(x,y,\xi))$. Using the Leibniz rule, we get a finite sum of terms of the type
\begin{equation}
c_{\alpha^{'},\beta^{'},\delta,\gamma}e^{ix\xi}(ix)^{\alpha^{'}}\partial^\delta_\xi(i\xi)^{\beta^{'}}\partial^{\alpha-\alpha^{'}-\delta}_\xi\partial^{\beta-\beta^{'}}_x\partial^{\gamma^{'}}_y a_\epsilon(x,y,\xi)\partial^\sigma u_\epsilon(y)\partial^{\gamma-\gamma^{'}-\sigma}\hat{\varphi}(\epsilon y)\epsilon^{|\gamma-\gamma^{'}-\sigma|} .
\end{equation}
Arguing as before we estimate the absolute value of (5.8) by
\begin{equation}
\begin{split}
&c(\alpha,\beta,\gamma,a,u,\varphi)\langle x\rangle^{|\alpha|+m_++m^{'}_+}\langle\xi\rangle^{|\beta|+m_+}\langle y\rangle^{m^{'}_+}\epsilon^{-N_a(\alpha,\beta,\gamma)}\epsilon^{-2N_u(\gamma)}\\
&\le c^{'}(\alpha,\beta,\gamma,a,u,\varphi)\langle x\rangle^{|\alpha|+\nu}\langle(y,\xi)\rangle^{\nu+|\beta|}\epsilon^{-N_a(\alpha,\beta,\gamma)-2N_u(\gamma)} .
\end{split}
\end{equation}
As a consequence we can claim that $(\partial^\beta_x b_\epsilon(x,y,\xi))_\epsilon\in\mathcal{A}_0^{\nu+|\beta|}(\mathbb{R}^n_y\times\mathbb{R}^n_\xi)$. The second part of the proof is obvious from (5.9).
\end{proof}
\ni Proposition 5.2 allows us to claim that 
\begin{equation}
A_\epsilon u_\epsilon(x):=\int_{\mathbb{R}^{2n}}\hskip-10pt e^{i(x-y)\xi}a_\epsilon(x,y,\xi)u_\epsilon(y)\widehat{\varphi_\epsilon}(y)\hskip2pt dy\dslash\xi
\end{equation}
is an element of $\mathcal{E}_M(\mathbb{R}^n)$. More precisely we have the following result.
\begin{theorem}
Let $(a_\epsilon)_\epsilon\in\ampg$. Then
\[
\begin{split}
(u_\epsilon)_\epsilon\in\mathcal {E}_{\tau}(\mathbb{R}^n)\quad &\Rightarrow\quad (A_\epsilon u_\epsilon)_\epsilon\in\mathcal {E}_{\tau}(\mathbb{R}^n),\\[0.3cm]
(u_\epsilon)_\epsilon\in\mathcal{N}_{\mathcal{S}}(\mathbb{R}^n)\quad &\Rightarrow\quad (A_\epsilon u_\epsilon)_\epsilon\in\mathcal{N}_{\mathcal{S}}(\mathbb{R}^n).
\end{split}
\]
Moreover if $(a_\epsilon)_\epsilon\in{\overline{N}}^m_{\Lambda,\rho}$
\[
(u_\epsilon)_\epsilon\in\mathcal {E}_{\tau}(\mathbb{R}^n)\quad \Rightarrow\quad (A_\epsilon u_\epsilon)_\epsilon\in\mathcal {N}_{\mathcal{S}}(\mathbb{R}^n).
\]
\end{theorem}
\begin{proof}
Since for any $\beta\in\mathbb{N}^{n}$, $(\partial^\beta_x b_\epsilon(x,y,\xi))_\epsilon\in\mathcal{A}_0^{\nu+|\beta|}(\mathbb{R}^n_y\times\mathbb{R}^n_\xi)$, Theorem 3.1 and in par\-ti\-cu\-lar Remark 1, allows us to conclude that there exists a positive constant $c>0$ such that for all $x\in\mathbb{R}^n$, and for all $\epsilon\in(0,1]$  
\begin{equation}
|\partial^\beta A_\epsilon u_\epsilon(x)|\le c\Vert\partial^\beta_x b_\epsilon(x,y,\xi)\Vert_{\overline{N}} ,
\end{equation}
with $\overline{N}$ equal to the least integer greater than $\nu+|\beta|+n+2$.
In order to estimate the right-hand side of (5.11), it suffices to return to the proof of Proposition 5.2, and more precisely to (5.9). Assuming $(u_\epsilon)_\epsilon\in\mathcal {E}_{\tau}(\mathbb{R}^n)$, we put in (5.9) the condition $|\alpha+\gamma|\le\overline{N}$. In this way we conclude that there exists a positive constant $c_{\overline{N}}$ such that for all $x\in\mathbb{R}^n$ and for all $\epsilon\in(0,1]$
\begin{equation}
\Vert\partial^\beta_x b_\epsilon(x,y,\xi)\Vert_{\overline{N}}\le c_{\overline{N}}\langle x\rangle^{\overline{N}+\nu}\epsilon^{-M} ,
\end{equation}
where $M\in\mathbb{N}$, depending on $\overline{N}$, can be chosen larger than $\overline{N}+\nu$. In this way $(A_\epsilon u_\epsilon)_\epsilon\in\mathcal {E}_{\tau}(\mathbb{R}^n)$. We consider now the case of $(u_\epsilon)_\epsilon\in\mathcal{N}_{\mathcal{S}}(\mathbb{R}^n)$. First we use integration by parts and the Leibniz rule and observe that 
\[
\begin{split}
&x^\alpha\partial^\beta A_\epsilon u_\epsilon(x)= 
i^{|\alpha|}\hskip-3pt\sum_{\gamma\le\beta}c_{\gamma}\hskip-3pt\int_{\mathbb{R}^{2n}}\hskip-12pt\partial^\alpha_\xi[e^{-iy\xi}(i\xi)^\gamma\partial^{\beta-\gamma}_x a_\epsilon(x,y,\xi)]e^{ix\xi}u_\epsilon(y)\widehat{\varphi_\epsilon}(y)dy\dslash\xi\\
&=\hskip-2pt i^{|\alpha|}\hskip-8pt\sum_{\substack{\gamma\le\beta\\ \delta\le\alpha\\ \sigma\le\alpha-\delta}}\hskip-7pt c_{\gamma,\delta,\sigma}\hskip-5pt\int_{\mathbb{R}^{2n}}\hskip-10pt e^{-iy\xi} e^{ix\xi}(-iy)^\delta\partial^\sigma(i\xi)^\gamma\partial^{\alpha-\delta-\sigma}_\xi\partial^{\beta-\gamma}_x a_\epsilon(x,y,\xi)u_\epsilon(y)\widehat{\varphi_\epsilon}(y)dy\dslash\xi .
\end{split}
\]
Now for every $\gamma\le\beta$, $\delta\le\alpha$, $\sigma\le\alpha-\delta$, 
\[
e^{ix\xi}(-iy)^\delta\partial^\sigma(i\xi)^\gamma\partial^{\alpha-\delta-\sigma}_\xi\partial^{\beta-\gamma}_x a_\epsilon(x,y,\xi)u_\epsilon(y)\widehat{\varphi_\epsilon}(y)\in{\mathcal{A}}^{\nu+|\beta|}_{0}(\mathbb{R}^n_y\times\mathbb{R}^n_{\xi}) .
\]
In fact since $(u_\epsilon)_\epsilon\in\mathcal{N}_{\mathcal{S}}(\mathbb{R}^n)$, for arbitrary $q\in\mathbb{N}$, $\eta,\mu\in\mathbb{N}^n$
\[
\begin{split}
&\biggl|\partial^\eta_\xi\partial^\mu_y[e^{ix\xi}(-iy)^\delta\partial^\sigma(i\xi)^\gamma\partial^{\alpha-\delta-\sigma}_\xi\partial^{\beta-\gamma}_x a_\epsilon(x,y,\xi)u_\epsilon(y)\widehat{\varphi_\epsilon}(y)]\biggr|\\
&\le c(a,u,\varphi,\eta,\mu,q)\langle x\rangle^{|\eta|+\nu}\langle(y,\xi)\rangle^{\nu+|\beta|}\epsilon^{q-N_a(\alpha,\beta,\eta,\mu)} .
\end{split}
\]
Then under the assumptions $|\eta+\mu|\le\overline{N}$, we obtain for any natural number $q$
\begin{equation}
|x^\alpha\partial^\beta A_\epsilon u_\epsilon(x)|\le c\langle x\rangle^{\overline{N}+\nu}\epsilon^{q-M} ,
\end{equation}
where $M$ depends on $\alpha,\beta,\overline{N}$ while the constant $c$ depends on $\alpha,\beta, \overline{N}$ and $q$. In conclusion choosing $q\ge M$ in (5.13) we have that
\begin{equation}
\begin{array}{cc}
\forall q\in\mathbb{N},\ \exists c>0:\ \forall x\in\mathbb{R}^n,\ \forall\epsilon\in(0,1],\\[0.2cm]
|x^\alpha\partial^\beta A_\epsilon u_\epsilon(x)|\le c\langle x\rangle^{\overline{N}+\nu}\epsilon^{q} .
\end{array}
\end{equation}
At this point, since $\overline{N}+\nu$ is independent of $\alpha$, the following statement holds:
\[
\begin{array}{cc}
\forall q\in\mathbb{N},\ \exists c>0:\ \forall x\in\mathbb{R}^n,\ \forall\epsilon\in(0,1],\\[0.2cm]
|x^\alpha\partial^\beta A_\epsilon u_\epsilon(x)|\le c\epsilon^q .
\end{array}
\]
In the case of $(a_\epsilon)_\epsilon\in{\overline{\mathcal{N}}}^{m}_{\Lambda,\rho}$, the same arguments lead, for $(u_\epsilon)_\epsilon\in{\mathcal{E}}_{\tau}(\mathbb{R}^n)$, to $(A_\epsilon u_\epsilon)_\epsilon\in{\mathcal{N}}_{\mathcal{S}}(\mathbb{R}^n)$. This result completes the proof.
\end{proof}
\begin{remark}
From the previous computations we get the additional information that $(A_\epsilon u_\epsilon)_\epsilon$ is a net of Schwartz functions. More precisely:
\[
\begin{array}{cc}
\forall\alpha ,\beta\in\mathbb{N}^n,\ \exists N\in\mathbb{N},\ \exists c>0:\ \forall x\in\mathbb{R}^n,\ \forall\epsilon\in(0,1],\\[0.2cm]
|x^\alpha\partial^\beta A_\epsilon u_\epsilon(x)|\le c\epsilon^{-N}.
\end{array}
\]
\end{remark}
\ni A consequence of Theorem 5.1 is that it enables us to give a natural definition of pseudo-dif\-fe\-ren\-tial operator acting on $\mathcal{G}_{\tau,\mathcal{S}}(\mathbb{R}^n)$.
\begin{definition}
Let $(a_\epsilon)_\epsilon\in\ampg$. We call pseudo-differential operator of amplitude $(a_\epsilon)_\epsilon$ the linear map $A:\mathcal{G}_{\tau,\mathcal{S}}(\mathbb{R}^n)\rightarrow\mathcal{G}_{\tau,\mathcal{S}}(\mathbb{R}^n)$ such that, for $u\in\mathcal{G}_{\tau,\mathcal{S}}(\mathbb{R}^n)$ with representative $(u_\epsilon)_\epsilon$, $Au$ is the generalized function having as representative $(A_\epsilon u_\epsilon)_\epsilon $ defined in (5.10).
\end{definition}
\noindent Note that Definition 5.1 is given with $(a_\epsilon)_\epsilon\in\ampg$, but in the sequel we shall be concerned mainly with regular amplitudes in $\overline{\mathcal{S}}^{m}_{\Lambda,\rho,N}$. Regular amplitudes allow to develop a complete theory, modelled on the classical one, for pseudo-differential operators acting on $\gts$.
\ni
From Theorem 5.1 it follows that if $(a_\epsilon)_\epsilon$ is a negligible amplitude, then $A$ is the operator identically zero. Therefore, if $(a_\epsilon)_\epsilon$ and $(b_\epsilon)_\epsilon$ belong to ${\overline{\mathcal{S}}}^m_{\Lambda,\rho}$ with $(a_\epsilon -b_\epsilon)_\epsilon\in{\overline{\mathcal{N}}}^{m}_{\Lambda,\rho}$, the corresponding pseudo-differential operators $A$ and $B$ coincide.\\
\ni Before introducing the definition of operators with $\mathcal{S}$-regular kernel, we present the following proposition concerning tempered generalized functions defined by integrals.
\begin{proposition}
Let $k=(k_\epsilon)_\epsilon+\mathcal{N}_{\tau}(\mathbb{R}^{2n})\in\mathcal{G}_{\tau}(\mathbb{R}^{2n})$ and $u=(u_\epsilon)_\epsilon+\mathcal{N}_{\tau}(\mathbb{R}^n)\in\mathcal{G}_{\tau}(\mathbb{R}^n)$. We have that:
\begin{itemize}
\item[i)] for every $x\in\mathbb{R}^n$, $k(x,\cdot):=(k_\epsilon(x,y))_\epsilon+\mathcal{N}_\tau(\mathbb{R}^n_y)$ belongs to $\mathcal{G}_\tau(\mathbb{R}^n)$ and the integral $\displaystyle\int_{\mathbb{R}^n}\hskip-4pt k(\cdot,y)u(y)dy$, defines an element $v$ of $\mathcal{G}_{\tau}(\mathbb{R}^n)$ with representative
\begin{equation}
v_\epsilon(x)=\int_{\mathbb{R}^n}\hskip-3pt k_\epsilon(x,y)u_\epsilon(y)\widehat{\varphi_\epsilon}(y)\hskip1pt dy ;
\end{equation}
\item[ii)] if $k\in\mathcal{G}^\infty_{\mathcal{S}}(\mathbb{R}^{2n})$ then $v\in\gss$.
\end{itemize}
The same results hold with $u\in\gts$.
\end{proposition}
\begin{proof}
For fixed $x\in\mathbb{R}^n$ it is immediate to verify that $(k_\epsilon)_\epsilon\in\mathcal{E}_{\tau}(\mathbb{R}^{2n})$ implies $(k_\epsilon(x,y))_\epsilon\in\mathcal{E}_{\tau}(\mathbb{R}^{n}_y)$ and in the same way $(k_\epsilon)_\epsilon\in\mathcal{N}_{\tau}(\mathbb{R}^{2n})$ implies $(k_\epsilon(x,y))_\epsilon\in\mathcal{N}_{\tau}(\mathbb{R}^n_y)$. As a consequence we have that $k(x,\cdot)$ is a generalized function in $\mathcal{G}_{\tau}(\mathbb{R}^n)$.\\
The integral in i) defines for every $x$ a generalized complex number $v(x)$. Our aim is to prove that:
\begin{itemize} 
\item[1)] $(k_\epsilon)_\epsilon\in\mathcal{E}_{\tau}(\mathbb{R}^{2n})$ and $ (u_\epsilon)_\epsilon\in\mathcal {E}_{\tau}(\mathbb{R}^n)$ implies $(v_\epsilon)_\epsilon\in\mathcal {E}_{\tau}(\mathbb{R}^n)$;
\item[2)] $(k_\epsilon)_\epsilon\in\mathcal{N}_{\tau}(\mathbb{R}^{2n})$ and $(u_\epsilon)_\epsilon\in\mathcal {E}_{\tau}(\mathbb{R}^n)$ implies  $(v_\epsilon)_\epsilon\in\mathcal{N}_{\tau}(\mathbb{R}^n)$;
\item[3)] $(k_\epsilon)_\epsilon\in\mathcal{E}_{\tau}(\mathbb{R}^{2n})$ and $(u_\epsilon)_\epsilon\in\mathcal{N}_{\tau}(\mathbb{R}^{n})$ implies  $(v_\epsilon)_\epsilon\in\mathcal{N}_{\tau}(\mathbb{R}^n)$.
\end{itemize}
At first
\begin{equation}
\partial^\alpha v_\epsilon(x)=\int_{\mathbb{R}^n}\hskip-3pt\partial^\alpha_x k_\epsilon(x,y)u_\epsilon(y)\widehat{\varphi_\epsilon}(y)\hskip1pt dy .
\end{equation}
In case 1, we have estimates of the type:
\begin{equation}
\begin{split}
|\partial^\alpha v_\epsilon(x)|&\le c(k,u,\alpha)\int_{\mathbb{R}^{n}}\hskip-3pt\langle(x,y)\rangle^{N_k(\alpha)}\langle y\rangle^{N_u}|\hat{\varphi}(\epsilon y)|dy\hskip1pt\epsilon^{-N_k(\alpha)-N_u}\\
&\le c^{'}(k,u,\alpha,\varphi)\langle x\rangle^{N_k(\alpha)}\int_{\mathbb{R}^{n}}\hskip-2pt\langle y\rangle^{-n-1}dy\hskip1pt\epsilon^{-2N_k(\alpha)-2N_u-n-1} .
\end{split}
\end{equation}
Choosing $M=2N_k(\alpha)+2N_u+n+1$, we obtain the result $|\partial^\alpha v_\epsilon(x)|\le c^{''}\langle x\rangle^M\epsilon^{-M}$.\\
Case 2: for any $q$ we get
\begin{equation}
\begin{split}
|\partial^\alpha v_\epsilon(x)|&\le c(k,u,\alpha,q)\int_{\mathbb{R}^{n}}\langle(x,y)\rangle^{N_k(\alpha)}\langle y\rangle^{N_u}|\hat{\varphi}(\epsilon y)|dy\hskip1pt\epsilon^{q-N_u}\\
&\le c^{'}(k,u,\alpha,q,\varphi)\langle x\rangle^{N_k(\alpha)}\int_{\mathbb{R}^{n}}\langle y\rangle^{-n-1}dy\hskip1pt\epsilon^{q-N_k(\alpha)-2N_u-n-1} .
\end{split}
\end{equation}
We omit the proof of the third case since it is analogous to the previous one.
It remains to consider the possibility that $k$ is an ${\mathcal{S}}$-regular ge\-ne\-ra\-li\-zed function. This means to prove that, replacing $\mathcal{E}_{\tau}(\mathbb{R}^{2n})$ and $\mathcal{N}_{\tau}(\mathbb{R}^{2n})$ by 
$\mathcal{E}^\infty_{\mathcal{S}}(\mathbb{R}^{2n})$ and $\mathcal{N}_{\mathcal{S}}(\mathbb{R}^{2n})$ respectively in 1), 2) and 3), we obtain $(v_\epsilon)_\epsilon\in\mathcal{E}^\infty_{\mathcal{S}}(\mathbb{R}^n)$ in 1) and $(v_\epsilon)_\epsilon\in\mathcal{N}_{\mathcal{S}}(\mathbb{R}^n)$ in 2) and 3).
Therefore, in case 1 we assume
 \[
\begin{split}
|x^\alpha\partial^{\beta}_x k_\epsilon(x,y)|&\le c_1\langle y\rangle^{-N_u-n-1}\epsilon^{-N_k} ,\\
|u_\epsilon(y)|&\le c_2\langle y\rangle^{N_u}\epsilon^{-N_u} ,
\end{split}
\]
and we get $|x^\alpha \partial^\beta v_\epsilon(x)|\le (k,u,\alpha,\beta,\varphi)\epsilon^{-N_k-N_u}$.  
In case 2, the estimates
\[
\begin{split}
|x^\alpha\partial^{\beta}_x k_\epsilon(x,y)|&\le c_1\langle y\rangle^{-N_u-n-1}{\epsilon}^q ,\\
|u_\epsilon(y)|&\le c_2\langle y\rangle^{N_u}\epsilon^{-N_u}
\end{split}
\]
lead, for arbitrary $q$, to $|x^\alpha \partial^\beta v_\epsilon(x)|\le c(k,u,\alpha,\beta,\varphi,q)\epsilon^{q-N_u}$ and we can argue in an analogous way for the third case.
These results complete the proof. In fact since $\mathcal{N}_{\mathcal{S}}(\mathbb{R}^n)\subset\mathcal{N}_{\tau}(\mathbb{R}^n)$ we already proved assertions $i)$ and $ii)$ with $u\in\gts$.
\end{proof}
\begin{definition}
A linear map $A:\mathcal{G}_{\tau,\mathcal{S}}(\mathbb{R}^n)\rightarrow\mathcal{G}_{\tau,\mathcal{S}}(\mathbb{R}^n)$ is called operator with $\mathcal{S}$-regular kernel iff there exists $k_A\in\mathcal{G}^\infty_{\mathcal{S}}(\mathbb{R}^{2n})$ such that for all $u$
\begin{equation}
Au=\int_{\mathbb{R}^n}\hskip-7pt k_{A}(\cdot,y)u(y)dy .
\end{equation}
\end{definition}
\begin{proposition}
Any operator with $\mathcal{S}$-regular kernel maps $\mathcal{G}_{\tau,\mathcal{S}}(\mathbb{R}^n)$ into $\gss$.
\end{proposition}
\begin{proof}
This is a simple consequence of Proposition 5.3.
\end{proof}
\ni We can consider the mapping property $\gts\to\gss$ as the definition of regularizing operators. Therefore, every operator with $\mathcal{S}$-regular kernel is regularizing. There exists an interesting characterization of operators with $\mathcal{S}$-regular kernel.
\begin{proposition}
$A$ is an operator with $\mathcal{S}$-regular kernel iff it is a pseudo-differential operator with smoo\-thing symbol.
\end{proposition}
\begin{proof}
From the definition of $A$ there exists $k_A\in\mathcal{G}^{\infty}_{\mathcal{S}}(\mathbb{R}^{2n})$ satisfying (5.19). At first we prove that, for any representative $(k_{A,\epsilon})_\epsilon$ of $k_A$
\begin{equation}
a_\epsilon(x,\xi)=e^{-ix\xi}\int_{\mathbb{R}^n}e^{iw\xi}k_{A,\epsilon}(x,w)dw
\end{equation}
belongs to $\mathcal{S}^{-\infty}$. From Proposition 4.3, we have to verify that 
\[
\begin{array}{cc}
\exists N\in\mathbb{N}:\ \forall\alpha,\beta\in\mathbb{N}^n,\\[0.2cm]
\displaystyle\sup_{\epsilon\in(0,1]}\epsilon^{N}\Vert z^\alpha\partial^\beta a_\epsilon\Vert_{L^\infty(\mathbb{R}^{2n})}<\infty . 
\end{array}
\]
In fact $(a_\epsilon)_\epsilon\in\mathcal{E}(\mathbb{R}^{2n})$ and, applying the Leibniz rule, $z^\alpha\partial^\beta a_\epsilon(z)=\xi^{\alpha_1} x^{\alpha_2}\partial^{\beta_1}_\xi\partial^{\beta_2}_x a_\epsilon(x,\xi)$ is a finite sum of terms of the type
\[
c_{\gamma_1,\gamma_2,\delta}\xi^{\alpha_1}x^{\alpha_2}e^{-ix\xi}(-i\xi)^{\gamma_2}\partial^\delta(-ix)^{\gamma_1}\hskip-5pt\int_{\mathbb{R}^n}\hskip-7pt e^{iw\xi}(iw)^{\beta_1-\gamma_1}\partial_x^{\beta_2-\gamma_2-\delta}k_{A,\epsilon}(x,w)dw .
\]
After integration by parts, we obtain terms of the kind
\begin{equation}
e^{-ix\xi}\int_{\mathbb{R}^n}e^{iw\xi}x^{\alpha_2}\partial^\delta(-ix)^{\gamma_1}\partial^{\alpha_1+\gamma_2}_w[(iw)^{\beta_1-\gamma_1}\partial^{\beta_2-\gamma_2-\delta}_xk_{A,\epsilon}(x,w)]dw .
\end{equation}
At this point, we easily conclude that the absolute value of $z^\alpha\partial^\beta a_\epsilon(z)$ can be estimated by $c\epsilon^{-N}$, where $c$ is a positive constant independent of $x$, $\xi$ and $\epsilon$, and $N$ appears in the definition of $k_{A,\epsilon}$. Now, from (5.19), for arbitrary $u\in\mathcal{G}_{\tau,\mathcal{S}}(\mathbb{R}^n)$, $Au$ has the following representative
\begin{equation}
\int_{\mathbb{R}^n}\hskip-4pt k_{A,\epsilon}(x,y)u_\epsilon(y)\widehat{\varphi_\epsilon}(y)\hskip1pt dy .
\end{equation}
By Fourier transform and anti-transform in $\S$
\begin{equation}
k_{A,\epsilon}(x,y)=\hskip-4pt\int_{\mathbb{R}^n}\hskip-5pt e^{i(x-y)\xi}\int_{\mathbb{R}^n}\hskip-5pt e^{-iw\xi}k_{A,\epsilon}(x,x-w)dw\hskip1pt\dslash\xi=\int_{\mathbb{R}^n}\hskip-5pt e^{i(x-y)\xi}a_\epsilon(x,\xi)\dslash\xi
\end{equation}
and, changing order in integration
\[
\int_{\mathbb{R}^n}\int_{\mathbb{R}^n}\hskip-5pt e^{i(x-y)\xi}a_\epsilon(x,\xi)\dslash\xi\hskip1pt u_\epsilon(y)\widehat{\varphi_\epsilon}(y)dy=\int_{\mathbb{R}^{2n}}\hskip-5pt e^{i(x-y)\xi}a_\epsilon(x,\xi)u_\epsilon(y)\widehat{\varphi_\epsilon}(y)\hskip1pt dy\dslash\xi .
\]
In conclusion we can claim that $A$ is a pseudo-differential operator with smoothing symbol and the necessity of the condition is shown.\\
For the converse implication assume now that $A$ is a pseudo-differential operator with smoothing symbol $(a_\epsilon)_\epsilon$. We want to prove that
\begin{equation}
k_{A,\epsilon}(x,y)=\int_{\mathbb{R}^n}\hskip-5pt e^{i(x-y)\xi}a_\epsilon(x,\xi)\dslash\xi
\end{equation}
is the representative of an $\mathcal{S}$-regular generalized function $k_A$ and that for any $u\in\mathcal{G}_{\tau,\mathcal{S}}(\mathbb{R}^n)$, $\displaystyle Au=\int_{\mathbb{R}^n}k_{A}(\cdot,y)u(y)dy$.\\
One easily proves that $(k_{A,\epsilon})_\epsilon\in\mathcal{E}(\mathbb{R}^{2n})$ and applying the Leibniz rule and integration by parts, we have
\[
x^\alpha\hskip-2pt y^\beta\hskip-2pt\partial^\gamma_x\partial^\delta_yk_{A,\epsilon}(x,y)=( -i)^{|\beta|}x^\alpha\hskip-5pt\sum_{\gamma^{'}\le\gamma}\hskip-5pt\binom{\gamma}{\gamma^{'}}\hskip-5pt\int_{\mathbb{R}^{n}}\hskip-9pt e^{-iy\xi}\partial^\beta_\xi[e^{ix\xi}(i\xi)^{\gamma^{'}}\hskip-5pt(-i\xi)^\delta\partial_x^{\gamma-\gamma^{'}}\hskip-5pt a_\epsilon(x,\xi)]\dslash\xi .
\]
We can estimate the absolute value of the sum by $c\epsilon^{-N}$, where $c$ is a positive constant independent of $x$, $y$ and $\epsilon$ and $N$ appears in the definition of $(a_\epsilon)_\epsilon\in\mathcal{S}^{-\infty}$. Finally, $Au$ has as representative
\[
\begin{split}
A_\epsilon u_\epsilon(x)&=\int_{\mathbb{R}^{2n}}e^{i(x-y)\xi}a_\epsilon(x,\xi)u_\epsilon(y)\widehat{\varphi_\epsilon}(y)\hskip1pt dy\dslash\xi\\
&=\hskip-3pt\int_{\mathbb{R}^n}\int_{\mathbb{R}^n}\hskip-5pt e^{i(x-y)\xi}a_\epsilon(x,\xi)\dslash\xi\hskip1pt u_\epsilon(y)\widehat{\varphi_\epsilon}(y)dy=\hskip-3pt\int_{\mathbb{R}^n}\hskip-5pt k_{A,\epsilon}(x,y)u_\epsilon(y)\widehat{\varphi_\epsilon}(y)dy .
\end{split}
\]
This equality concludes the proof.
\end{proof}
\begin{remark}
Every representative $(k_{A,\epsilon})_\epsilon$ of $k_{A}$ defines a symbol $(a_\epsilon)_\epsilon\in{\mathcal{S}}^{-\infty}$ in (5.20). It is clear from the proof above that the difference between two symbols obtained in this way belongs to ${\mathcal{N}}^{-\infty}$.
\end{remark}
\ni In the sequel we study the mapping properties of pseudo-dif\-fe\-ren\-tial operators with regular symbol in more detail.
\begin{proposition}
Let $A$ be a pseudo-differential operator with regular symbol $(a_\epsilon)_\epsilon$ in $\mathcal{S}_{\Lambda,\rho,N}^{m}$. Then $A$ maps ${\mathcal{G}}_{\tau,\mathcal{S}}(\mathbb{R}^n)$ into ${\mathcal{G}}_{\tau,\mathcal{S}}(\mathbb{R}^n)$ and ${\mathcal{G}}^\infty_{\mathcal{S}}(\mathbb{R}^n)$ into ${\mathcal{G}}^\infty_{\mathcal{S}}(\mathbb{R}^n)$.
\end{proposition}
\begin{proof}
It remains to prove only the second part of the claim.\\
Let $u\in\gss$. Since $\mathcal{F}_\varphi(u)=(\widehat{u_\epsilon})_\epsilon+\mathcal{N}_\mathcal{S}(\mathbb{R}^n)$,  
$\displaystyle\big(\int_{\mathbb{R}^n}\hskip-7pt e^{ix\xi}a_\epsilon(x,\xi)\widehat{u_\epsilon}(\xi)\dslash\xi\big)_\epsilon$ is a representative of $Au$. From Proposition 2.3, $(\widehat{u_\epsilon})_\epsilon\in{\mathcal{E}}^\infty_{\mathcal{S}}(\mathbb{R}^n)$, 
and, using integration by parts, $x^\alpha\partial^\beta A_\epsilon u_\epsilon(x)$ is a finite sum of terms of the type
\begin{equation}
c_{\delta,\gamma,\sigma}\int_{\mathbb{R}^n}\hskip-9pt e^{ix\xi}\partial^\delta(i\xi)^\gamma\partial^\sigma_\xi\partial^{\beta-\gamma}_x a_\epsilon(x,\xi)\partial^{\alpha-\delta-\sigma}_\xi\widehat{u_\epsilon}(\xi)\dslash\xi .
\end{equation}
At this point, using the definition of regular symbol and the properties established previously, we conclude that
\begin{equation}
\begin{split}
|e^{ix\xi}\partial^\delta(i\xi)^\gamma\partial^\sigma_\xi\partial^{\beta-\gamma}_x a_\epsilon(x,\xi)\partial^{\alpha-\delta-\sigma}_\xi\widehat{u_\epsilon}(\xi)|&\le c(a,\alpha,\beta)\langle\xi\rangle^{|\beta|}\Lambda(x,\xi)^{m_+}|\partial^{\alpha-\delta-\sigma}_\xi\widehat{u_\epsilon}(\xi)|\epsilon^{-N}\\
&\le c(a,\alpha,\beta)\langle x\rangle^{m_+}\langle\xi\rangle^{|\beta|+m_+}|\partial^{\alpha-\delta-\sigma}_\xi\widehat{u_\epsilon}(\xi)|\epsilon^{-N}\\
&\le c^{'}(a,u,\alpha,\beta)\langle x\rangle^{m_+}\langle\xi\rangle^{-n-1}\epsilon^{-N-M} ,
\end{split}
\end{equation}
where $M$ appears in the definition of $(\widehat{u_\epsilon})_\epsilon$ and is independent of derivatives. In other words, we can state that there exists $N^{'}\in\mathbb{N}$ such that for all $\alpha,\beta\in\mathbb{N}^{n}$
\vskip-5pt
\begin{equation}
\sup_{\epsilon\in(0,1]}\epsilon^{N^{'}}\Vert\langle x\rangle^{-m_+} x^\alpha\partial^\beta A_\epsilon u_\epsilon\Vert_{L^\infty(\mathbb{R}^n)}<\infty .
\end{equation}
(5.27) implies the existence of a natural number $N^{'}$ such that for all $\beta\in\mathbb{N}^n$ and $p\in\mathbb{N}$, $\displaystyle\sup_{\epsilon\in(0,1]}\epsilon^{N^{'}}\Vert\langle x\rangle^{p}\partial^\beta A_\epsilon u_\epsilon\Vert_{L^\infty(\mathbb{R}^n)}<\infty$ and then we obtain $(A_\epsilon u_\epsilon)_\epsilon\in\mathcal{E}^\infty_{\mathcal{S}}(\mathbb{R}^n)$.
\end{proof}
\ni We conclude this section proving that Definition 5.1 is independent, in the weak sense, of the choice of the mollifier $\varphi$.
\begin{proposition}
Let $\varphi_1$ and $\varphi_2$ be two mollifiers satisfying property (2.1). Then for all $(a_\epsilon)_\epsilon\in\amprN$ and $u\in\gts$
\begin{equation}
\begin{split}
&\hskip2pt \biggl(\int_{\mathbb{R}^{2n}}\hskip-10pt e^{i(x-y)\xi}a_\epsilon(x,y,\xi)u_\epsilon(y)\widehat{\varphi_{1,\epsilon}}(y)\hskip2pt dy\dslash\xi\biggr)_\epsilon+\mathcal{N}_{\mathcal{S}}({\mathbb{R}^n})\\
&=_{g.t.d.}\biggl(\int_{\mathbb{R}^{2n}}\hskip-10pt e^{i(x-y)\xi}a_\epsilon(x,y,\xi)u_\epsilon(y)\widehat{\varphi_{2,\epsilon}}(y)\hskip2pt dy\dslash\xi\biggr)_\epsilon +\mathcal{N}_{\mathcal{S}}(\mathbb{R}^n) .
\end{split}
\end{equation}
\end{proposition}
\begin{proof}
It is sufficient to write for $f\in\mathcal{S}(\mathbb{R}^n)$
\begin{equation}
\begin{split}
&\int_{\mathbb{R}^n}f(x)\int_{\mathbb{R}^{2n}}\hskip-10pt e^{i(x-y)\xi}a_\epsilon(x,y,\xi)u_\epsilon(y)(\widehat{\varphi_{1,\epsilon}}-\widehat{\varphi_{2,\epsilon}})(y)\hskip2pt dy\dslash\xi\ dx\\
&=\int_{\mathbb{R}^n}\int_{\mathbb{R}^{2n}}\hskip-10pt e^{i(x-y)\xi}a_\epsilon(x,y,\xi)f(x)\hskip2pt dx\dslash\xi\ u_\epsilon(y)(\widehat{\varphi_{1,\epsilon}}-\widehat{\varphi_{2,\epsilon}})(y)dy\\
&=\int_{\mathbb{R}^n}\int_{\mathbb{R}^{2n}}\hskip-10pt e^{i(x-y)\xi}a_\epsilon(y,x,-\xi)f(y)\hskip2pt dy\dslash\xi\ u_\epsilon(x)(\widehat{\varphi_{1,\epsilon}}-\widehat{\varphi_{2,\epsilon}})(x)dx .
\end{split}
\end{equation}
Now repeating the same arguments in the proof of Theorem 5.1, we have that for all $x\in\mathbb{R}^n$ and for all $\epsilon\in(0,1]$
\begin{equation}
\biggl|x^\alpha\partial_x^\beta\biggr(\int_{\mathbb{R}^{2n}}\hskip-10pt e^{i(x-y)\xi}a_\epsilon(y,x,-\xi)f(y)\hskip2pt dy\dslash\xi\biggl)\biggl|\le c\epsilon^{-N} .
\end{equation}
By definition of the mollifier for all $\alpha\in\mathbb{N}^n$, $\partial^\alpha(\widehat{\varphi_{1,\epsilon}}-\widehat{\varphi_{2,\epsilon}})(0)=0$. Therefore, for arbitrary $q\in\mathbb{N}$
\begin{equation}
\begin{split}
&\biggl|\int_{\mathbb{R}^n}\int_{\mathbb{R}^{2n}}\hskip-10pt e^{i(x-y)\xi}a_\epsilon(y,x,-\xi)f(y)\hskip2pt dy\dslash\xi\ u_\epsilon(x)(\widehat{\varphi_{1,\epsilon}}-\widehat{\varphi_{2,\epsilon}})(x)dx\biggr|\\
&=\biggl|\int_{\mathbb{R}^n}\int_{\mathbb{R}^{2n}}\hskip-10pt e^{i(x-y)\xi}a_\epsilon(y,x,-\xi)f(y)\hskip2pt dy\dslash\xi\ u_\epsilon(x)\hskip-5pt\sum_{|\alpha|=q+1}\frac{\partial^\alpha(\widehat{\varphi_1}-\widehat{\varphi_2})(\epsilon\theta x)}{\alpha !}(\epsilon x)^\alpha dx\biggr|\\
&\le c\epsilon^{q+1-M} ,
\end{split}
\end{equation}
where $M$ depends on the representative of $u\in\gts$. This estimate completes the proof.
\end{proof}
\ni It is easy to see that the equality in $\gts$ cannot be expected in $(5.28)$. In fact, choosing $(a_\epsilon)_\epsilon$ and $(u_\epsilon)_\epsilon$ identically 1 and taking $\widehat{\varphi_1}(\xi)=1$ for $|\xi|\le 1$, $\widehat{\varphi_1}(\xi)=0$ for $|\xi|\ge 2$ and $\widehat{\varphi_2}(\xi)=1$ for $|\xi|\le 3$, we get $\sup_{\xi\in\mathbb{R}^n}|(\widehat{\varphi_{1,\epsilon}}-\widehat{\varphi_{2,\epsilon}})(\xi)|\ge 1$. Then  $(\widehat{\varphi_{1,\epsilon}})_\epsilon+\mathcal{N}_\mathcal{S}(\mathbb{R}^n)=_{g.t.d.}(\widehat{\varphi_{2,\epsilon}})_\epsilon+\mathcal{N}_\mathcal{S}(\mathbb{R}^n)$ but $(\widehat{\varphi_{1,\epsilon}}-\widehat{\varphi_{2,\epsilon}})_\epsilon\notin\mathcal{N}_\mathcal{S}(\mathbb{R}^n)$.
\section{Alternative definitions}
In this section we propose other possible definitions of pseudo-differential operator, investigating, in some particular cases, the relationships with De\-fi\-ni\-tion 5.1.
\begin{proposition}
Let $(a_\epsilon)_\epsilon\in\ampg$. For $(u_\epsilon)_\epsilon\in\mathcal {E}_{\tau}(\mathbb{R}^n)$ we define
\begin{equation}
{\widetilde{A}}_\epsilon u_\epsilon(x)=\int_{\mathbb{R}^{2n}}\hskip-7pt e^{i(x-y)\xi}a_\epsilon(x,y,\xi)u_\epsilon(y)\widehat{\varphi_\epsilon}(y)\widehat{\varphi_\epsilon}(\xi)\hskip2pt dy\dslash\xi .
\end{equation}
Then 
\[
\begin{split}
(u_\epsilon)_\epsilon\in\mathcal {E}_{\tau}(\mathbb{R}^n)\quad &\Rightarrow\quad ({\widetilde{A}}_\epsilon u_\epsilon)_\epsilon\in\mathcal {E}_{\tau}(\mathbb{R}^n),\\[0.2cm]
(u_\epsilon)_\epsilon\in\mathcal{N}_{\mathcal{S}}(\mathbb{R}^n)\quad &\Rightarrow\quad({\widetilde{A}}_\epsilon u_\epsilon)_\epsilon\in\mathcal{N}_{\mathcal{S}}(\mathbb{R}^n).
\end{split}
\]
\end{proposition}
\begin{proof}
At first we observe that the integral in (6.1) is absolutely convergent. Since for any $\beta\in\mathbb{N}^n$
\[
\partial^\beta{\widetilde{A}}_\epsilon u_\epsilon(x)=\sum_{\gamma\le\beta}\binom{\beta}{\gamma}\int_{\mathbb{R}^{2n}}\hskip-7pt e^{i(x-y)\xi}(i\xi)^{\gamma}\partial^{\beta-\gamma}_xa_\epsilon(x,y,\xi)u_\epsilon(y)\widehat{\varphi_\epsilon}(y)\widehat{\varphi_\epsilon}(\xi)\hskip2pt dy\dslash\xi ,
\]
we can estimate each  $(i\xi)^{\gamma}\partial^{\beta-\gamma}_xa_\epsilon(x,y,\xi)u_\epsilon(y)\widehat{\varphi_\epsilon}(y)\widehat{\varphi_\epsilon}(\xi)$ directly. In this way we obtain for $(u_\epsilon)_\epsilon\in\mathcal {E}_{\tau}(\mathbb{R}^n)$
\begin{equation}
\begin{split}
&|(i\xi)^{\gamma}\partial^{\beta-\gamma}_xa_\epsilon(x,y,\xi)u_\epsilon(y)\widehat{\varphi_\epsilon}(y)\widehat{\varphi_\epsilon}(\xi)|\\
&\le c(\beta,a,u)\langle\xi\rangle^{|\beta|}\Lambda(x,\xi)^{m_+}\langle x-y\rangle^{m^{'}_+}\langle y\rangle^{N_u}|\widehat{\varphi_\epsilon}(y)\widehat{\varphi_\epsilon}(\xi)|\epsilon^{-N_{a}(\beta)-N_u}\\
&\le c(\beta,a,u)\langle\xi\rangle^{|\beta|+m_+}\langle x\rangle^{\nu}\langle y\rangle^{m^{'}_++N_u}|\widehat{\varphi_\epsilon}(y)\widehat{\varphi_\epsilon}(\xi)|\epsilon^{-N_{a}(\beta)-N_u}\\
&\le c^{'}(\beta,a,u,\varphi)\langle x\rangle^{\nu}\langle\xi\rangle^{-n-1}\langle y\rangle^{-n-1}\epsilon^{-N_{a}(\beta)-2N_u-|\beta|-\nu-2n-2} ,
\end{split}
\end{equation}
where, as usual, $N_a(\beta)$ comes from the definition of the amplitude and $N_u$ from the definition of $(u_\epsilon)_\epsilon$. Choosing $M=N_{a}(\beta)+2N_u+|\beta|+\nu+2n+2$ we have for all $x\in\mathbb{R}^n$ $|\partial^\beta\widetilde{A}_\epsilon u_\epsilon(x)|\le c^{''}\langle x\rangle^M\epsilon^{-M}$ or in other words $({\widetilde A}_\epsilon u_\epsilon)_\epsilon\in{\mathcal{E}}_{\tau}(\mathbb{R}^n)$. Let us assume now $(u_\epsilon)_\epsilon\in\mathcal{N}_{\mathcal{S}}(\mathbb{R}^n)$. Our starting point is observing that $x^\alpha\partial^\beta{\widetilde{A}}_\epsilon u_\epsilon(x)$ is a finite sum of terms of this type 
\begin{equation}
\int_{\mathbb{R}^{2n}}\hskip-12pt e^{-iy\xi}e^{ix\xi}(-iy)^\delta\partial^\sigma(i\xi)^\gamma\partial^\eta_\xi\partial^{\beta-\gamma}_x a_\epsilon(x,y,\xi)\partial^{\alpha-\delta-\sigma-\eta}\widehat{\varphi_\epsilon}(\xi)u_\epsilon(y)\widehat{\varphi_\epsilon}(y)\hskip2pt dy\dslash\xi ,
\end{equation}
multiplied for some constants independent of $\epsilon$ and variables. Moreover we can write 
\begin{equation}
\begin{split}
&\biggl|(-iy)^\delta\partial^\sigma(i\xi)^\gamma\partial^\eta_\xi\partial^{\beta-\gamma}_x a_\epsilon(x,y,\xi)\partial^{\alpha-\delta-\sigma-\eta}\widehat{\varphi_\epsilon}(\xi)u_\epsilon(y)\widehat{\varphi_\epsilon}(y)\biggr|\\
&\le c(\alpha,\beta,a)\langle x\rangle^{\nu}\langle y\rangle^{|\alpha|+m^{'}_+}|u_\epsilon(y)||\widehat{\varphi_\epsilon}(y)|\langle\xi\rangle^{|\beta|+m_+}|\partial^{\alpha-\delta-\sigma-\eta}\widehat{\varphi}(\epsilon\xi)|\epsilon^{-N_{a}(\alpha,\beta)}\\
&\le c(\alpha,\beta,a,u,\varphi,q)\langle x\rangle^\nu\langle y\rangle^{-n-1}\langle\xi\rangle^{-n-1}\epsilon^{q}\epsilon^{-N_a(\alpha,\beta)-|\beta|-m_+-n-1},
\end{split}
\end{equation}
where it is important to note the independence of $\nu=m_++m^{'}_+$ of $\alpha$. In conclusion
\begin{equation}
\begin{array}{cc}
\forall\alpha,\beta\in\mathbb{N}^n,\ \forall q\in\mathbb{N},\ \exists c>0:\ \forall x\in\mathbb{R}^n,\ \forall\epsilon\in(0,1],\\[0.2cm]
|x^\alpha\partial^\beta\widetilde{A}_\epsilon u_\epsilon(x)|\le c\langle x\rangle^\nu\epsilon^q
\end{array}
\end{equation}
or in other words $({\widetilde{A}}_\epsilon u_\epsilon)_\epsilon\in\mathcal{N}_{\mathcal{S}}(\mathbb{R}^n)$.
\end{proof}
\begin{remark}
One can easily prove the previous proposition considering the integral (6.1) as an oscillatory integral and applying Theorem 5.1.
\end{remark}
\ni We introduce as a straightforward consequence of Proposition 6.1, the following definition.
\begin{definition}
Let $(a_\epsilon)_\epsilon\in\ampg$. $\widetilde{A}:{\mathcal{G}}_{\tau,\mathcal{S}}(\mathbb{R}^n)\to{\mathcal{G}}_{\tau,\mathcal{S}}(\mathbb{R}^n)$ is the linear operator which maps  $u\in{\mathcal{G}}_{\tau,\mathcal{S}}(\mathbb{R}^n)$, with representative $(u_\epsilon)_\epsilon$, into the generalized function $\widetilde{A}u$ with representative $({\widetilde{A}}_\epsilon u_\epsilon)_\epsilon$ defined in (6.1).
\end{definition}
\begin{remark}
The definition of $\widetilde{A}$ consists in the iterated application of the Co\-lom\-be\-au-Fou\-ri\-er tran\-sform and anti-transform in $\mathcal{G}_{\tau,\mathcal{S}}(\mathbb{R}^n)$. In fact 
\[
\widetilde{A}_\epsilon u_\epsilon(x)=\mathcal{F}^{\ast}_{\varphi,\xi\to x}\mathcal{F}_{\varphi,y\to\xi}\big( a_\epsilon(x,y,\xi)u_\epsilon(y)\big).
\]
We remark that the amplitude $(a_\epsilon)_\epsilon\in\ampg$ and the generalized function $u=(u_\epsilon)_\epsilon+\mathcal{N}_{\mathcal{S}}(\mathbb{R}^n)$ define, for all $(x,\xi)\in\mathbb{R}^{2n}$,  $(a_\epsilon(x,y,\xi)u_\epsilon(y))_\epsilon+\mathcal{N}_{\mathcal{S}}(\mathbb{R}^n_y)\in\mathcal{G}_{\tau,\mathcal{S}}(\mathbb{R}^n_y)$.
\end{remark} 
\noindent Let us analyse now the relationships between $A$ and $\widetilde{A}$, when $(a_\epsilon)_\epsilon$ is a regular amplitude in $\overline{\mathcal{S}}^m_{\Lambda,\rho,N}$.
\begin{proposition}
Let $(a_\epsilon)_\epsilon\in\overline{\mathcal{S}}^m_{\Lambda,\rho,N}$. Then the operators $A$ and $\widetilde{A}$ are equal in the weak sense, i.e. for all $u\in\mathcal{G}_{\tau,\mathcal{S}}(\mathbb{R}^n)$
\begin{equation}
Au=_{g.t.d.}\widetilde{A}u .
\end{equation}
\end{proposition}
\begin{proof}
According to Definition 2.5, we have to prove that for all $f\in\S$
\begin{equation}
\int_{\mathbb{R}^n}\hskip-5pt(Au-\widetilde{A}u)(x)\imath(f)(x)dx=0\quad\quad\text{in}\quad\overline{\mathbb{C}}.
\end{equation}
Since $\biggl(\displaystyle\int_{\mathbb{R}^n}\hskip-5pt(A_\epsilon u_\epsilon(x)-{\widetilde{A}}_\epsilon u_\epsilon(x))f(x)(\widehat{\varphi_\epsilon}(x)-1)dx\biggr)_\epsilon\in\mathcal{N}_o$, we can choose a representative of the integral in (6.7) of the following form
\begin{equation}
\int_{\mathbb{R}^n}\int_{\mathbb{R}^{2n}}\hskip-8pt e^{i(x-y)\xi}a_\epsilon(x,y,\xi)u_\epsilon(y)\widehat{\varphi_\epsilon}(y)(\widehat{\varphi_\epsilon}(\xi)-1) \hskip1pt dy\dslash\xi\ f(x)dx .
\end{equation}
Now changing order in integration we have
\begin{equation}
\int_{\mathbb{R}^n}u_\epsilon(y)\int_{\mathbb{R}^n}e^{-iy\xi}\int_{\mathbb{R}^n}e^{ix\xi}a_\epsilon(x,y,\xi)f(x)dx\ (\widehat{\varphi_\epsilon}(\xi)-1)\dslash\xi\ \widehat{\varphi_\epsilon}(y)dy .
\end{equation}
We study in detail $\displaystyle g_\epsilon(y,\xi):=\int_{\mathbb{R}^n}e^{ix\xi}a_\epsilon(x,y,\xi)f(x)dx$. For arbitrary $\alpha\in\mathbb{N}^n$
\[
(i\xi)^\alpha g_\epsilon(y,\xi)=(-1)^{|\alpha|}\sum_{\beta\le\alpha}\binom{\alpha}{\beta}\int_{\mathbb{R}^n}e^{ix\xi}\partial^\beta_x a_\epsilon(x,y,\xi)\partial^{\alpha-\beta}_x f(x)dx
\]
and then
\begin{equation}
\begin{split}
|(i\xi)^\alpha g_\epsilon(y,\xi)|&\le\sum_{\beta\le\alpha}c(a,\beta)\int_{\mathbb{R}^n}\Lambda(x,\xi)^{m_+}\langle x-y\rangle^{m^{'}_+}|\partial^{\alpha-\beta}f(x)|dx\hskip2pt \epsilon^{-N}\\
&\le \langle \xi\rangle^{m_+}\langle y\rangle^{m^{'}_+}\sum_{\beta\le\alpha}c(a,\beta)\int_{\mathbb{R}^n}\langle x\rangle^\nu|\partial^{\alpha-\beta}f(x)|dx\hskip2pt \epsilon^{-N} .
\end{split}
\end{equation}
As a consequence for all $p\in\mathbb{N}$ there exists a positive constant $c$ such that for all $y\in\mathbb{R}^n$, $\xi\in\mathbb{R}^n$ and $\epsilon\in(0,1]$
\begin{equation}
|g_\epsilon(y,\xi)|\le c\langle\xi\rangle^{-p}\langle y\rangle^{m^{'}_+}\epsilon^{-N} .
\end{equation}
In order to estimate (6.8), we use (6.11) and Taylor's formula applied to $\widehat\varphi$ at 0. We obtain for arbitrary $q\in\mathbb{N}$
\begin{equation}
\begin{split}
&\biggl|\int_{\mathbb{R}^n}\hskip-8pt u_\epsilon(y)\int_{\mathbb{R}^n}\hskip-10pt e^{-iy\xi}g_\epsilon(y,\xi)(\widehat{\varphi_\epsilon}(\xi)-1)\dslash\xi\ \widehat{\varphi_\epsilon}(y)dy\biggr|\\
&=\biggl|\int_{\mathbb{R}^n}\hskip-10pt u_\epsilon(y)\int_{\mathbb{R}^n}\hskip-10pt e^{-iy\xi}g_\epsilon(y,\xi)\hskip-5pt\sum_{\gamma=q+1}\frac{\partial^\gamma\widehat{\varphi}(\epsilon\theta\xi)}{\gamma !}(\epsilon\xi)^\gamma \dslash\xi\ \widehat{\varphi_\epsilon}(y)dy\biggr|\\
&\le c(g,\varphi,q)\int_{\mathbb{R}^n}|u_\epsilon(y)|\langle y\rangle^{m^{'}_+}|\widehat{\varphi_\epsilon}(y)|dy\ \int_{\mathbb{R}^n}\langle\xi\rangle^{-n-1}\dslash\xi\ \epsilon^{q+1-N}\\
&\le c(u,g,\varphi,q)\ \epsilon^{q+1-N-2N_u-m^{'}_+-n-1} .
\end{split}
\end{equation}
Since $N_u$, $N$ and $m^{'}_+$ do not depend on $q$, (6.12) allows us to conclude that  
\begin{equation}
\begin{array}{cc}
\forall q\in\mathbb{N},\ \exists c>0:\ \forall\epsilon\in(0,1],\\[0.2cm]
\biggl|\displaystyle\int_{\mathbb{R}^n}(A_\epsilon u_\epsilon(x)-{\widetilde{A}}_\epsilon u_\epsilon(x))f(x)dx\biggr|\le c\epsilon^{q}.
\end{array}
\end{equation}
\end{proof}
\ni We complete this section studying pseudo-differential operators acting on $\mathcal{G}_{\tau,\mathcal{S}}(\mathbb{R}^n)$ with amplitude in $\mathcal{S}^{m}_{\Lambda,\rho,0}$ independent of $\epsilon$ as in \cite{bog, bogII}. For simplicity we call these amplitudes classical. It is well known that for all $f\in\S$, the oscillatory integral 
\[
Af(x)=\int_{\mathbb{R}^{2n}}\hskip-8pt e^{i(x-y)\xi}a(x,y,\xi)f(y)\hskip2pt dy\dslash\xi 
\]
defines a continuous linear map from $\S$ to $\S$ and it extends to a continuous map from $\mathcal{S}^{'}(\mathbb{R}^n)$ to $\mathcal{S}^{'}(\mathbb{R}^n)$. In detail if $w\in\mathcal{S}^{'}(\mathbb{R}^n)$ then $\langle Aw,f\rangle=\langle w,{\ }^tAf\rangle$, where
\[
{\ }^tAf(y)=\int_{\mathbb{R}^{2n}}\hskip-8pt e^{i(x-y)\xi}a(x,y,\xi)f(x)\hskip2pt dx\dslash\xi .
\]
We want to compare the definition of $A$ on $\S$ and $\mathcal{S}^{'}(\mathbb{R}^n)$ with Definition 5.1 on $\mathcal{G}_{\tau,\mathcal{S}}(\mathbb{R}^n)$, introduced in Section 5 of this paper.
\begin{proposition}
Let $a\in\ampc$ be a classical amplitude. Then for all $u\in\mathcal{G}_{\tau,\mathcal{S}}(\mathbb{R}^n)$ and $f\in\S$
\begin{equation}
\int_{\mathbb{R}^n}\hskip-3pt Au(x)\imath(f)(x)dx=\int_{\mathbb{R}^n}\hskip-3pt u(x)\imath({\ }^tAf)(x)dx .
\end{equation}
\end{proposition}
\begin{proof}
Let us consider $u\in\gts$. We can choose as representative of the left-hand side of (6.14)
\begin{equation}
\begin{split}
&\int_{\mathbb{R}^n}\int_{\mathbb{R}^{2n}}\hskip-8pt e^{i(x-y)\xi}a(x,y,\xi)u_\epsilon(y)\widehat{\varphi_\epsilon}(y)\hskip2pt dy\dslash\xi\hskip2pt f(x)dx\\
&=\int_{\mathbb{R}^n}\hskip-5pt u_\epsilon(y)\int_{\mathbb{R}^{2n}}\hskip-8pt e^{i(x-y)\xi}a(x,y,\xi)f(x)\hskip2pt dx\dslash\xi\hskip2pt \widehat{\varphi_\epsilon}(y)dy=\int_{\mathbb{R}^n}\hskip-5pt u_\epsilon(y){\ }^tAf(y)\widehat{\varphi_\epsilon}(y)dy .
\end{split}
\end{equation}
The last one is a representative of $\displaystyle\int_{\mathbb{R}^n}\hskip-3pt u(x)\imath({\ }^tAf)(x)dx$.
\end{proof}
\begin{corollary}
For all classical amplitudes $a$, the corresponding pseudo-differential operator $A$ in Definition 5.1 maps the factor $\gts/=_{g.t.d.}$ into itself.
\end{corollary}
\begin{corollary}
For all classical amplitudes $a$, the corresponding pseudo-differential o\-pe\-ra\-tor\\ $A:\gts/=_{g.t.d.}\to\gts/=_{g.t.d.}$ is an extension of the classical one defined on $\mathcal{S}^{'}(\mathbb{R}^n)$, i.e. for every $w\in\mathcal{S}^{'}(\mathbb{R}^n)$
\begin{equation}
A(\imath(w))=_{g.t.d.}\imath(Aw) .
\end{equation}
\end{corollary}
\begin{proof}
This proof is simply obtained combining Proposition 6.3 with the equality $$\displaystyle\int_{\mathbb{R}^n}\hskip-3pt \imath(w)(x)\imath(f)(x)dx=\langle w,f\rangle+{\mathcal{N}}_o$$ valid for $w\in\mathcal{S}^{'}(\mathbb{R}^n)$ and $f\in\S$.
\end{proof}
\ni The previous result can be improved on $\S$.
\begin{proposition}
For all classical amplitudes $a$, the corresponding pseudo-differential o\-pe\-ra\-tor $A$ on $\gts$ is an extension of the classical one defined on $\S$, i.e. for every $f\in\S$
\begin{equation}
A(\imath(f))=\imath(Af)\qquad\quad\text{in}\ \ \gts .
\end{equation}
\end{proposition}
\begin{proof}
We want to show that $v_\epsilon(x)\hskip-3pt:=\hskip-3pt\displaystyle\int_{\mathbb{R}^{2n}}\hskip-10pt e^{i(x-y)\xi}a(x,y,\xi)f(y)(\widehat{\varphi_\epsilon}(y)-1)\hskip2pt dy\dslash\xi$, 
defines an element of $\mathcal{N}_{\mathcal{S}}(\mathbb{R}^n)$.
For all $\alpha,\beta\in\mathbb{N}^n$, $x^\alpha\partial^\beta_x v_\epsilon(x)$ is a finite sum of terms of this kind
\begin{equation}
c_{\gamma,\delta,\sigma}\hskip-5pt\int_{\mathbb{R}^{2n}}\hskip-10pt e^{-iy\xi} e^{ix\xi}(-iy)^\delta\partial^\sigma(i\xi)^\gamma\partial^{\alpha-\delta-\sigma}_\xi\partial^{\beta-\gamma}_x a(x,y,\xi)f(y)(\widehat{\varphi_\epsilon}(y)-1)dy\dslash\xi ,
\end{equation}
where for $\gamma\le\beta$, $\delta\le\alpha$, $\sigma\le\alpha-\delta$
\begin{equation}
\big(e^{ix\xi}(-iy)^\delta\partial^\sigma(i\xi)^\gamma\partial^{\alpha-\delta-\sigma}_\xi\partial^{\beta-\gamma}_x a(x,y,\xi)f(y)(\widehat{\varphi_\epsilon}(y)-1)\big)_\epsilon\in\mathcal{A}^{\nu+|\beta|}_0(\mathbb{R}^n_y\times\mathbb{R}^n_\xi) .
\end{equation}
Now fixing $\overline{N}$ as in the proof of Theorem 5.1, we obtain for $|\eta+\mu|\le\overline{N}$ and arbitrary $q\in\mathbb{N}$
\begin{equation}
\begin{split}
&|\partial^\eta_\xi\partial^\mu_y\big(e^{ix\xi}(-iy)^\delta\partial^\sigma(i\xi)^\gamma\partial^{\alpha-\delta-\sigma}_\xi\partial^{\beta-\gamma}_x a(x,y,\xi)f(y)(\widehat{\varphi_\epsilon}(y)-1)\big)|\\
&\le c\langle x\rangle^{\overline{N}+\nu}\langle(y,\xi)\rangle^{\nu+|\beta|}\epsilon^q ,
\end{split}
\end{equation}
where $\epsilon^q$ comes from Taylor's formula applied to $\widehat{\varphi_\epsilon}$ and its derivatives at 0. In conclusion from (6.20) it follows that
\begin{equation}
\begin{array}{cc}
\forall\beta\in\mathbb{N}^n,\ \exists \overline{N}\in\mathbb{N}:\ \forall\alpha\in\mathbb{N}^n,\ \forall q\in\mathbb{N},\ \exists c>0:\ \forall x\in\mathbb{R}^n,\ \forall\epsilon\in(0,1],\\[0.2cm]
|x^\alpha\partial^\beta v_\epsilon(x)|\le c\langle x\rangle^{\overline{N}+\nu}\epsilon^q .
\end{array}
\end{equation}
Since $\overline{N}$ does not depend on $\alpha$, we obtain that $(v_\epsilon)_\epsilon\in\mathcal{N}_{\mathcal{S}}(\mathbb{R}^n)$.
\end{proof}
\section{$\theta$-symbols and product}
In this section we introduce a generalization of Weyl symbols and we study the product of pseudo-differential operators. We follow in the proofs the arguments of Boggiatto, Buzano, Rodino in \cite{bog, bogII}. 
\begin{proposition}
Let $\theta\in\mathbb{R}^n$ and $(a_\epsilon)_\epsilon\in\amprN$. There exists a symbol $(b_{\theta,\epsilon})_\epsilon\in\symrN$ such that for all $(u_\epsilon)_\epsilon\in\mathcal {E}_{\tau}(\mathbb{R}^n)$
\begin{equation}
A_\epsilon u_\epsilon(x)=\int_{\mathbb{R}^{2n}}\hskip-7pt e^{i(x-y)\xi}b_{\theta,\epsilon}((1-\theta)x+\theta y,\xi)u_\epsilon(y)\widehat{\varphi_\epsilon}(y)\hskip2pt dy\dslash\xi .
\end{equation}
In detail 
\begin{equation}
b_{\theta,\epsilon}(x,\xi)=\int_{\mathbb{R}^{2n}}\hskip-7pt e^{-iy\eta}a_\epsilon(x+\theta y,x-(1-\theta)y,\xi-\eta)\hskip2pt dy\dslash\eta
\end{equation}
and it has the following asymptotic expansion 
\begin{equation}
(b_{\theta,\epsilon})_\epsilon\sim\sum_{\beta,\gamma}\frac{(-1)^{|\beta|}}{\beta !\gamma !}\theta^\beta(1-\theta)^\gamma\big(\big(\partial_\xi^{\beta+\gamma}D_x^\beta D_y^\gamma a_\epsilon\big)\vert_{x=y}\big)_\epsilon .
\end{equation}
\end{proposition}
\begin{proof}
We begin by showing that for $\theta\in\mathbb{R}^n$
\begin{equation}
(a_\epsilon(x+\theta y,x-(1-\theta)y,\xi-\eta))_\epsilon\in\mathcal{A}^\nu_{0,N}(\mathbb{R}^n_y\times\mathbb{R}^n_\eta),
\end{equation}
where as usual $\nu=m_++m^{'}_+$.\\
\ni In detail for $\alpha,\beta\in\mathbb{N}^n$, $\partial^\alpha_\eta\partial^\beta_y(a_\epsilon(x+\theta y,x-(1-\theta)y,\xi-\eta))$ is a finite sum of terms of the type
\[
c(\alpha,\gamma)\partial^{\alpha}_\xi\partial^\gamma_x\partial^{\beta-\gamma}_y a_\epsilon(x+\theta y,x-(1-\theta)y,\xi-\eta)\theta^\gamma(1-\theta)^{\beta-\gamma}
\]
and then for all $x,y,\xi,\eta\in\mathbb{R}^n$, for all $\epsilon\in(0,1]$
\[
|\partial^\alpha_\eta\partial^\beta_y a_\epsilon(x+\theta y,x-(1-\theta)y,\xi-\eta)|\le\hskip-2pt c(\alpha,\beta,\theta)\Lambda(x+\theta y,\xi-\eta)^m\langle y\rangle^{m^{'}}\hskip-3pt\epsilon^{-N}.
\]
At this point using the definition of weight function and (4.3), with $z=(x+\theta y,\xi-\eta)$, $\zeta=(x,\xi)$, $s=m_+$, we conclude that
\begin{equation}
|\partial^\alpha_\eta\partial^\beta_y a_\epsilon(x+\theta y,x-(1-\theta)y,\xi-\eta)|\le c(\alpha,\beta,\theta)\Lambda(x,\xi)^{m_+}\langle(y,\eta)\rangle^{\nu}\epsilon^{-N}.
\end{equation}
Hence the integral in (7.2) makes sense as an oscillatory integral. The proof of smoothness of $b_{\theta,\epsilon}(x,\xi)$ for every $\epsilon$, is left to the reader because it is an easy application of Proposition 3.1. We want to prove instead that $(b_{\theta,\epsilon})_\epsilon\in\symrN$.   
Using integration by parts, we write
\begin{equation}
\begin{split}
&\partial^\alpha_\xi\partial^\beta_x b_{\theta,\epsilon}(x,\xi)=\int_{\mathbb{R}^{2n}}\hskip-7pt e^{-iy\eta}a_{\alpha,\beta,\theta,\epsilon}(x,y,\xi,\eta)\hskip2pt dy\dslash\eta\\
&=\int_{\mathbb{R}^{2n}}\hskip-10pt e^{-iy\eta}\langle y\rangle^{-2M_1}(1-\Delta_\eta)^{M_1}\{\langle\eta\rangle^{-2M_2}(1-\Delta_y)^{M_2}a_{\alpha,\beta,\theta,\epsilon}(x,y,\xi,\eta)\}\hskip2pt dy\dslash\eta ,
\end{split}
\end{equation}
where
\[
a_{\alpha,\beta,\theta,\epsilon}(x,y,\xi,\eta)=\sum_{\gamma\le\beta}\binom{\beta}{\gamma}(\partial^\alpha_\xi\partial^{\beta-\gamma}_x\partial^\gamma_y a_\epsilon)(x+\theta y,x-(1-\theta)y,\xi-\eta).
\]
Since $a_{\alpha,\beta,\theta,\epsilon}(x,y,\xi,\eta)$ is estimated by
\begin{equation}
\begin{split}
&c(\alpha,\beta)\Lambda(x+\theta y,\xi-\eta)^{m}\langle y\rangle^{m^{'}}\big(1+\Lambda(x+\theta y,\xi-\eta)\langle y\rangle^{-m^{'}}\big)^{-\rho|\alpha+\beta|}\epsilon^{-N}\\
&\le c(\alpha,\beta)\Lambda(x+\theta y,\xi-\eta)^{m-\rho|\alpha+\beta|}\langle y\rangle^{m^{'}+\rho m^{'}|\alpha+\beta|}\epsilon^{-N},
\end{split}
\end{equation}
we obtain that
\begin{equation}
\begin{split}
&|\langle y\rangle^{-2M_1}(1-\Delta_\eta)^{M_1}\{\langle\eta\rangle^{-2M_2}(1-\Delta_y)^{M_2}a_{\alpha,\beta,\theta,\epsilon}(x,y,\xi,\eta)\}|\\
&\le c(\alpha,\beta,M_1,M_2)\langle y\rangle^{m^{'}+\rho m^{'}|\alpha+\beta|-2M_1}\langle \eta\rangle^{-2M_2}\Lambda(x+\theta y,\xi-\eta)^{m-\rho|\alpha+\beta|}\epsilon^{-N}.
\end{split}
\end{equation}
Under the assumptions $2M_1\ge m^{'}+\rho m^{'}|\alpha+\beta|+|m|+\rho|\alpha+\beta|+n+1$, $2M_2\ge |m|+\rho|\alpha+\beta|+n+1$, and, choosing $\zeta=(x,\xi)$, $z=(x+\theta y,\xi-\eta)$, $s=m-\rho|\alpha+\beta|$ in (4.3), we conclude that there exists a constant $c>0$ such that for all $x,\xi\in\mathbb{R}^n$, for all $\epsilon\in(0,1]$
\[
|\partial^\alpha_\xi\partial^\beta_x b_{\theta,\epsilon}(x,\xi)|\le c\Lambda(x,\xi)^{m-\rho|\alpha+\beta|}\epsilon^{-N}.
\]
We omit to prove the equality
\begin{equation}
A_\epsilon u_\epsilon(x)=\int_{\mathbb{R}^{2n}}\hskip-10pt e^{i(x-y)\xi}b_{\theta,\epsilon}((1-\theta)x+\theta y,\xi)u_\epsilon(y)\widehat{\varphi_\epsilon}(y)\hskip2pt dy\dslash\xi ,
\end{equation}
because for fixed $\epsilon$ it suffices to repeat the classical proofs in \cite{bog} p.42 and \cite{bogII} p.14-15. Finally it remains to verify (7.3).\\
\ni At first we observe that $((\partial^{\beta+\gamma}_\xi D_x^\beta D_y^\gamma a_\epsilon(x,y,\xi))\vert_{x=y})_\epsilon$ belongs to ${\mathcal{S}}^{m-2\rho|\beta+\gamma|}_{\Lambda,\rho,N}$. In this way the formal series in (7.3) is an element of $F{\mathcal{S}}^{m}_{\Lambda,\rho,N}$ with $m_j=m-2\rho j$ and $N_j=N$. Now we consider $b_{\theta,\epsilon}(x,\xi)$ and we expand $a_\epsilon(x+\theta y,x-(1-\theta)y,\xi-\eta)$ with respect to $y$ at $y=0$. We have for $M\ge 1$
\begin{equation}
\begin{split}
&a_\epsilon(x+\theta y,x-(1-\theta)y,\xi-\eta)\\
&=\sum_{|\beta+\gamma|<M}\frac{(-1)^{|\gamma|}}{\beta !\gamma !}\theta^\beta(1-\theta)^\gamma\partial^\beta_x\partial^\gamma_y a_\epsilon(x,x,\xi-\eta) y^{\beta+\gamma}\hskip-2pt + \hskip-2pt r_{M,\epsilon}(x,y,\xi,\eta),
\end{split}
\end{equation}
where $r_{M,\epsilon}(x,y,\xi,\eta)$ is the following sum
\begin{equation}
\sum_{|\beta+\gamma|=M}\hskip-7pt\frac{(-1)^{|\gamma|}}{\beta !\gamma !}\theta^\beta (1-\theta)^\gamma\hskip-4pt\int_{0}^{1}\hskip-3pt(1-t)^{M-1}\partial^\beta_x\partial^\gamma_y a_\epsilon(x+\theta ty,x-(1-\theta)ty,\xi-\eta)dt\ y^{\beta+\gamma}.
\end{equation}
Applying integration by parts and Proposition 3.4, we conclude that
\begin{equation}
b_{\theta,\epsilon}(x,\xi)-\hskip-10pt\sum_{|\beta+\gamma|<M}\hskip-7pt\frac{(-1)^{|\beta|}}{\beta !\gamma !}\theta^\beta(1-\theta)^\gamma(\partial^{\beta+\gamma}_\xi D^\beta_xD^\gamma_y a_\epsilon)(x,x,\xi)
=\int_{\mathbb{R}^{2n}}e^{-iy\eta} r_{M,\epsilon}(x,y,\xi,\eta)\hskip2pt dy\dslash\eta .
\end{equation}
In order to complete the proof we prove that the last integral in (7.12) defines an element of ${\mathcal{S}}^{m-2\rho M}_{\Lambda,\rho,N}$. The crucial point is to observe that
\begin{equation}
\begin{split}
&\int_{\mathbb{R}^{2n}}e^{-iy\eta}y^{\beta+\gamma}\int_{0}^{1}\hskip-3pt(1-t)^{M-1}\partial^\beta_x\partial^\gamma_y a_\epsilon(x+\theta ty,x-(1-\theta)ty,\xi-\eta)dt\hskip2pt dy\dslash\eta\\
&=(-1)^{|\beta+\gamma|}\int_{0}^{1}(1-t)^{M-1}\hskip-3pt\int_{\mathbb{R}^{2n}}\hskip-10pt e^{-iy\eta}\partial^{\beta+\gamma}_\xi D^\beta_x D^\gamma_y a_\epsilon(x+\theta ty,x-(1-\theta)ty,\xi-\eta)\hskip2pt dy\dslash\eta\hskip1pt dt .
\end{split}
\end{equation}
Since $|\beta+\gamma|=M$ repeating previous arguments we conclude that
\[
s_{\beta,\gamma,t,\theta,\epsilon}(x,\xi)=\int_{\mathbb{R}^{2n}}\hskip-5pt e^{-iy\eta}\partial^{\beta+\gamma}_\xi D^\beta_x D^\gamma_y a_\epsilon(x+\theta ty,x-(1-\theta)ty,\xi-\eta)dy\dslash\eta 
\]
belongs to ${\mathcal{S}}^{m-2\rho M}_{\Lambda,\rho,N}$ with uniform estimates with respect to $t\in[0,1]$. This allows us to claim that $\biggl(\displaystyle\int_{\mathbb{R}^{2n}}\hskip-10pt e^{-iy\eta}r_{M,\epsilon}(x,y,\xi,\eta)dy\dslash\eta\biggr)_\epsilon$ is a symbol in ${\mathcal{S}}^{m-2\rho M}_{\Lambda,\rho,N}$.
\end{proof}
\begin{remark}
From the Schwartz kernel theorem, as in the statement of Theorem 4.4 in \cite{bog} and Theorem 5.5 in \cite{bogII}, we can say that given $(a_\epsilon)_\epsilon\in\overline{\mathcal{S}}^{m}_{\Lambda,\rho,N}$, $(b_{\theta,\epsilon})_\epsilon$ is the unique symbol in ${\mathcal{S}}^{m}_{\Lambda,\rho,N}$ such that the equality (7.1) is pointwise valid, for all $(f_\epsilon)_\epsilon\in{{\mathcal{S}}(\mathbb{R}^n)}^{(0,1]}$ in place of $(u_\epsilon\widehat{\varphi_\epsilon})_\epsilon$.
\end{remark}
\ni Choosing $\theta=(0,...,0)$ we can always write a pseudo-differential operator $A$ with regular amplitude $(a_\epsilon(x,y,\xi))_\epsilon\in\amprN$ as a pseudo-differential operator with symbol $(b_{0,\epsilon}(x,\xi))_\epsilon\in\symrN$. As a consequence of Proposition 5.6, we obtain that $A$ maps $\gss$ into $\gss$.\\
Finally, we complete the discussion of $\theta$-symbols, with the generalization to our context, of Theorem 5.1 in \cite{bog}. The easy proof is left to the reader.
\begin{proposition}
If $(b_{\theta_1,\epsilon})_\epsilon$ and $(b_{\theta_2,\epsilon})_\epsilon$ are respectively $\theta_1$ and $\theta_2$-symbols, according to (7.2), of a pseudo-differential operator $A$ with regular amplitude $(a_\epsilon)_\epsilon\in\amprN$, then
\begin{equation}
(b_{\theta_2,\epsilon})_\epsilon\sim\sum_{\alpha}\frac{1}{\alpha !}(\theta_1-\theta_2)^\alpha(\partial^\alpha_\xi D^\alpha_x b_{\theta_1,\epsilon})_\epsilon .
\end{equation}
\end{proposition}
\noindent Now let us consider two pseudo-differential operators $A^{'}$ and $A^{''}$, with regular amplitudes $(a^{'}_{\epsilon})_\epsilon\in\overline{\mathcal{S}}^{m^{'}}_{\Lambda,\rho,N^{'}}$ and $(a^{''}_\epsilon)_\epsilon\in\overline{\mathcal{S}}^{m^{''}}_{\Lambda,\rho,N^{''}}$ respectively. In order to study the composition $A^{'}A^{''}$, we write $A^{'}$ in terms of the 0-symbol $(b^{'}_{0,\epsilon})_\epsilon$ and $A^{''}$ in terms of the 1-symbol $(b^{''}_{1,\epsilon})_\epsilon$. More precisely, since for $(u_\epsilon)_\epsilon\in\mathcal {E}_{\tau}(\mathbb{R}^n)$
\begin{equation}
\begin{split}
A^{'}_\epsilon u_\epsilon(x)&=\int_{\mathbb{R}^{2n}}\hskip-10pt e^{i(x-y)\xi}b^{'}_{0,\epsilon}(x,\xi)u_\epsilon(y)\widehat{\varphi_\epsilon}(y)\hskip2pt dy\dslash\xi=\int_{\mathbb{R}^{n}}e^{ix\xi}b^{'}_{0,\epsilon}(x,\xi)\mathcal{F}_\varphi u_\epsilon(\xi)\dslash\xi\\
&=\int_{\mathbb{R}^n}e^{ix\xi}b^{'}_{0,\epsilon}(x,\xi)\big(u_\epsilon\widehat{\varphi_\epsilon}\big)\widehat{\ }(\xi)\dslash\xi ,
\end{split}
\end{equation}
we can state that $A^{'}A^{''}:\gts\to\gts$ is a linear operator, defined on an arbitrary representative $(u_\epsilon)_\epsilon$ of $u\in\gts$ by
\begin{equation}
A^{'}_\epsilon A^{''}_\epsilon u_\epsilon(x)=\int_{\mathbb{R}^n}e^{ix\xi}b^{'}_{0,\epsilon}(x,\xi)\big(A^{''}_\epsilon u_\epsilon\widehat{\varphi_\epsilon}\big)\widehat{\ }(\xi)\dslash\xi ,
\end{equation}
where
\begin{equation}
A^{''}_\epsilon u_\epsilon(x)=\int_{\mathbb{R}^{2n}}\hskip-10pt e^{i(x-y)\xi}b^{''}_{1,\epsilon}(y,\xi)u_\epsilon(y)\widehat{\varphi_\epsilon}(y)\hskip2pt dy\dslash\xi .
\end{equation}
Our aim is to prove that there exists a pseudo-differential operator $(A^{'}A^{''})_1:\gts\to\gts$ such that for all $u\in\gts$, 
\[
A^{'}A^{''}(u)=_{g.t.d.}(A^{'}A^{''})_1(u) .
\]
We begin with the following proposition.
\begin{proposition}
Under the previous assumptions, we define for an arbitrary $(u_\epsilon)_\epsilon$ in $\mathcal {E}_{\tau}(\mathbb{R}^n)$
\begin{equation}
(A^{'}_\epsilon A^{''}_\epsilon)_{1}u_\epsilon(x)=\int_{\mathbb{R}^n}e^{ix\xi}b^{'}_{0,\epsilon}(x,\xi)\widehat{A^{''}_\epsilon u_\epsilon}(\xi)\dslash\xi .
\end{equation}
Then $(u_\epsilon)_\epsilon\in\mathcal {E}_{\tau}(\mathbb{R}^n)$ implies $((A^{'}_\epsilon A^{''}_\epsilon)_{1}u_\epsilon)_\epsilon\in\mathcal {E}_{\tau}(\mathbb{R}^n)$ and $(u_\epsilon)_\epsilon\in\mathcal{N}_{\mathcal{S}}(\mathbb{R}^n)$ implies $((A^{'}_\epsilon A^{''}_\epsilon)_{1}u_\epsilon)_\epsilon\in\mathcal{N}_{\mathcal{S}}(\mathbb{R}^n)$.
\end{proposition}
\noindent The proof will be based on the following preparatory result.
\begin{lemma}
$(u_\epsilon)_\epsilon\in{\mathcal{E}}_\tau(\mathbb{R}^n)$ implies the following statement
\begin{equation}
\begin{array}{cc}
\forall\alpha,\beta\in\mathbb{N}^n,\ \exists N\in\mathbb{N},\ \exists c>0:\ \forall\xi\in\mathbb{R}^n,\ \forall\epsilon\in(0,1],\\[0.2cm]
|\xi^\alpha\partial^\beta \widehat{A^{''}_\epsilon u_\epsilon}(\xi)|\le c\epsilon^{-N}.
\end{array}
\end{equation}
Further, $(u_\epsilon)_\epsilon\in\mathcal{N}_{\mathcal{S}}(\mathbb{R}^n)$ implies $(\widehat{A^{''}_\epsilon u_\epsilon})_\epsilon\in\mathcal{N}_{\mathcal{S}}(\mathbb{R}^n)$.
\end{lemma}
\begin{proof}
Theorem 5.1 and in particular Remark 5, guarantee for all $\alpha,\beta\in\mathbb{N}^n$, the existence of a natural number $N$ such that
\begin{equation}
\sup_{\epsilon\in(0,1]}\epsilon^{N}\Vert x^\alpha\partial^\beta A^{''}_\epsilon u_\epsilon\Vert_{L^\infty(\mathbb{R}^n)}<\infty .
\end{equation}
Now it suffices to write 
\begin{equation}
\xi^\alpha\partial^\beta\widehat{A^{''}_\epsilon u_\epsilon}(\xi) 
=(-i)^{|\alpha|}\sum_{\gamma\le\alpha}\binom{\alpha}{\gamma}\int_{\mathbb{R}^n}\hskip-4pt e^{-ix\xi}\partial^\gamma (-ix)^\beta\partial^{\alpha-\gamma}A^{''}_\epsilon u_\epsilon(x)dx
\end{equation}
and apply (7.20), for obtaining the first part of our claim. If $(u_\epsilon)_\epsilon\in{\mathcal{N}}_{\mathcal{S}}(\mathbb{R}^n)$ then $(A^{''}_\epsilon u_\epsilon)_\epsilon\in{\mathcal{N}}_{\mathcal{S}}(\mathbb{R}^n)$, so the assertion $(\widehat{A^{''}_\epsilon u_\epsilon})_\epsilon\in\mathcal{N}_{\mathcal{S}}(\mathbb{R}^n)$ follows naturally from the definition of the ideal ${\mathcal{N}}_{\mathcal{S}}(\mathbb{R}^n)$.
\end{proof}
\begin{proof}[Proof of Proposition 7.3]
For arbitrary $\alpha\in\mathbb{N}^n$ and $(u_\epsilon)_\epsilon\in\mathcal {E}_{\tau}(\mathbb{R}^n)$
\begin{equation}
\partial^\alpha((A^{'}_\epsilon A^{''}_\epsilon)_1u_\epsilon)(x)=\sum_{\beta\le\alpha}\binom{\alpha}{\beta}\int_{\mathbb{R}^n}\hskip-4pt e^{ix\xi}(i\xi)^\beta\partial^{\alpha-\beta}_x b^{'}_{0,\epsilon}(x,\xi)\widehat{A^{''}_\epsilon u_\epsilon}(\xi)\dslash\xi .
\end{equation}
Now using Lemma 7.1, we have that
\begin{equation}
\begin{split}
|e^{ix\xi}(i\xi)^\beta\partial^{\alpha-\beta}_x b^{'}_{0,\epsilon}(x,\xi)\widehat{A^{''}_\epsilon u_\epsilon}(\xi)|&\le c(\alpha,\beta,b^{'}_0)\langle x\rangle^{m^{'}_+}\langle\xi\rangle^{m^{'}_++|\beta|}|\widehat{A^{''}_\epsilon u_\epsilon}(\xi)|\epsilon^{-N^{'}}\\
&\le c(\alpha,\beta,b^{'}_0,u)\langle x\rangle^{m^{'}_+}\langle \xi\rangle^{-n-1}\epsilon^{-N^{'}-N(\beta)}.
\end{split}
\end{equation}
Therefore, $(u_\epsilon)_\epsilon\in\mathcal {E}_{\tau}(\mathbb{R}^n)$ implies $((A^{'}_\epsilon A^{''}_\epsilon)_1u_\epsilon)_\epsilon\in\mathcal {E}_{\tau}(\mathbb{R}^n)$. Now we assume  $(u_\epsilon)_\epsilon\in\mathcal{N}_{\mathcal{S}}(\mathbb{R}^n)$. $x^\alpha\partial^\beta(A^{'}_\epsilon A^{''}_\epsilon)_1u_\epsilon(x)$ is a finite sum of terms of the type
\begin{equation}
c(\gamma,\delta,\sigma) 
\int_{\mathbb{R}^n}\hskip-8pt e^{ix\xi}\partial^\delta(i\xi)^\gamma\partial^\sigma_\xi\partial^{\beta-\gamma}_x b^{'}_{0,\epsilon}(x,\xi)\partial^{\alpha-\delta-\sigma}_\xi\widehat{A^{''}_\epsilon u_\epsilon}(\xi)\dslash\xi .
\end{equation}
Lemma 7.1 allows us to conclude that for any $q\in\mathbb{N}$
\begin{equation}
\begin{split}
|e^{ix\xi}\partial^\delta(i\xi)^\gamma\partial^\sigma_\xi\partial^{\beta-\gamma}_x b^{'}_{0,\epsilon}(x,\xi)\partial^{\alpha-\delta-\sigma}_\xi\widehat{A^{''}_\epsilon u_\epsilon}(\xi)|&\le c(\alpha,\beta,b_0^{'})\langle x\rangle^{m^{'}_+}\langle \xi\rangle^{|\beta|+m^{'}_+}|\partial^{\alpha-\delta-\sigma}_\xi\widehat{A^{''}_\epsilon u_\epsilon}(\xi)|\epsilon^{-N^{'}}\\
&\le c(\alpha,\beta,b^{'}_0,u,q)\langle x\rangle^{m^{'}_+}\langle\xi\rangle^{-n-1}\epsilon^{q-N^{'}}.
\end{split}
\end{equation}
Summarizing for all $\alpha,\beta\in\mathbb{N}^n$, for all $q\in\mathbb{N}$, there exists a positive constant $c$ such that, for every $x\in\mathbb{R}^n$ , for every $\epsilon\in(0,1]$
\begin{equation}
|x^\alpha\partial^{\beta}(A^{'}_\epsilon A^{''}_\epsilon)_1u_\epsilon(x)|\le c\langle x\rangle^{m^{'}_+}\epsilon^{q-N^{'}}.
\end{equation}
Since $m^{'}_+$ is independent of $\alpha$ we have that $((A^{'}_\epsilon A^{''}_\epsilon)_1u_\epsilon)_\epsilon\in\mathcal{N}_{\mathcal{S}}(\mathbb{R}^n)$. 
\end{proof}
\ni From Proposition 7.3, we can define $(A^{'}A^{''})_1$ as a linear operator acting on $\gts$. In par\-ti\-cu\-lar we easily prove the following result.
\begin{proposition}
$(A^{'}A^{''})_1:\gts\to\gts$ is a pseudo-differential operator with re\-gu\-lar amplitude 
$(b^{'}_{0,\epsilon}(x,\xi)b^{''}_{1,\epsilon}(y,\xi))_\epsilon\in\overline{\mathcal{S}}^{m^{'}+m^{''}}_{\Lambda,\rho,N^{'}+N^{''}}$.  Its $\theta$-symbol $(b_{\theta,\epsilon}(x,\xi))_\epsilon$ has the asymptotic expansion
\begin{equation}
(b_{\theta,\epsilon})_\epsilon\sim\sum_{\substack{\beta,\gamma,\delta,\sigma\\ \delta+\sigma=\beta+\gamma}}\hskip-4pt\frac{(-1)^{|\beta|}(\beta+\gamma)!}{\beta! \gamma !\delta! \sigma !}\theta^\beta(1-\theta)^\gamma \big(\partial^\delta_\xi D^\beta_x b^{'}_{0,\epsilon}\big)_\epsilon\big(\partial^\sigma_\xi D^\gamma_x b^{''}_{1,\epsilon}\big)_\epsilon
\end{equation}
and in particular
\begin{equation}
(b_{0,\epsilon})_\epsilon\sim\sum_{\alpha}\frac{1}{\alpha !}\big(\partial^\alpha_\xi b^{'}_{0,\epsilon})_\epsilon\big(D^\alpha_x b^{''}_{0,\epsilon})_\epsilon .
\end{equation}
\end{proposition}
\begin{proof}
From (7.18)
\[
(A^{'}_\epsilon A^{''}_\epsilon)_1u_\epsilon(x)=\int_{\mathbb{R}^n}\hskip-4pt e^{ix\xi}b^{'}_{0,\epsilon}(x,\xi)\widehat{A^{''}_\epsilon u_\epsilon}(\xi)\dslash\xi .
\]
Since
\begin{equation}
\widehat{A^{''}_\epsilon u_\epsilon}(\xi)=\int_{\mathbb{R}^n}\hskip-4pt e^{-iy\xi}b^{''}_{1,\epsilon}(y,\xi)u_\epsilon(y)\widehat{\varphi_\epsilon}(y)dy ,
\end{equation}
we can write
\begin{equation}
\begin{split}
(A^{'}_\epsilon A^{''}_\epsilon)_1u_\epsilon(x)&=\int_{\mathbb{R}^n}\hskip-4pt e^{ix\xi}b^{'}_{0,\epsilon}(x,\xi)\int_{\mathbb{R}^n}\hskip-4pt e^{-iy\xi}b^{''}_{1,\epsilon}(y,\xi)u_\epsilon(y)\widehat{\varphi_\epsilon}(y)dy\hskip1pt \dslash\xi\\
&=\int_{\mathbb{R}^{2n}}\hskip-10pt e^{i(x-y)\xi}b^{'}_{0,\epsilon}(x,\xi)b^{''}_{1,\epsilon}(y,\xi)u_\epsilon(y)\widehat{\varphi_\epsilon}(y)\hskip2pt dy\dslash\xi .
\end{split}
\end{equation}
We know that $(b^{'}_{0,\epsilon}(x,\xi))_\epsilon\in\overline{\mathcal{S}}^{m^{'}}_{\Lambda,\rho,N^{'}}$, $(b^{''}_{1,\epsilon}(y,\xi))_\epsilon\in\overline{\mathcal{S}}^{m^{''}}_{\Lambda,\rho,N^{''}}$, and then, from Proposition 4.8, point ii), the amplitude $(b^{'}_{0,\epsilon}(x,\xi)b^{''}_{1,\epsilon}(y,\xi))_\epsilon$ belongs to $\overline{\mathcal{S}}^{m^{'}+m^{''}}_{\Lambda,\rho,N^{'}+N^{''}}$. Proposition 7.1 leads to the following formula
\begin{equation}
(b_{\theta,\epsilon}(x,\xi))_\epsilon\sim \sum_{\substack{\beta,\gamma,\delta,\sigma\\ \delta+\sigma=\beta+\gamma}}\hskip-4pt\frac{(-1)^{|\beta|}(\beta+\gamma)!}{\beta! \gamma !\delta! \sigma !}\theta^\beta(1-\theta)^\gamma(\partial^\delta_\xi D^\beta_x b^{'}_{0,\epsilon}(x,\xi))_\epsilon(\partial^\sigma_\xi D^\gamma_x b^{''}_{1,\epsilon}(x,\xi))_\epsilon
\end{equation}
and therefore, putting $\theta=0$ in (7.31), we have
\begin{equation}
(b_{0,\epsilon}(x,\xi))_\epsilon\sim\sum_{\substack{\gamma,\delta,\sigma\\ \delta+\sigma=\gamma}}\frac{1}{\delta! \sigma !}(\partial^\delta_\xi b^{'}_{0,\epsilon}(x,\xi))_\epsilon(\partial^\sigma_\xi D^{\gamma}_x b^{''}_{1,\epsilon}(x,\xi))_\epsilon .
\end{equation}
Choosing $\theta_2=1$ and $\theta_1=0$ in (7.14), from Proposition 7.2 we obtain an asymptotic expansion of $(b^{''}_{1,\epsilon})_\epsilon$ in terms of the derivatives of $(b^{''}_{0,\epsilon})_\epsilon$, which substituted in (7.32) gives us, by reordering as in \cite{shu} p.28, the assertion $(b_{0,\epsilon})_\epsilon\sim \displaystyle\sum_{\delta}\frac{1}{\delta !}(\partial^\delta_\xi b^{'}_{0,\epsilon})_\epsilon(D^\delta_x b^{''}_{0,\epsilon})_\epsilon$.
\end{proof}
\ni We conclude this section with the following result of weak equality.
\begin{theorem}
The linear operators $A^{'}A^{''}$ and $(A^{'}A^{''})_1$ are equal in the weak sense, i.e. for all $u\in\gts$ 
\begin{equation}
A^{'}A^{''}(u)=_{g.t.d.}(A^{'}A^{''})_1u .
\end{equation}
\end{theorem}
\begin{proof}
We want to check that for all $f\in\mathcal{S}(\mathbb{R}^n)$
\begin{equation}
h_\epsilon=\int_{\mathbb{R}^n}\hskip-7pt f(x)\biggl[\int_{\mathbb{R}^n}\hskip-5pt e^{ix\xi}b^{'}_{0,\epsilon}(x,\xi)\widehat{A^{''}_\epsilon u_\epsilon\widehat{\varphi_\epsilon}}(\xi)\dslash\xi-\int_{\mathbb{R}^n}\hskip-5pt e^{ix\xi}b^{'}_{0,\epsilon}(x,\xi)\widehat{A^{''}_\epsilon u_\epsilon}(\xi)\dslash\xi\biggr]dx
\end{equation}
defines an element of $\mathcal{N}_o$. At first we change order in integration; in this way
\begin{equation}
h_\epsilon=\int_{\mathbb{R}^n}(\widehat{A^{''}_\epsilon u_\epsilon\widehat{\varphi_\epsilon}}-\widehat{A^{''}_\epsilon u_\epsilon})(\xi)\int_{\mathbb{R}^n}\hskip-5pt e^{ix\xi}b^{'}_{0,\epsilon}(x,\xi)f(x)dx\hskip3pt \dslash\xi .
\end{equation}
We study $g_{\epsilon}(\xi):=\displaystyle\int_{\mathbb{R}^n}\hskip-5pt e^{ix\xi}b^{'}_{0,\epsilon}(x,\xi)f(x)dx$ in some detail. By arguments similar to the ones used in the proof of Proposition 6.2, we obtain that for all $\alpha ,\beta\in\mathbb{N}^n$, there exists a positive constant $c$ such that 
\begin{equation}
\forall\xi\in\mathbb{R}^n,\ \forall\epsilon\in(0,1],\qquad\qquad\qquad |\xi^\alpha\partial^\beta g_\epsilon(\xi)|\le c\langle \xi\rangle^{m^{'}_+}\epsilon^{-N^{'}}.\qquad 
\end{equation}
where $m^{'}_+$ and $N^{'}$ appear in the  definition of $b^{'}_{0,\epsilon}$ and they are independent of the derivatives. From the properties of the classical Fourier transform on $\S$, and Taylor's formula applied to $\widehat{\varphi}$ at $0$, we can write for arbitrary $q\in\mathbb{N}$
\begin{equation}
\begin{split}
h_\epsilon &=\int_{\mathbb{R}^n}(\widehat{A^{''}_\epsilon u_\epsilon\widehat{\varphi_\epsilon}}-\widehat{A^{''}_\epsilon u_\epsilon})(\xi)g_\epsilon(\xi)\dslash\xi=\int_{\mathbb{R}^n}A^{''}_\epsilon u_\epsilon(y)(\widehat{\varphi_\epsilon}(y)-1)\widehat{g_\epsilon}(y)\dslash y\\
&=\sum_{|\gamma|=q+1}\int_{\mathbb{R}^n}\hskip-5pt A^{''}_\epsilon u_\epsilon(y)\frac{\partial^\gamma\widehat{\varphi}(\epsilon\theta y)}{\gamma !}(\epsilon y)^\gamma \widehat{g_\epsilon}(y)\dslash y .
\end{split}
\end{equation}
In order to estimate $(h_\epsilon)_\epsilon$, we observe that as a consequence of (7.36), $\sup_{\epsilon\in(0,1]}\epsilon^{N^{'}}\Vert y^\alpha\widehat{g_\epsilon}\Vert_{L^\infty(\mathbb{R}^n)}$ is finite for all $\alpha\in\mathbb{N}^n$. At this point, recalling that $(A^{''}_\epsilon u_\epsilon)_\epsilon\in{\mathcal{E}}_\tau(\mathbb{R}^n)$, and in particular Remark 5 following Theorem 5.1, we conclude
\begin{equation}
|h_\epsilon|\le c\int_{\mathbb{R}^n}\langle y\rangle^{q+1}\langle y\rangle^{-n-2-q}dy\ \epsilon^{q+1-N_u-N^{'}},
\end{equation}
where $N_u$ and $N^{'}$ do not depend on $q$ and $c$ does not depend on $\epsilon$. This estimate completes the proof.
\end{proof}
\begin{remark}
From the previous proof it is clear that the weak equality in (7.33) remains valid even if we use different mollifiers in the definition of $A^{'}$ and $A^{''}$.
\end{remark}
\ni We show by means of an example that weak equality in (7.33) cannot be strengthened to equality. Let $a^{'}$ and $a^{''}$ classical amplitudes identically equal to $1$. Then
\begin{equation}
A^{''}_\epsilon u_\epsilon(x)=\int_{\mathbb{R}^{2n}}\hskip-10pt e^{i(x-y)\xi}u_\epsilon(y)\widehat{\varphi_\epsilon}(y)\hskip2pt dy\dslash\xi= u_\epsilon(x)\widehat{\varphi_\epsilon}(x)
\end{equation}
and
\begin{equation}
A^{'}_\epsilon A^{''}_\epsilon u_\epsilon(x)=\int_{\mathbb{R}^n}e^{ix\xi}\widehat{A^{''}_\epsilon u_\epsilon\widehat{\varphi_\epsilon}}(\xi)\dslash\xi= u_\epsilon(x)(\widehat{\varphi_\epsilon}(x))^2 .
\end{equation}
Since, from Proposition 3.4, $\displaystyle b_{0}(x,\xi)=\int_{\mathbb{R}^{2n}}\hskip-7pt e^{-iy\eta}dy\dslash\eta=1$,  
\begin{equation}
(A^{'}_\epsilon A^{''}_\epsilon)_1u_\epsilon(x)=\int_{\mathbb{R}^n}\hskip-5pt e^{ix\xi}\widehat{u_\epsilon\widehat{\varphi_\epsilon}}(\xi)\dslash\xi=u_\epsilon(x)\widehat{\varphi_\epsilon}(x).
\end{equation}
It is easy to prove that, taking $u_\epsilon$ identically equal to $1$,  $(\widehat{\varphi_\epsilon}^2-\widehat{\varphi_\epsilon})_\epsilon\notin \mathcal{N}_{\mathcal{S}}(\mathbb{R}^n)$.
\begin{proposition}
The linear operators $A^{'}A^{''}$ and $(A^{'}A^{''})_1$ coincide on $\gss$.
\end{proposition}
\begin{proof}
From Proposition 5.6, since $(b^{''}_{1,\epsilon})_\epsilon$ is a regular symbol, $(u_\epsilon)_\epsilon\in\mathcal{E}^\infty_{\mathcal{S}}(\mathbb{R}^n)$ implies $(A^{''}_\epsilon u_\epsilon)_\epsilon\in\mathcal{E}^\infty_{\mathcal{S}}(\mathbb{R}^n)$ and as a consequence $(A^{''}_\epsilon u_\epsilon\widehat{\varphi_\epsilon}-A^{''}_\epsilon u_\epsilon)_\epsilon\in\mathcal{N}_{\mathcal{S}}(\mathbb{R}^n)$. This allows us to conclude that 
\[
A^{'}_\epsilon A^{''}_\epsilon u_\epsilon(x)-(A^{'}_\epsilon A^{''}_\epsilon)_1u_\epsilon(x)=\int_{\mathbb{R}^{2n}}e^{i(x-y)\xi}b^{'}_{0,\epsilon}(x,\xi)A^{''}_\epsilon u_\epsilon(y)(\widehat{\varphi_\epsilon}(y)-1)\, dy\dslash\xi
\]
is an element of $\mathcal{N}_\mathcal{S}(\mathbb{R}^n)$.
\end{proof}
\section{Global hypoellipticity and results of regularity}
In this section we consider a special set of regular symbols and their corresponding pseudo-differential operators. It turns out that this kind of symbols allows the construction of a parametrix acting on $\gts$. The following definition is modelled on the classical one presented in \cite{bog, bogII}.
\begin{definition}
A symbol $(a_\epsilon)_\epsilon\in\symrN$ is called hypoelliptic if there exist $l\le m$ and $R>0$ such that the following statements hold:
\begin{itemize}
\item[i)] $\exists c>0$:\ $\forall z=(x,\xi)\in\mathbb{R}^{2n}$, $|z|\ge R$,\ $\forall\epsilon\in(0,1]$
\begin{equation}
|a_\epsilon(z)|\ge c\Lambda(z)^l\epsilon^N ;
\end{equation}
\item[ii)] $\forall\gamma\in\mathbb{N}^{2n}$,\ $\exists c_\gamma>0$:\ $\forall z\in\mathbb{R}^{2n}$, $|z|\ge R$,\ $\forall\epsilon\in(0,1]$
\begin{equation}
{|\partial^\gamma a_\epsilon(z)|}\le c_\gamma{|a_\epsilon(z)|}\Lambda(z)^{-\rho|\gamma|} .
\end{equation}
\end{itemize}
\end{definition}
\noindent If $l=m$, $(a_\epsilon)_\epsilon$ is called an \it{elliptic symbol}\rm .\\
We denote the set of all $(a_\epsilon)_\epsilon$ satisfying Definition 8.1 with $H{\mathcal{S}}^{m,l}_{\Lambda,\rho,N}$, while for the set of elliptic symbols we use the notation $E{\mathcal{S}}^{m}_{\Lambda,\rho,N}$.\\
\ni Observe now, that (8.1) implies $a_\epsilon(z)\neq 0$ for $|z|\ge R$ and $\epsilon\in(0,1]$, so that $a^{-1}_\epsilon(z)$ is well defined in $\mathcal{E}[\{|z|>R\}]$; multiplying by $\psi\in\mathcal{C}^\infty(\mathbb{R}^{2n})$, $\psi(z)=0$ for $|z|\le R$, $\psi(z)=1$ for $|z|\ge 2R$, we get $(\psi(z)a_\epsilon^{-1}(z))_\epsilon\in\mathcal{E}[\mathbb{R}^{2n}]$. In the sequel we denote $(\psi(z)a_\epsilon^{-1}(z))_\epsilon$ by $(p_{0,\epsilon})_\epsilon$.
\begin{proposition}
We have that
\begin{itemize}
\item[i)] if $(a_\epsilon)_\epsilon\in H{\mathcal{S}}^{m,l}_{\Lambda,\rho,N}$ then $(p_{0,\epsilon})_\epsilon\in H{\mathcal{S}}^{-l,-m}_{\Lambda,\rho,N}$;
\item[ii)] if $(a_\epsilon)_\epsilon\in H{\mathcal{S}}^{m,l}_{\Lambda,\rho,N}$ then $(p_{0,\epsilon}\partial^\gamma a_\epsilon)_\epsilon\in\mathcal{S}^{-\rho|\gamma|}_{\Lambda,\rho,0}$ for every $\gamma$;
\item[iii)] if $(a_\epsilon)_\epsilon\in H{\mathcal{S}}^{m,l}_{\Lambda,\rho,N}$, $(b_\epsilon)_\epsilon\in H\mathcal{S}^{m^{'},l^{'}}_{\Lambda,\rho,N^{'}}$ then $(a_\epsilon b_\epsilon)_\epsilon\in H\mathcal{S}^{m+m^{'},l+l^{'}}_{\Lambda,\rho,N+N^{'}}$.
\end{itemize}
\end{proposition}
\begin{proof}
$i)$ From (8.1) we obtain that there exists a constant $c>0$ such that for all $z\in\mathbb{R}^{2n}$ and $\epsilon\in(0,1]$
\begin{equation}
|p_{0,\epsilon}(z)|\le c\Lambda(z)^{-l}\epsilon^{-N},
\end{equation}
while for $|z|\ge 2R$ and $\epsilon\in(0,1]$
\begin{equation}
|p_{0,\epsilon}(z)|\ge c\Lambda(z)^{-m}\epsilon^{N}.
\end{equation}
Now we want to prove the following statement:
\begin{equation}
\begin{array}{cc}
\forall\gamma\in\mathbb{N}^{2n},\ \exists c_\gamma>0:\ \forall z\in\mathbb{R}^{2n},\ |z|\ge R,\  \forall\epsilon\in(0,1],\\[0.2cm]
|\partial^\gamma a^{-1}_\epsilon(z)|\le c_\gamma |a_\epsilon^{-1}(z)|\Lambda(z)^{-\rho|\gamma|}.
\end{array}
\end{equation}
The case $\gamma=0$ is obvious. So we assume that (8.5) is valid for $|\gamma|\le M$ and we want to verify the same assertion for $|\gamma|\le M+1$. At first we differentiate the equation
\begin{equation}
a_\epsilon(z)a_\epsilon^{-1}(z)=1 ,
\end{equation}
for $|z|\ge R$ and $\epsilon\in(0,1]$. We obtain applying the Leibniz rule
\begin{equation}
a_\epsilon(z)\partial^\gamma a_\epsilon^{-1}(z)=-\sum_{\substack{\alpha+\beta=\gamma\\ \beta<\gamma}}\frac{\gamma !}{\alpha !\beta !}\partial^\alpha a_\epsilon(z)\partial^\beta a_\epsilon^{-1}(z) . 
\end{equation}
As a consequence
\begin{equation}
\frac{\partial^\gamma a_\epsilon^{-1}(z)}{a^{-1}_\epsilon(z)}= -\sum_{\substack{\alpha+\beta=\gamma\\ \beta<\gamma}}\frac{\gamma !}{\alpha !\beta !}\frac{\partial^\alpha a_\epsilon(z)}{a_\epsilon(z)}\frac{\partial^\beta a_\epsilon^{-1}(z) }{a^{-1}_\epsilon(z)} .
\end{equation}
We estimate the left-hand side of (8.8) using the hypothesis of hypoellipticity of $(a_\epsilon)_\epsilon$ and the induction hypothesis. In this way (8.5) holds. This result easily implies
\begin{equation}
|\partial^\gamma p_{0,\epsilon}(z)|\le c'_\gamma\Lambda(z)^{-l-\rho|\gamma|}\epsilon^{-N},\qquad\qquad\qquad z\in\mathbb{R}^n,\ \epsilon\in(0,1]
\end{equation}
and 
\begin{equation}
\quad\ |\partial^\gamma p_{0,\epsilon}(z)|\le c_\gamma |p_{0,\epsilon}(z)|\Lambda(z)^{-\rho|\gamma|},\quad\qquad\qquad |z|\ge 2R,\ \epsilon\in(0,1].\quad
\end{equation}
Collecting (8.3), (8.4), (8.9) and (8.10) we conclude that $(p_{0,\epsilon})_\epsilon\in H{\mathcal{S}}^{m,l}_{\Lambda,\rho,N}$.\\
\ni $ii)$ We write
\[
\partial^\delta(p_{0,\epsilon}(z)\partial^\gamma a_\epsilon(z))=\sum_{\alpha+\beta=\delta}\frac{\delta !}{\alpha !\beta !}\partial^\alpha p_{0,\epsilon}(z)\partial^{\beta+\gamma}a_\epsilon(z).
\]
From (8.5) and (8.2) we obtain for all $z\in\mathbb{R}^{2n}$
\begin{equation}
\begin{split}
|\partial^\delta(p_{0,\epsilon}(z)\partial^\gamma a_\epsilon(z))|&\le\sum_{\alpha+\beta=\delta}\frac{\delta !}{\alpha !\beta !}c_\alpha 1_{[R,+\infty)}(|z|)|a^{-1}_{\epsilon}(z)|\Lambda(z)^{-\rho|\alpha|}c_{\beta +\gamma}|a_{\epsilon}(z)|\Lambda(z)^{-\rho|\beta+\gamma|}\\
&\le c_\delta \Lambda(z)^{-\rho|\gamma+\delta|},
\end{split}
\end{equation}
where $1_{[R,+\infty)}$ is the characteristic function of the interval $[R,+\infty)$.\\
$iii)$ The conclusion follows easily from a direct application of (8.1), (8.2) and the Leibniz formula.
\end{proof}
\begin{theorem}
Let $A$ be a pseudo-differential operator with hypoelliptic symbol $(a_\epsilon)_\epsilon\hskip-2pt\in\hskip-2pt H\mathcal{S}^{m,l}_{\Lambda,\rho,N}$. Then there exists a pseudo-differential operator $P$ with symbol $(p_\epsilon)_\epsilon\in\mathcal{S}^{-l}_{\Lambda,\rho,N}$ such that
\begin{equation}
\begin{split}
PA&=_{g.t.d.} I+R_1 ,\\
AP& =_{g.t.d.} I+R_2 ,
\end{split}
\end{equation}
where $R_1$ and $R_2$ are operators with $\mathcal{S}$-regular kernel.
\end{theorem}
\noindent $P$ is called a \it{parametrix}\rm\ of $A$.
\begin{proof}
We will find a suitable symbol $(p_\epsilon)_\epsilon$ by an asymptotic expansion and then applying Theorem 4.1. At first from Proposition 8.1, point i), $(p_{0,\epsilon})_\epsilon\in H\mathcal{S}^{-l,-m}_{\Lambda,\rho,N}$. Following the construction proposed in \cite{wong}, we define for $k\ge 1$
\begin{equation}
p_{k,\epsilon}=-\biggl\{\sum_{\substack{|\gamma|+j=k\\ j<k}}\frac{(-i)^{|\gamma|}}{\gamma !}\partial^\gamma_x a_\epsilon\partial^\gamma_\xi p_{j,\epsilon}\biggr\}p_{0,\epsilon} .
\end{equation}
We prove by induction that $(p_{k,\epsilon})_\epsilon\in\mathcal{S}^{-l-2\rho k}_{\Lambda,\rho,N}$. For $k=1$ we have
\begin{equation}
p_{1,\epsilon}=-\biggl\{\sum_{|\gamma|=1}\frac{(-i)^{|\gamma|}}{\gamma !}\partial^\gamma_x a_\epsilon\partial^\gamma_\xi p_{0,\epsilon}\biggr\} p_{0,\epsilon}
\end{equation}
and since $(\partial^\gamma_x a_\epsilon p_{0,\epsilon})_\epsilon\in\mathcal{S}^{-\rho|\gamma|}_{\Lambda,\rho,0}$ (Prop. 8.1, point ii), $(\partial^\gamma_\xi p_{0,\epsilon})_\epsilon\in\mathcal{S}^{-l-\rho|\gamma|}_{\Lambda,\rho,N}$, we conclude that $(p_{1,\epsilon})_\epsilon\in\mathcal{S}^{-l-2\rho}_{\Lambda,\rho,N}$. Now we assume that for $j\le k$, $(p_{j,\epsilon})_\epsilon\in\mathcal{S}^{-l-2\rho j}_{\Lambda,\rho,N}$. We write
\[
p_{k+1,\epsilon}=-\biggl\{\sum_{\substack{|\gamma|+j=k+1\\ j<k+1}}\frac{(-i)^{|\gamma|}}{\gamma !}\partial^\gamma_x a_\epsilon\partial^\gamma_\xi p_{j,\epsilon}\biggr\}p_{0,\epsilon} .
\]
$(\partial^\gamma_x a_\epsilon p_{0,\epsilon})_\epsilon\in\mathcal{S}^{-\rho|\gamma|}_{\Lambda,\rho,0}$ and, by induction hypothesis, $(\partial^\gamma_\xi p_{j,\epsilon})_\epsilon\in\mathcal{S}^{-l-2\rho j-\rho|\gamma|}_{\Lambda,\rho,N}$. As a consequence $(p_{k+1,\epsilon})_\epsilon\in\mathcal{S}^{-l-2\rho(k+1)}_{\Lambda,\rho,N}$. At this point
\begin{equation}
\sum_{j=0}^{\infty}(p_{j,\epsilon})_\epsilon\in F\mathcal{S}^{-l}_{\Lambda,\rho,N}
\end{equation}
and Theorem 4.1 allows us to find a symbol $(p_\epsilon)_\epsilon\in\mathcal{S}^{-l}_{\Lambda,\rho,N}$ such that $(p_\epsilon)_\epsilon\sim\sum_{j}(p_{j,\epsilon})_\epsilon$.\\
Let $P$ be the pseudo-differential operator with symbol $(p_\epsilon)_\epsilon$. Let us consider the composition $PA$. From Proposition 7.4 and Theorem 7.1
\begin{equation}
PA=_{g.t.d.} (PA)_1 ,
\end{equation}
where $(PA)_1$ is a pseudo-differential operator with symbol $(b_{0,\epsilon})_\epsilon\in\mathcal{S}^{m-l}_{\Lambda,\rho,2N}$ and
\begin{equation}
(b_{0,\epsilon})_\epsilon\sim\sum_{\alpha}\frac{1}{\alpha !}(\partial^\alpha_\xi p_\epsilon )_\epsilon(D^\alpha_x a_\epsilon)_\epsilon .
\end{equation}
We want to show that $(b_{0,\epsilon}-1)_\epsilon\in\mathcal{S}^{-\infty}_{\Lambda,\rho,2N}$. At first from the definition of asymptotic expansion it follows that for all $M\in\mathbb{N}$, $M\neq 0$
\begin{equation}
\big(b_{0,\epsilon}-\sum_{|\alpha|<M}\frac{1}{\alpha !}\partial^\alpha_\xi p_\epsilon D^\alpha_x a_\epsilon\big)_\epsilon\in\mathcal{S}^{m-l-2\rho M}_{\Lambda,\rho,2N} .
\end{equation}
Introducing in (8.18) the asymptotic expansion of $(p_\epsilon)_\epsilon$ we have
\begin{equation}
\begin{split}
&b_{0,\epsilon}-\sum_{|\alpha|<M}\frac{1}{\alpha !}\partial^\alpha_\xi p_\epsilon D^\alpha_x a_\epsilon\\
&=b_{0,\epsilon}-\sum_{|\alpha|<M}\frac{1}{\alpha !}D^\alpha_x a_\epsilon\sum_{j=0}^{M-1}\partial^\alpha_\xi p_{j,\epsilon}-\sum_{|\alpha|<M}\frac{1}{\alpha !}\partial^\alpha_\xi r_{M,\epsilon}D^\alpha_x a_\epsilon .
\end{split}
\end{equation}
Since $(\partial^\alpha_\xi r_{M,\epsilon}D^\alpha_x a_\epsilon)_\epsilon\in\mathcal{S}^{m-l-2\rho(M+|\alpha|)}_{\Lambda,\rho,2N}$, we obtain combining (8.19) and (8.18), that
\begin{equation}
\biggl(b_{0,\epsilon}-\sum_{|\alpha|<M}\frac{1}{\alpha !}D^\alpha_x a_\epsilon\sum_{j=0}^{M-1}\partial^\alpha_\xi p_{j,\epsilon}\biggr)_\epsilon\in\mathcal{S}^{m-l-2\rho M}_{\Lambda,\rho,2N}
\end{equation}
Now we observe that 
\begin{equation}
\begin{split}
\sum_{|\alpha|<M}\frac{1}{\alpha !}D^\alpha_x a_\epsilon\sum_{j=0}^{M-1}\partial^\alpha_\xi p_{j,\epsilon}= p_{0,\epsilon}a_\epsilon &+\sum_{k=1}^{M-1}\biggl\{ p_{k,\epsilon}a_\epsilon+\sum_{\substack{|\alpha|+j=k\\ j<k}}\frac{1}{\alpha !}\partial^\alpha_\xi p_{j,\epsilon}D^\alpha_x a_\epsilon\biggr\}\\
&+\sum_{\substack{|\alpha|+j\ge M\\ |\alpha|<M, j<M}}\frac{1}{\alpha !}\partial^\alpha_\xi p_{j,\epsilon}D^\alpha_x a_\epsilon .
\end{split}
\end{equation}
We recall (8.13) and the equality, $p_{0,\epsilon}a_\epsilon=1$, for $|z|\ge 2R$ and $\epsilon\in(0,1]$. As a consequence, under the assumption $|z|\ge 2R$, we obtain from (8.21)
\begin{equation}
\sum_{|\alpha|<M}\frac{1}{\alpha !}D^\alpha_x a_\epsilon\sum_{j=0}^{M-1}\partial^\alpha_\xi p_{j,\epsilon}= 1+\sum_{\substack{|\alpha|+j\ge M\\ |\alpha|<M, j<M}}\frac{1}{\alpha !}\partial^\alpha_\xi p_{j,\epsilon}D^\alpha_x a_\epsilon ,
\end{equation}
where the sum on the right-hand side of (8.22) is an element of $\mathcal{S}^{m-l-2\rho M}_{\Lambda,\rho,2N}$. Due to the properties of $(p_{0,\epsilon})_\epsilon$ and continuity over compact sets of the functions involved in (8.22), we can omit the hypothesis $|z|\ge 2R$. In conclusion
\begin{equation}
\biggl(\sum_{|\alpha|<M}\frac{1}{\alpha !}D^\alpha_x a_\epsilon\sum_{j=0}^{M-1}\partial^\alpha_\xi p_{j,\epsilon}- 1\biggr)_\epsilon\in\mathcal{S}^{m-l-2\rho M}_{\Lambda,\rho,2N} .
\end{equation}
Combining (8.20) with (8.23), we conclude that for all natural $M\neq 0$, the difference $(b_{0,\epsilon}-1)_\epsilon$ belongs to $\mathcal{S}^{m-l-2\rho M}_{\Lambda,\rho,2N}$ and then  $(b_{0,\epsilon}-1)_\epsilon\in\mathcal{S}^{-\infty}_{\Lambda,\rho,2N}$.\\
\ni In order to complete the proof we observe that the pseudo-differential operator having as symbol 1, is equal in the weak sense to the identity. In fact for all $u\in\gts$, as a consequence of Proposition 6.2, 
\begin{equation}
\begin{split}
&\int_{\mathbb{R}^{2n}}\hskip-10pt e^{i(x-y)\xi}u_\epsilon(y)\widehat{\varphi_\epsilon}(y)\hskip2pt dy\dslash\xi+\mathcal{N}_{\mathcal{S}}(\mathbb{R}^n)\\
&=_{g.t.d.}\int_{\mathbb{R}^{2n}}\hskip-10pt e^{i(x-y)\xi}u_\epsilon(y)\widehat{\varphi_\epsilon}(y)\widehat{\varphi_\epsilon}(\xi)\hskip2pt dy\dslash\xi+\mathcal{N}_{\mathcal{S}}(\mathbb{R}^n)\\
&=_{g.t.d.}\mathcal{F}^{\ast}_\varphi\mathcal{F}_\varphi u =_{g.t.d.}u .
\end{split}
\end{equation}
Therefore, $PA=_{g.t.d.}(PA)_1=_{g.t.d.}I+R_1$, where $R_1$ is an o\-pe\-ra\-tor with $\mathcal{S}$-regular kernel. Analogously we can construct $(q_\epsilon)_\epsilon\in\mathcal{S}^{-l}_{\Lambda,\rho,N}$ and the corresponding pseudo-differential operator $Q$ such that $AQ=_{g.t.d.}I$ modulo some operator with $\mathcal{S}$-regular kernel. Since, as in the classical theory, the difference $P-Q$ has smoothing symbol, we conclude that there exists an operator $R_2$ with $\mathcal{S}$-regular kernel such that $AP=_{g.t.d.}I+R_2$.
\end{proof}
\ni We conclude this section with the typical result of regularity obtained by the existence of a parametrix.
\begin{theorem}
Let $A$ be a pseudo-differential operator with hypoelliptic symbol $(a_\epsilon)_\epsilon\hskip-2pt\in\hskip-2pt H\mathcal{S}^{m,l}_{\Lambda,\rho,N}$. If $u\in\gts$, $v\in\gss$ and 
\begin{equation}
Au=v ,
\end{equation}
then $u$ is equal in the weak sense to a generalized function in $\gss$.
\end{theorem}
\begin{proof}
We consider the parametrix $P$ of $A$. (8.25) implies
\begin{equation}
PA(u)=Pv .
\end{equation}
Since $PA=_{g.t.d.}I+R_1$, where $R_1$ has $\mathcal{S}$-regular kernel
\begin{equation}
(I+R_1)u=_{g.t.d.}Pv
\end{equation}
and then
\begin{equation}
u=_{g.t.d.}-R_1u+Pv ,
\end{equation}
with $-R_1u+Pv\in\gss$. In fact $-R_1u\in\gss$, since $R_1$ is an operator with $\mathcal{S}$-regular kernel, and $Pv\in\gss$, because $v\in\gss$ and $P$ maps $\gss$ into $\gss$. 
\end{proof}
\ni Inspired by \cite{bog, bogII} we can call a linear map A from $\gts$ into $\gts$ satisfying the assertion of  Theorem 8.2, \it{globally hypoelliptic in the weak (or g.t.d.) sense}\rm. 
\begin{proposition}
Let $a$ be a classical symbol belonging to $H\mathcal{S}^{m,l}_{\Lambda,\rho,0}$ and let $A$ be the corresponding pseudo-differential operator. If
\begin{equation}
Au=_{g.t.d.}v ,
\end{equation}
where $u\in\gts$ and $v\in\gss$, then $u$ is equal in the weak sense to a generalized function in $\gss$.
\end{proposition}
\begin{proof}
It is sufficient to observe that if $A$ is defined with a classical symbol then the operator $P$ involved in (8.12) has a classical symbol too. Therefore, from Corollary 6.1 and Theorem 8.1 we obtain that $Pv=_{g.t.d.}PAu=_{g.t.d.}u+R_1u$, where $Pv$ and $R_1u$ belong to $\gss$.
\end{proof}
\ni The following proposition shows the consistency with the classical regularity result mentioned in the introduction.
\begin{proposition}
Let $a$ be a classical symbol belonging to $H\mathcal{S}^{m,l}_{\Lambda,\rho,0}$, let $A$ be the corresponding pseudo-differential operator and $u$ and $v$ tempered distributions. If
\begin{equation}
A(\imath(u))=_{g.t.d.}\imath(v) ,
\end{equation}
where $\imath(v)\in\gss$ then $u\in\S$.
\end{proposition}
\begin{proof}
As in the proof above, we can write
\[
\imath(u)=_{g.t.d.}P(\imath(v))-R_1(\imath(u)),
\]
where $P$ and $R_1$ have classical symbol of order $-l$ and $-\infty$ respectively, and from Theorem 2.1, $v$ belongs to $\S$. Since from Corollary 6.2, $P(\imath(v))-R_1(\imath(u))=_{g.t.d.}\imath(Pv-R_1u)$, we conclude that $\imath(u)$ is equal in the weak sense to a function in $\S$. This means that $u\in\S$.
\end{proof}
\ni Let us finally give some examples of hypoelliptic symbols. The classical hypoelliptic and elliptic symbols introduced in \cite{bog, bogII} can be considered as elements of $H\mathcal{S}^{m,l}_{\Lambda,\rho,0}$ and $E\mathcal{S}^{m}_{\Lambda,\rho,0}$ independent of $\epsilon$. Moreover if $a(z)$ is a hypoelliptic or elliptic symbol as in \cite{bog, bogII}, we can easily construct a symbol satisfying Definition 8.1, writing $a_\epsilon(z)=\epsilon^b a(z)$ with $b\in\mathbb{R}$. In order to obtain examples with increasing generality we first prove the following proposition.
\begin{proposition}
Let $(a_\epsilon)_\epsilon\in H\mathcal{S}^{m,l}_{\Lambda,\rho,N}$ and $(b_\epsilon)_\epsilon\in\mathcal{S}^{m^{'}}_{\Lambda,\rho,N^{'}}$ with $m^{'}<l$. Then the symbol $a_\epsilon(z)+\epsilon^{N+N^{'}} b_\epsilon(z)$ belongs to $H\mathcal{S}^{m,l}_{\Lambda,\rho,N}$.
\end{proposition}
\begin{proof}
Let us first assume $(b_\epsilon)_\epsilon\in\mathcal{S}^{m^{'}}_{\Lambda,\rho,0}$. We begin by observing that since $(\epsilon^N b_\epsilon)_\epsilon\in\mathcal{S}^{m^{'}}_{\Lambda,\rho,0}\subset\mathcal{S}^{m}_{\Lambda,\rho,N}$ then $a_\epsilon(z)+\epsilon^N b_\epsilon(z)$ belongs to $\mathcal{S}^{m}_{\Lambda,\rho,N}$. Before checking the estimate $(8.1)$ and $(8.2)$, we remark that $m^{'}<l$ implies the statement
\[
\forall c>0,\ \exists R>0:\ \forall z\in\mathbb{R}^{2n}, |z|\ge R,\quad\quad\quad\quad \Lambda(z)^{m^{'}}\le c\Lambda(z)^l .\quad\quad\quad\quad\quad\quad\quad\quad\quad
\]
This result combined with the definitions of $(a_\epsilon)_\epsilon\in H\mathcal{S}^{m,l}_{\Lambda,\rho,N}$ and $(b_\epsilon)_\epsilon\in\mathcal{S}^{m^{'}}_{\Lambda,\rho,0}$ allows us to infer the existence of a constant $R^{'}>0$ such that for all $z\in\mathbb{R}^{2n}$ with $|z|\ge R^{'}$ and for all $\epsilon\in(0,1]$
\begin{equation}
\label{stima1}
\begin{split}
|a_\epsilon(z)+\epsilon^N b_\epsilon(z)|\ge |a_\epsilon(z)|-\epsilon^N|b_\epsilon(z)|\ge c_1\epsilon^N\Lambda(z)^l-c_2\epsilon^N\Lambda(z)^{m^{'}}&\ge\epsilon^N(c_1\Lambda(z)^l-\frac{c_1}{2}\Lambda(z)^l)\\
&=\frac{c_1}{2}\epsilon^N\Lambda(z)^l .
\end{split}
\end{equation}
We consider now $\gamma\in\mathbb{N}^{2n}$, $\gamma\neq 0$. Under the assumptions $|z|\ge R^{'}$, $\epsilon\in(0,1]$
\begin{equation}
\label{stima2}
|\partial^\gamma(a_\epsilon(z)+\epsilon^N b_\epsilon(z))|\le c(|a_\epsilon(z)|\Lambda(z)^{-\rho|\gamma|}+\epsilon^N\Lambda(z)^{m^{'}-\rho|\gamma|})=c\Lambda(z)^{-\rho|\gamma|}(|a_\epsilon(z)|+\epsilon^N\Lambda(z)^{m^{'}}).
\end{equation}
From \eqref{stima1} it follows that
\begin{equation}
\label{stima3}
\begin{split}
\frac{|a_\epsilon(z)|+\epsilon^N\Lambda(z)^{m^{'}}}{|a_\epsilon(z)+\epsilon^N b_\epsilon(z)|}\le 1+\frac{|\epsilon^N b_\epsilon(z)|+\epsilon^N\Lambda(z)^{m^{'}}}{|a_\epsilon(z)+\epsilon^N b_\epsilon(z)|}&\le 1+\frac{\epsilon^N(|b_\epsilon(z)|+\Lambda(z)^{m^{'}})}{\epsilon^N\frac{c_1}{2}\Lambda(z)^l}\\
&\le 1+\frac{\epsilon^N(c_2+1)\Lambda(z)^{m^{'}}}{\epsilon^N\frac{c_1}{2}\Lambda(z)^l}\le c^{'},\quad\quad \text{for\ $|z|\ge R^{'}$}
\end{split}
\end{equation}
and therefore, substituing \eqref{stima3} in \eqref{stima2}, we obtain that 
\begin{equation}
|\partial^\gamma(a_\epsilon(z)+\epsilon^N b_\epsilon(z))|\le c^{''}\Lambda(z)^{-\rho|\gamma|}|a_\epsilon(z)+\epsilon^N b_\epsilon(z)|,\quad\quad\text{for\ $|z|\ge R^{'},\ \epsilon\in(0,1]$}.
\end{equation}
The proof for a general $(b_\epsilon)_\epsilon\in\mathcal{S}^{m^{'}}_{\Lambda,\rho,N^{'}}$ with $m^{'}<l$ is now immediate. In fact if $(b_\epsilon)_\epsilon\in\mathcal{S}^{m^{'}}_{\Lambda,\rho,N^{'}}$ then $(\epsilon^{N^{'}}b_\epsilon)_\epsilon\in\mathcal{S}^{m^{'}}_{\Lambda,\rho,0}$ and for the symbol $a_\epsilon(z)+\epsilon^N(\epsilon^{N^{'}}b_\epsilon(z))$ we can repeat the previous arguments.
\end{proof}
\begin{remark}
$H\mathcal{S}^{m,l}_{\Lambda,\rho,N}$ is a subset of $H\mathcal{S}^{m,l}_{\Lambda,\rho,M}$ for $M\ge N$. As a consequence if $(a_\epsilon)_\epsilon\in H\mathcal{S}^{m,l}_{\Lambda,\rho,N}$ and $(b_\epsilon)_\epsilon\in\mathcal{S}^{m^{'}}_{\Lambda,\rho, N^{'}}$ with $m^{'}<l$, the symbol $a_\epsilon(z)+\epsilon^{M+N^{'}}b_\epsilon(z)$ belongs to $H\mathcal{S}^{m,l}_{\Lambda,\rho,M}$ for $M\ge N$.
\end{remark}
\ni Classically we find interesting examples of hypoelliptic symbols considering polynomials of the form $\sum_{\alpha\in\mathcal{A}}c_\alpha z^\alpha$ where $\mathcal{A}$ is a finite subset of $\mathbb{N}^{2n}$ and $c_\alpha\in\mathbb{C}$. In the sequel we collect some classical results and examples, referring for details to \cite{bog, bogII}.
\begin{example}\bf{Standard elliptic polynomials}\rm\\
Fix an integer $m\ge 1$. Defining the weight function $\Lambda(z)=\langle z\rangle$, we have that  $a(z)=\sum_{|\alpha|\le m}c_\alpha z^\alpha$ belongs to $E\mathcal{S}^{m}_{\Lambda,1,0}$ iff $\sum_{|\alpha|=m}c_\alpha z^\alpha\neq 0$ for $z\neq 0$. For $z=(x,\xi)\in\mathbb{R}^2$, simple examples of standard elliptic polynomials are given by
\[
x^m+i\xi^m ,\quad\quad\quad\quad  m\ \text{positive integer}
\]
and 
\[
\quad\quad x^m+\xi^m ,\quad\quad\quad\quad m\ \text{even positive integer}.
\]
\end{example}
\begin{example}\bf{Quasi-elliptic polynomials}\rm\\
Fix a $2n$-tuple $M=(M_1,...,M_{2n})$ of positive integers.\\
Let us write $\mu=\max_{j}M_j$ and $m=(m_1,...,m_{2n})$ with $m_j=\mu /M_j$ for $j=1,...,2n$. Choosing the weight function $\Lambda(z)=(1+\sum_{j=1}^{2n}z_j^{2M_j})^{\frac{1}{2\mu}}$, $a(z)=\sum_{\alpha\cdot m\le\mu}c_\alpha z^\alpha$ belongs to $E\mathcal{S}^{\mu}_{\Lambda,1,0}$ iff $\sum_{\alpha\cdot m=\mu}c_\alpha z^\alpha\neq 0$ for $z\neq 0$.\\
A simple example of a quasi-elliptic polynomial is given for $z=(x,\xi)\in\mathbb{R}^2$ by
\[
\xi^h+rx^k ,
\]
where $r\in\mathbb{C}$ with $\Im r\neq 0$, $h,k$ are positive integers, $M=(k,h)$ and $\Lambda(x,\xi)=(1+x^{2k}+\xi^{2h})^{\frac{1}{2\max{(k,h)}}}$.
\end{example}
\ni At this point we may also consider $a_\epsilon(z)=\epsilon^N a(z)$ with $a(z)=\sum_{\alpha\in\mathcal{A}}c_\alpha z^\alpha$
a standard elliptic or quasi-elliptic symbol. For a suitable weight function $\Lambda$ and a suitable order $k$ we obtain an element of $E\mathcal{S}^k_{\Lambda,1,N}$.\\
Finally we look for hypoelliptic symbols in the more general family of polynomials $\sum_{\alpha\in\mathcal{A}}c_{\alpha,\epsilon}z^\alpha$ having coefficients $(c_{\alpha,\epsilon})_\epsilon\in\mathcal{E}_{o,M}$. We start with the following generic situation.
\begin{example}
Let $\sum_{\alpha\in\mathcal{A}}c_\alpha z^\alpha$ be a hypoelliptic symbol in $H\mathcal{S}^{m,l}_{\Lambda,\rho, 0}$ with coefficients $c_\alpha\in\mathbb{C}$. Let $\sum_{\alpha\in\mathcal{A}^{'}}c^{'}_{\alpha,\epsilon} z^\alpha$ be a polynomial with coefficients $(c^{'}_{\alpha,\epsilon})_\epsilon\in\mathcal{E}_{o,M}$, which belongs to $\mathcal{S}^{m^{'}}_{\Lambda,\rho,N^{'}}$ for certain $m^{'}\in\mathbb{R}$ and $N^{'}\in\mathbb{N}$. If $m^{'}<l$ we know from Proposition 8.4 that 
\begin{equation}
\epsilon^{N}\sum_{\alpha\in\mathcal{A}}c_\alpha z^\alpha+\epsilon^{N+N^{'}}\sum_{\alpha\in\mathcal{A}^{'}}c^{'}_{\alpha,\epsilon} z^\alpha
\end{equation}
is an element of $H\mathcal{S}^{m,l}_{\Lambda,\rho,N}$.
\end{example}
\ni In order to present concrete examples we need a technical lemma.
\begin{lemma}
{\ }
\begin{itemize}
\item[i)] Let $\sum_{|\alpha|\le m^{'}}c^{'}_{\alpha,\epsilon}z^\alpha$ be a polynomial with coefficients $(c^{'}_{\alpha,\epsilon})_\epsilon\in\mathcal{E}_{o,M}$, i.e.
\begin{equation}
\forall\alpha\in\mathbb{N}^n,\ |\alpha|\le m^{'},\ \exists N^{'}_\alpha\in\mathbb{N},\ \exists c^{'}_\alpha>0:\ \forall\epsilon\in(0,1],\qquad |c^{'}_{\alpha,\epsilon}|\le c^{'}_\alpha\epsilon^{-N^{'}_\alpha}.\qquad\quad
\end{equation}
Then $\big(\sum_{|\alpha|\le m^{'}}c^{'}_{\alpha,\epsilon}z^\alpha\big)_\epsilon\in\mathcal{S}^{m^{'}}_{\Lambda,1,N^{'}}$, where $\Lambda(z)=\langle z\rangle$ and $N^{'}\hskip-2pt=\max_{|\alpha|\le m^{'}}N^{'}_\alpha$.
\item[ii)] Let $M\in\mathbb{N}^{2n}$, $m\in({\mathbb{R}}^+)^{2n}$, $\mu\in\mathbb{N}$ be given as in Example 2. Let $\sum_{\alpha\cdot m\le\mu^{'}}c^{'}_{\alpha,\epsilon}z^\alpha$ be a polynomial with coefficients $(c^{'}_{\alpha,\epsilon})_\epsilon\in\mathcal{E}_{o,M}$ as in (8.36) .
Then $(\sum_{\alpha\cdot m\le\mu^{'}}c^{'}_{\alpha,\epsilon}z^\alpha)_\epsilon\in\mathcal{S}^{\mu^{'}}_{\Lambda,1,N^{'}}$, where $\Lambda(z)=(1+\sum_{j=1}^{2n}z_j^{2M_j})^{\frac{1}{2\mu}}$ and $N^{'}\hskip-4pt=\hskip-2pt\max_{\alpha\cdot m\le\mu^{'}}N^{'}_\alpha$.
\end{itemize}\quad
\end{lemma}
\begin{proof}
In the first case it is sufficient to apply Proposition 4.4. We can easily prove the second statement recalling that for $\Lambda(z)=(1+\sum_{j=1}^{2n}z_j^{2M_j})^{\frac{1}{2\mu}}$ and for $\alpha\in\mathbb{N}^{2n}$, $|z^\alpha|\prec\Lambda(z)^{\alpha\cdot m}$ (see \cite{bog}, Lemma 8.1, Example 8.2).
\end{proof}
\ni Combining the previous lemma with Proposition 8.4 and Example 3 we arrive at the following final example. 
\begin{example}
Let $a$ be a standard elliptic polynomial in $E\mathcal{S}^{m}_{\Lambda,1,0}$ with $\Lambda(z)=\langle z\rangle$ and let\\ $\sum_{|\alpha|\le m^{'}}c^{'}_{\alpha,\epsilon}z^\alpha$ be a polynomial with coefficients in $\mathcal{E}_{o,M}$ and $m^{'}<m$. There exists a natural number $N^{'}$ such that for all $N\in\mathbb{N}$, the symbol
\[
\epsilon^{N}a(z)+\epsilon^{N+N^{'}}\sum_{|\alpha|\le m^{'}}c^{'}_{\alpha,\epsilon}z^\alpha
\]
belongs to $E\mathcal{S}^{m}_{\Lambda,1,N}$. Analogously if $a$ is a quasi-elliptic polynomial in $E\mathcal{S}^{\mu}_{\Lambda,1,0}$ with $\Lambda(z)=(1+\sum_{j=1}^{2n}z_j^{2M_j})^{\frac{1}{2\mu}}$ and  $\sum_{\alpha\cdot m\le\mu^{'}}c^{'}_{\alpha,\epsilon}z^{\alpha}$ is a polynomial with coefficients in $\mathcal{E}_{o,M}$ and $\mu^{'}<\mu$, then there exists $N^{'}\in\mathbb{N}$ such that for all $N\in\mathbb{N}$ 
\[
\epsilon^{N}a(z)+\epsilon^{N+N^{'}}\sum_{\alpha\cdot m\le\mu^{'}}c^{'}_{\alpha,\epsilon}z^\alpha
\]
belongs to $E\mathcal{S}^{\mu}_{\Lambda,1,N}$. 
\end{example}
\ni In conclusion we consider a partial differential operator $\sum_{(\alpha,\beta)\in\mathcal{A}}c_{\alpha,\beta}x^\alpha{D}^\beta$ with coefficients $c_{\alpha,\beta}\in\overline{\mathbb{C}}$ and ${D}^\beta =(-i)^{|\beta|}\partial^\beta$. It is a linear map from $\gts$ into $\gts$. As proved in Proposition 4.4 for any weight function $\Lambda$, there exists $r\in\mathbb{R}$ such that the polynomial $\sum_{(\alpha,\beta)\in\mathcal{A}}c_{\alpha,\beta}x^\alpha\xi^\beta\in\overline{\mathbb{C}}[x,\xi]$ can be considered as an element of the factor $\mathcal{S}^r_{\Lambda,1,N}/\mathcal{N}^r_{\Lambda,1}$, for a suitable $N$ depending on the coefficients $c_{\alpha,\beta}$. In the sequel we shall investigate the re\-gu\-la\-ri\-ty properties of such a partial differential operator using the tools provided by the pseudo-differential calculus.  
\begin{lemma}
Let $\sum_{(\alpha,\beta)\in\mathcal{A}}c_{\alpha,\beta}x^\alpha{D}^\beta$ be a partial differential operators with coefficients in $\overline{\mathbb{C}}$, let $(a_\epsilon(x,\xi))_\epsilon:=\big(\sum_{(\alpha,\beta)\in\mathcal{A}}c_{\alpha,\beta,\epsilon}x^\alpha{\xi}^\beta)_\epsilon\in\mathcal{S}^{m^{'}}_{\Lambda,\rho,N^{'}}$ be a symbol obtained from the polynomial $\sum_{(\alpha,\beta)\in\mathcal{A}}c_{\alpha,\beta}x^\alpha{\xi}^\beta$ and $A$ the corresponding pseudo-differential operator. For all pseudo-dif\-fe\-ren\-tial operators $P$ with regular symbol $(p_\epsilon)_\epsilon\in\mathcal{S}^{m^{''}}_{\Lambda,\rho,N^{''}}$ and for all $u\in\gts$ the following weak equality holds:
\begin{equation}
P\big(\sum_{(\alpha,\beta)\in\mathcal{A}}c_{\alpha,\beta}x^\alpha{D}^\beta u\big)=_{g.t.d.}PAu .
\end{equation}
\end{lemma}
\begin{proof}
At first note that $(A_\epsilon u_\epsilon)_\epsilon =(\sum_{(\alpha,\beta)\in\mathcal{A}}c_{\alpha,\beta,\epsilon}x^\alpha{D}^\beta(u_\epsilon\widehat{\varphi_\epsilon}))_\epsilon$ is a representative of $Au$. We have to show that for all $f\in\S$, 
\[
\sum_{(\alpha,\beta)\in\mathcal{A}}\int_{\mathbb{R}^n}f(x)\int_{\mathbb{R}^{2n}}e^{i(x-y)\xi}p_\epsilon(x,\xi)c_{\alpha,\beta,\epsilon}y^\alpha[D^\beta u_\epsilon-D^\beta(u_\epsilon\widehat{\varphi_\epsilon})](y)\widehat{\varphi_\epsilon}(y)\, dy\dslash\xi\, dx
\]
defines an element of $\mathcal{N}_o$. Changing order in integration we get for each $(\alpha,\beta)\in\mathcal{A}$
\begin{equation}
\label{coeff}
c_{\alpha,\beta,\epsilon}\int_{\mathbb{R}^n}\int_{\mathbb{R}^{2n}}e^{i(x-y)\xi}p_\epsilon(x,\xi)f(x)\, dx\,\dslash\xi\, y^\alpha D^\beta(u_\epsilon-u_\epsilon\widehat{\varphi_\epsilon})(y)\widehat{\varphi_\epsilon}(y) dy ,
\end{equation}
where $(g_\epsilon(y))_\epsilon:=(\int_{\mathbb{R}^{2n}}e^{i(x-y)\xi}p_\epsilon(x,\xi)f(x)\, dx\dslash\xi)_\epsilon$ belongs to $\mathcal{E}^\infty_{\mathcal{S}}(\mathbb{R}^n)$ since $(p_\epsilon)_\epsilon$ is a regular symbol. Using integration by parts we rewrite \eqref{coeff} as 
\begin{equation}
\label{coeff1}
(-1)^{|\beta|}c_{\alpha,\beta,\epsilon}\int_{\mathbb{R}^n}D^\beta(g_\epsilon(y)y^\alpha\widehat{\varphi_\epsilon}(y))u_\epsilon(y)(\widehat{\varphi_\epsilon}(y)-1)\, dy .
\end{equation}
Taylor's formula applied to $\widehat{\varphi}$ at $0$, combined with the properties of $(g_\epsilon)_\epsilon\in\mathcal{E}^\infty_{\mathcal{S}}(\mathbb{R}^n)$, gives us the estimates characterizing $\mathcal{N}_o$ in \eqref{coeff1}.
\end{proof}
\begin{proposition}
Let $\sum_{(\alpha,\beta)\in{\mathcal{A}}}c_{\alpha,\beta}x^\alpha D^\beta$ be a partial differential operator with coefficients in $\overline{\mathbb{C}}$. If there exists a representative $(a_\epsilon(x,\xi))_\epsilon:=(\sum_{(\alpha,\beta)\in\mathcal{A}}c_{\alpha,\beta,\epsilon}x^\alpha\xi^\beta)_\epsilon$ of $\sum_{(\alpha,\beta)\in\mathcal{A}}c_{\alpha,\beta}x^\alpha\xi^\beta$ belonging to the set $H\mathcal{S}^{m,l}_{\Lambda,\rho,N}$ of hypoelliptic symbols, then 
\begin{equation}
\sum_{(\alpha,\beta)\in{\mathcal{A}}}\hskip-5pt c_{\alpha,\beta}x^\alpha D^\beta u=v ,
\end{equation}
where $u\in\gts$ and $v\in\gss$, implies that $u$ is equal in the weak sense to a generalized function in $\gss$.
\end{proposition}
\begin{proof}
Let $A$ be the pseudo-differential operator with symbol $(a_\epsilon)_\epsilon$ and let $P$ be a parametrix of $A$. From Lemma 8.2, $Pv=P(\sum_{(\alpha,\beta)\in{\mathcal{A}}}c_{\alpha,\beta}x^\alpha D^\beta u)=_{g.t.d.}PAu$, where $Pv\in\gss$. We complete the proof recalling that from Theorem 8.1 there exists an operator $R_1$ with $\mathcal{S}$-regular kernel such that $Pv=_{g.t.d.}u+R_1u$.
\end{proof}
\ni \bf{Acknowledgements.}\rm\ The author is thankful to Dr. $\rm{G\ddot{u}nther}$ $\rm{H\ddot{o}rmann}$, Prof. Michael Oberguggenberger and  Prof. Luigi Rodino for helpful discussions during the preparation of the paper.

\end{document}